\documentclass[preprint,12pt,authoryear]{elsarticle}

\usepackage{amssymb}
\usepackage{subfigure}
\usepackage{amsfonts}
\usepackage[ruled,vlined,linesnumbered]{algorithm2e} 
\usepackage{multirow}
\usepackage{xcolor}
\usepackage{graphicx}
\usepackage{epstopdf}
\usepackage{hyperref}

\journal{Computers \& Chemical Engineering}

\begin{document}

\begin{frontmatter}



\title{A Multilevel Coordinate Search Algorithm for Well Placement, Control and Joint Optimization}


\author[add_upc]{Xiang Wang\corref{co_author_main}}
\cortext[co_author_main]{Principal corresponding author}
\ead{xiangwangdr@gmail.com}

\author[add_mun]{Ronald D. Haynes\corref{co_author}}
\cortext[co_author]{Corresponding author}
\ead{rhaynes@mun.ca}

\author[add_upc]{Qihong Feng}
\ead{fengqihong@126.com}

\address[add_upc]{School of Petroleum Engineering, China University of Petroleum (East China), Qingdao, Shandong, China 266580}
\address[add_mun]{Department of Mathematics \& Statistics, Memorial University of Newfoundland, St. John's, NL, Canada A1C 5S7}

\begin{abstract}
Determining optimal well placements and controls are two important tasks in oil field development. These problems are computationally expensive, nonconvex, and contain multiple optima. The practical solution of these problems require efficient and robust algorithms. In this paper, the multilevel coordinate search (MCS) algorithm is applied for well placement and control optimization problems. MCS is a derivative-free algorithm that combines global and local search. Both synthetic and real oil fields are considered. The performance of MCS is compared to generalized pattern search (GPS), particle swarm optimization (PSO), and covariance matrix adaptive evolution strategy (CMA-ES) algorithms. Results show that the MCS algorithm is strongly competitive, and outperforms for the joint optimization problem and with a limited computational budget. The effect of parameter settings for MCS are compared for the test examples. For the joint optimization problem we compare the performance of the simultaneous and sequential procedures and show the utility of the latter.
\end{abstract}


\end{frontmatter}



\section{Introduction}
\label{sec:1_intro}

Determining the optimal well locations and controls in an oil field is a challenging task. The decision is hard since the reservoir performance is affected by geological, engineering, economical and other parameters \citep{tavallali_optimal_2013,knudsen_shut-based_2013,shakhsi-niaei_optimal_2014}. Optimization algorithms provide a systematic way to solve this problem. By using optimization algorithms, a quality solution can be achieved automatically and hence reduce the risk in decision-making. Well placement and control optimization generally are computationally expensive and nonconvex, and not every optimization algorithm is appropriate for these problems. Therefore, finding and applying algorithms that are efficient and robust is one of most important tasks in solving well placement and control optimization problems.

In this work, we introduce and apply the multilevel coordinate search (MCS) algorithm for the problems of optimizing well placement, well control, and joint placement with control. MCS,  introduced by \citep{huyer_global_1999}, is a global optimization algorithm and is designed to handle the complex topography and multimodality of the multidimensional nonlinear objective functions without requiring excessive computing resources. Rios \citep{rios_derivative-free_2013} completed a systematic comparison using a test set of 502 problems and found that MCS outperforms the other 21 derivative-free algorithms tested (see Table \ref{tab:solver}). Though MCS has shown its superiority in benchmark and real world problems \citep{huyer_global_1999,rios_derivative-free_2013,lambot_global_2002}, to the best of our knowledge, it has not been applied to the optimization of oil field development. We compare MCS, generalized pattern search (GPS), particle swarm optimization (PSO), and covariance matrix adaptive evolution strategy (CMA-ES) in four typical test cases from the field of optimal reservoir production development. Our results demonstrate that MCS is strongly competitive and outperforms the other algorithms in most cases.

\begin{table}[htb]
\caption{Derivative-free solvers considered by \citep{rios_derivative-free_2013}.}
\label{tab:solver}
\centering
\begin{tabular}{lll}
\hline\noalign{\smallskip}
Solver & Version & Language\\
\noalign{\smallskip}\hline\noalign{\smallskip}
ASA  &  26.30  &  C   \\
BOBYQA  &  2009  &  Fortran   \\
CMA-ES  &  3.26beta  &  Matlab  \\
DAKOTA/DIRECT  &  4.2  &  C++   \\
DAKOTA/EA  &  4.2  &  C++  \\
DAKOTA/PATTERN  &  4.2  &  C++  \\
DAKOTA/SOLIS-WETS  &  4.2  &  C++  \\
DFO  &  2.0  &  Fortran   \\
FMINSEARCH  &  1.1.6.2  &  Matlab  \\
GLOBAL  &  1.0  &  Matlab \\
HOPSPACK  &  2.0  &  C++  \\
IMFIL  &  1.01  &  Matlab  \\
MCS  &  2.0  &  Matlab  \\
NEWUOA  &  2004  &  Fortran  \\
NOMAD  &  3.3  &  C++  \\
PSWARM  &  1.3  &  Matlab  \\
SID-PSM  &  1.1  &  Matlab  \\
SNOBFIT  &  2.1  &  Matlab  \\
TOMLAB/GLCCLUSTER  &  7.3  &  Matlab  \\
TOMLAB/LGO  &  7.3  &  Matlab  \\
TOMLAB/MULTIMIN  &  7.3  &  Matlab  \\
TOMLAB/OQNLP  &  7.3  &  Matlab  \\
\noalign{\smallskip}\hline
\end{tabular}
\end{table}

Oil field development optimization has two main sub-problems: well placement optimization, and well control optimization. These two problems are often treated separately \citep{oliveira_adaptive_2014,bouzarkouna_well_2012,wang_production_2009,brouwer_dynamic_2004}. Recently, there has been increasing focus on optimizing well placement and control jointly \citep{forouzanfar_covariance_2015,humphries_simultaneous_2013,isebor_generalized_2014}. Well placement problems aim to optimize the locations of injection and production wells. The location of each vertical well is parametrized by its  plane coordinates $(x,y)$, which are usually integers in the reservoir simulator. Well control problems focus on optimizing production scheduling. The optimization variables are often the time-varying bottom hole pressures (BHPs) or the flow rates for each well. The joint problem optimizes well placement and control parameters simultaneously. Thus, the joint problems are more complex and challenging with an increase in the number and type of variables \citep{isebor_generalized_2014}. 

In the past, a number of algorithms have been devised and analysed for both separate and joint problem of well placement and control optimization. These algorithms fall into two categories: gradient-based methods and derivative-free methods. Applications of gradient-based methods to oil field problems have been presented in many papers \citep{volkov_effect_2014,wang_production_2009,brouwer_dynamic_2004,zandvliet_adjoint-based_2008,sarma_efficient_2006,zhou_optimal_2013}. These methods take advantage of the gradient information to guide their search. 
The gradient of the objective function can be obtained by using adjoint-based techniques \citep{brouwer_dynamic_2004,sarma_efficient_2006,zandvliet_adjoint-based_2008,volkov_effect_2014}, or may be approximated by using numerical methods such as finite differences \citep{wang_production_2009,zhou_optimal_2013}. 
The adjoint method, developed in the 1970s \citep{chen_new_1974,chavent_identification_1974}, is widely used for assisted history matching \citep{wu_conditioning_1999,li_history_2003} and well production optimization \citep{asheim_maximization_1988,zakirov_optimizing_1996,brouwer_dynamic_2004}.
Gradient based methods have some drawbacks for the well placement and control problem; these problems are nonconvex and generally contain multiple optima. For some problems, particularly well placement, the optimization surface can be very rough, which results in discontinuous gradients \citep{ciaurri_derivative-free_2011}.
However, the gradient-based methods are often the most efficient methods especially for the optimal well control problem \citep{zhao_maximization_2013,handels_adjoint-based_2007,vlemmix_adjoint-based_2009,wang_optimal_2007,forouzanfar_joint_2014}.

For the joint well placement and control optimization problem, two procedures are proposed and studied. The first one is a simultaneous procedure, which optimizes over all well locations and control parameters simultaneously. The second one is a sequential procedure, that decouples the joint problem into the well placement optimization subproblem and the well control placement optimization subproblem. The simultaneous procedure ensures that the best solution exists somewhere in the search space. But it may be difficult to find the global optima because the search space may be very large and rough. The sequential procedure divides the optimization variables into two smaller groups and optimizes separately. For each subproblem, the search space is smaller than the simultaneous one, but it can not ensure the best solution exists in the search space because the optimal location depends on how the well is operated and vice-versa. 
There are several papers \citep{li_simultaneous_2012,bellout_joint_2012,isebor_derivative-free_2014} which demonstrate that the simultaneous procedure is superior to the sequential approach. In \citep{humphries_simultaneous_2013,humphries_joint_2015}, however, they found that a sequential procedure was competitive and even preferable to the simultaneous approach in some cases. To test this further, we do a further investigation of the effectiveness of these two procedures using a joint placement and control optimization example.

Many black-box, derivative-free methods have been used in oil field problems \citep{merlini_giuliani_derivative-free_2015}. Typical algorithms include genetic algorithms (GA) \citep{alqahtani_computational_2014,bouzarkouna_well_2012}, particle swarm optimization (PSO) \citep{onwunalu_development_2009,onwunalu_new_2011}, generalized pattern search (GPS) \citep{asadollahi_production_2014,isebor_constrained_2009}, covariance matrix adaptation strategy (CMA-ES) \citep{bouzarkouna_well_2012,forouzanfar_covariance_2015} and hybrid approaches \citep{isebor_generalized_2014,humphries_joint_2015}. These algorithms can be classified as either deterministic or stochastic, and provide global or local search. The stochastic/global approaches have also been hybridized with deterministic/local search techniques. These algorithms only require the value of objective function and involve no explicit gradient calculations. 
For smooth objective functions, the
derivative-free methods generally need more function evaluations to converge to a solution than the gradient-based ones. 
However, most of the derivative-free algorithms parallelize naturally and easily, which make their efficiency satisfactory \citep{ciaurri_derivative-free_2011}, and indeed these methods are less likely to become trapped in local optima. 

We are particularly interested in using the multilevel coordinate search (MCS) algorithm for the following reasons: 1) it combines a global search with a local search, which leads to a quicker convergence than many methods that operate only at the global level. 2) it is an intermediate between heuristic methods that find the global optimum only with high probability and methods that guarantee to find a optimum with a required accuracy. 3) it does not need analytic or numerical derivatives. 4) it is guaranteed to converge if the objective is continuous in the neighbourhood of a global minimizer, no additional smoothness properties are required. 5) the algorithm parameters in MCS have a clear meaning and are easy to choose. 6) it has proved itself in benchmark tests and many real world problems  \citep{huyer_global_1999}. Based on these features, we believe that MCS has great potential to solve oil field optimization problems, which are nonconvex, nonlinear, and contain many local optima and discontinuities.

In this paper, we apply MCS to optimization problems of varying complexity in terms of the number and type of optimization variables, the dimension and size of the reservoir models, and the number of wells. We investigate the effect of the algorithmic parameters (initialization list type, number of levels, and the effect of local search) on the optimization results. We propose a detailed comparison between MCS and three other popular derivative-free algorithms (GPS, PSO, and CMA-ES).

This paper is organized as follows: Section \ref{sec:2_problem} describes the formulation of the optimization problems. Section \ref{sec:3_solver} gives an overview of the optimization algorithms considered.
In Section \ref{sec:4_exam_resu} we describe our numerical experiments and the corresponding results. Finally, in Section \ref{sec:6_concl} we provide some concluding remarks.

\section{Problem Formulation}
\label{sec:2_problem}


\subsection{General Problem Statement}
\label{sec:21_genproblem}

The objective functions for general oil field development optimization problems are often the net present value (NPV) or cumulative oil production \citep{sarma_impl_2005,wang_production_2009,isebor_generalized_2014,oliveira_adaptive_2014}. We use NPV as the objective function for all our work. NPV accounts for revenue from the oil and gas produced, and for the cost of handling water production and injection. The NPV is defined as
\begin{equation}\label{eq:npv}
\mathrm{NPV}=\sum\limits_{k=1}^{N_t} \left[ {\frac{\Delta {t_k}}{{\left( 1+b \right)}^{\frac{t_k}{\tau }}}} \left( \sum\limits_{i=1}^{N_p}{r_{gp}q_{gp}^{i,k}}+\sum\limits_{i=1}^{N_p}{r_{op}q_{op}^{i,k}}-\sum\limits_{i=1}^{N_p}{c_{wp}q_{wp}^{i,k}}-\sum\limits_{i=1}^{N_i}{c_{wi}q_{wi}^{i,k}}    \right) \right],
\end{equation}
where $q_{gp}^{i,k}$, $q_{op}^{i,k}$, $q_{wp}^{i,k}$ and $q_{wi}^{i,k}$ are the flow rates of the gas, oil, water produced and water injected for well $i$ at time step $k$, respectively; $r_{gp}$ and $r_{op}$ are the gas and oil revenue; $c_{wp}$ and $c_{wi}$ are the costs of produced water and the costs of injected water. $N_t$ is total number of time steps, $t_k$ is the time at the end of $k$th time step, $\Delta t_k$ is $k$th time step size, $\tau$ provides the appropriate normalization for $t_k$, e.g., $\tau = 365$ days, and $b$ is the fractional discount rate.

The quantities $q_{gp}^{i,k}$, $q_{op}^{i,k}$, $q_{wp}^{i,k}$ are functions of the dynamic state variables $\bf x$ (e.g., pressure, saturation), the geological and fluid variables $\bf m$ (e.g., permeability, porosity, viscosity and density of each phase, etc.), the well placement vector $\bf v$, and the control vector $\bf u$. Hence equation (\ref{eq:npv}) can be written as
\begin{equation}\label{eq:npv1}
\mathrm{NPV}=J(\bf x,\bf m,\bf v,\bf u),
\end{equation}
where $J$ is a nonlinear response to the input variables -- $J$ depends on the output from the numerical solution of a system of nonlinear PDEs describing the reservoir response.  

For the well placement and control optimization problem, we wish to maximize the net present value $J$ by adjusting the placement vector $\bf v$ and the well control vector $\bf u$ subject to the constraints
\begin{equation}\label{eq:constr1}
\bf g(\bf x,\bf m, \bf v, \bf u)=0,
\end{equation}
\begin{equation}\label{eq:constr2}
\bf c(\bf x,\bf m, \bf v, \bf u)\leq 0,
\end{equation}
\begin{equation}\label{eq:constr3}
\bf{A_{eq,v} v}= \bf{b_{eq,v}},\quad \bf{A_{eq,u} u}= \bf{b_{eq,u}},
\end{equation}
\begin{equation}\label{eq:constr4}
\bf{A_{ieq,v} v}\leq \bf{b_{ieq,v}},\quad \bf{A_{ieq,u} u}\leq \bf{b_{ieq,u}},
\end{equation}
\begin{equation}\label{eq:constr5}
\bf{v}_{lb}\leq \bf v \leq \bf{v}_{ub},\quad \bf{u}_{ub}\leq \bf u \leq \bf{u}_{ub}.
\end{equation}

Equation (\ref{eq:constr1}) represents the reservoir simulation equations. This constraint ensures that the governing reservoir flow equations are satisfied. Equation (\ref{eq:constr2}) describes the  nonlinear constraint functions (e.g., distance between wells, well water cuts, facility constraints such as field-level injection limits, etc.). Equations 
(\ref{eq:constr3}--\ref{eq:constr5}) are the equality, inequality, and bound constraints of well placement vector and well control vector, respectively. These constraints are related to, for example, the capacity limitations of the wells and/or facilities. 

Due to the complexity of the well placement and control optimization problem, the objective function is not convex, it may have many minima, maxima, and saddlepoints. Furthermore, the objective surface is very rough, it is therefore hard to extract gradient information. 

\begin{figure}[htb]
  \centering 
  \subfigure[NPV surface map]{ 
    \label{fig:subfig:npvsurf} 
    \includegraphics[width=0.55\textwidth]{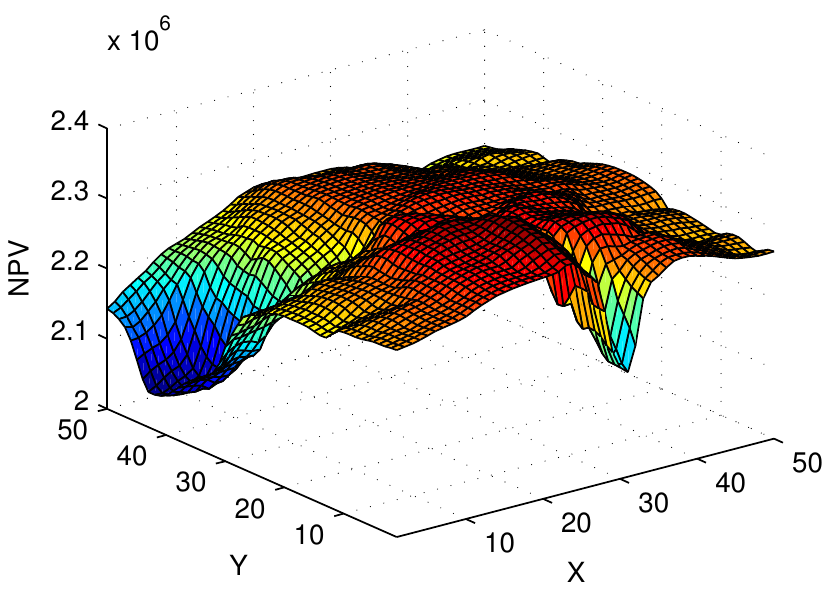}} 
  \subfigure[NPV contour map]{ 
    \label{fig:subfig:npvcont} 
    \includegraphics[width=0.35\textwidth]{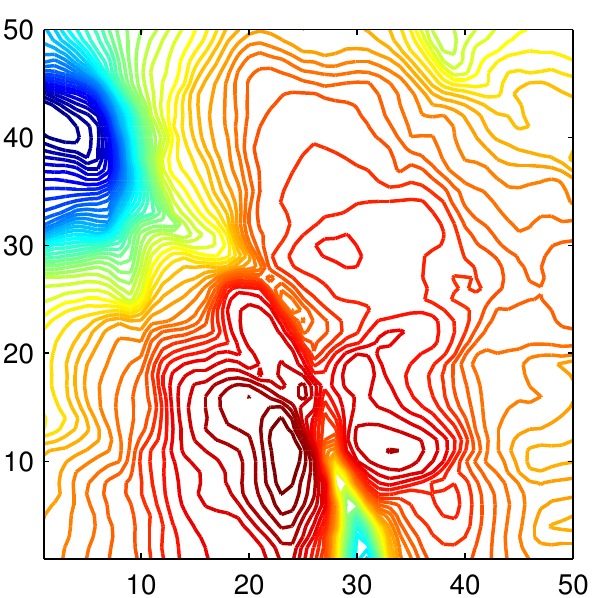}}  
  \caption{NPV surface map and contour map.} 
  \label{fig:npvmap} 
\end{figure}

To illustrate this, a two-dimensional heterogeneous reservoir model is selected. The model uses 50$\times$50 grid blocks with four producers, one in each of the four corners. We optimize the location of one injection well. The NPV surface map is generated numerically by putting a single injection well at each grid block and computing the NPV of production from the reservoir. The NPV surface map and contour map are given in Fig.\ \ref{fig:npvmap}. From the figure we can see that the NPV surface is noisy,  and includes at least five local optima. This optimization problem has only two variables; it can be regarded as the simplest problem in well placement and control optimization. For more complex problems, the objective surface has far more roughness with more local optima.  \citep{bangerth_on_2006} and \citep{forouzanfar_optimization_2012} also discuss the roughness of the NPV surface.

\subsection{Well placement optimization} 
\label{sec:22placementop}

In the well placement optimization problem, we seek to determine the optimal locations for a specified number of vertical wells, or the optimal trajectories for a specified number of 3-D angled wells. The optimization problem studied here can be expressed as

\begin{eqnarray}
\max\limits& \ \mathrm{NPV}=J(\bf v) \\
\mathrm{subject \ to} & \ \mathbf{v_{lb}} \leq \mathbf{v} \leq \mathbf{v_{ub}} .
\end{eqnarray}

For vertical wells, the location of each well is given by its plane coordinates $(x,y)$. Thus the total number of variables for well placement vector $\bf v$ is $2N$, where $N=N_p+N_i$ is the total number of wells, and $N_p$ and $N_i$ are the number of production and injection wells, respectively.

The placement of each of 3-D angled well is parameterized by 6 variables, $x$, $y$, $z$, $l$, $\theta$, and $\phi$ \citep{yeten_optimization_2002,bouzarkouna_well_2012,humphries_joint_2015}. 

\begin{figure}[htb]
  \centering 
  \subfigure[Vertical well]{ 
    \label{fig:subfig:vwell} 
    \includegraphics[width=0.45\textwidth]{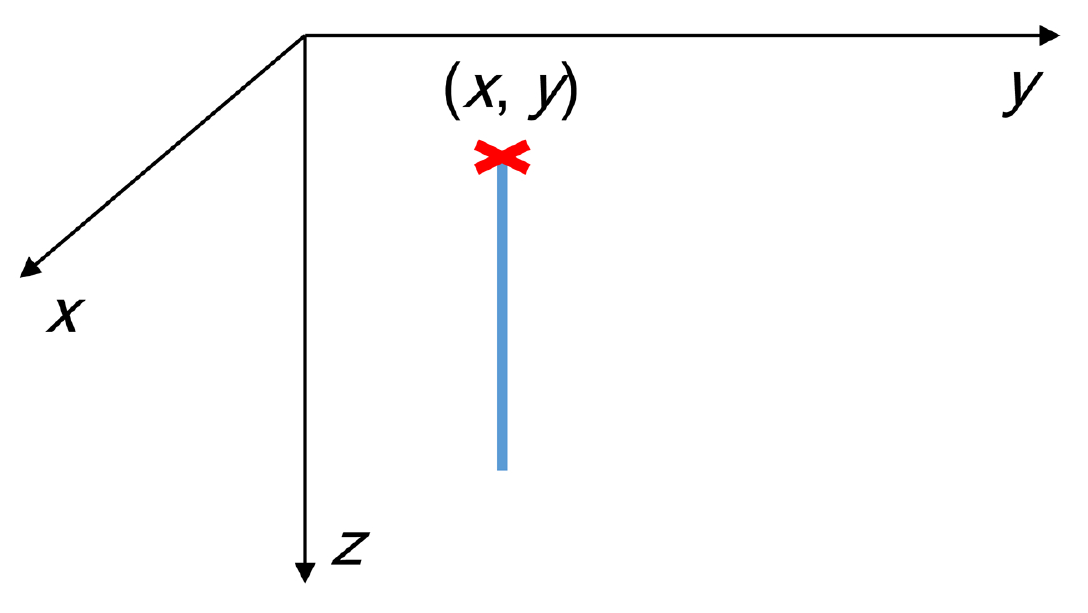}} 
  \subfigure[3-D angled well]{ 
    \label{fig:subfig:awell} 
    \includegraphics[width=0.45\textwidth]{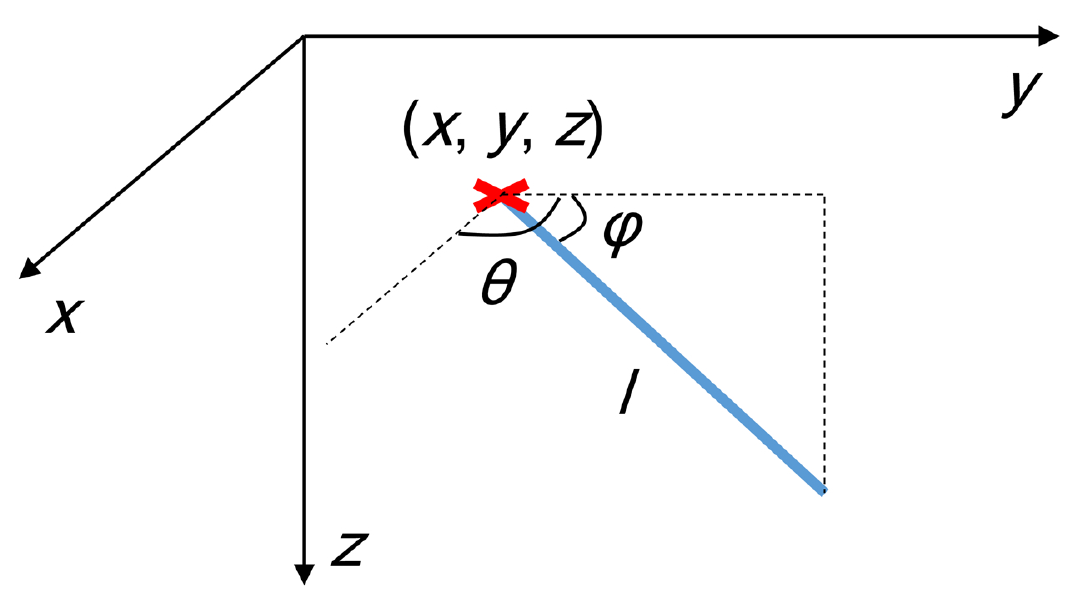}} 
  \caption{Types of well and parameterizations.} 
  \label{fig:wellpara} 
\end{figure}

The coordinates of the well heel are given by $x$ and $y$, and the depth of the well heel is given by $z$. $l$ is the length of the well. $\theta$ is the angle of the well in the $x$-$y$ plane, and $\phi$ is the angle of the well makes with the horizontal plane. For this 3-D angled well placement optimization problem, the total number of variables is $6N$. The parameterizations for vertical wells and 3-D angled wells are illustrated in Fig.\ \ref{fig:wellpara}. 

As in Section \ref{sec:21_genproblem}, the general oilfield development optimization problem includes a set of constraints. 
For our well placement and control optimization problem, we assume only bound constraints are imposed explicitly.
The governing reservoir flow equations, equation (\ref{eq:constr1}), are always satisfied since we use a reservoir simulator to calculate the objective function value. 
Some complex constraints, for example, the minimum distance between wells may be naturally enforced since, if two wells are very close, the NPV will generally not be high.
Alternatively, one can impose the minimum distance as explicit constraint as in \citep{li_variable-control_2012,emerick_well_2009}. The effect of imposing such constraints explicitly or implicitly on the optimization needs further research. In this paper we choose not to explore this topic. 
We also require that wells are not outside the reservoir boundary. For reservoirs with irregular boundaries, this is a nonlinear constraint. If a vertical well is placed outside the reservoir, or the 3-D angled well is drilled at inactive gridblocks, it can not produce oil. This leads to a low NPV and will be eliminated by the optimization algorithm. 

Some complex constraints can be transformed to bound constraints. The placement of a 3-D angled well can be parameterized either by $(x, y, z, l, \theta, \phi)$ or by $(x_1, y_1, z_1, x_2, y_2, z_2)$, where $(x_1, y_1, z_1)$ and $(x_2, y_2, z_2)$ are the heel point and the toe point of well, respectively \citep{emerick_well_2009}. 
Using the latter approach, the constraints of the maximum and minimum well length can be written as
\begin{equation}
a\leq (x_2-x_1)^2+(y_2-y_1)^2+(z_2-z_1)^2 \leq b,
\end{equation}
which are two nonlinear inequality constraints.  
For the first approach, the maximum and minimum well length can be stated as bound constraints.
Choosing how to enforce constraints requires some experience, and are one of the main difficulties in formulating and solving optimization problems.

Most of the optimization algorithms considered in this paper were originally designed for either the unconstrained optimization problem or problems with bound constraints only. While constraints can decrease the size of the search space, treating constraints explicitly does often burden the computation. 

The coordinates of each well are real variables in actual oilfields, but are usually treated as integers (i.e. gridblock locations ) in reservoir simulators. The well coordinate variables in our work are continuous but will be rounded before we pass them to the simulator to evaluate the objective function. This leads a discontinuity in NPV as the well moves across the boundary between two gridblocks. Optimization becomes problematic when there are discontinuities in objective function surface while using gradient-based algorithms. The derivative-free algorithms considered in this paper, MCS, GPS, PSO, and CMA-ES, do not need this gradient information and hence, can be used in this case. However, restricting the real locations to integer values does introduce an error in the optimized location. The error depends on the gridblock size, and is acceptable in an engineering sense. \citep{sarma_efficient_2008} and \citep{forouzanfar_optimization_2012} give ways to deal with this discontinuity in the NPV while optimizing well placement.




\subsection{Well control optimization} 
\label{sec:23contolop}

The well control optimization problem aims to determine the optimal time-varying well setting for each of the production and injection wells. The optimization problem can be stated as follows: 
 
\begin{eqnarray}
\max\limits& \ \mathrm{NPV}=J(\bf u) \\
\mathrm{subject \ to} & \ \mathbf{u_{lb}} \leq \mathbf{u} \leq \mathbf{u_{ub}} ,
\end{eqnarray}
where $\mathbf{u}$ denotes the well control vector. Each well can be controlled either by well rate or by bottom hole pressure (BHP). The time-varying well controls are represented by piecewise functions in time with $N_t$ intervals. The number of variables for this problem is $N_tN$.

As in well placement optimization problem, only bound constraints are considered explicitly here. 



\subsection{Joint well placement and control optimization} 
\label{sec:24jointop}


The joint problem optimizes both well locations and controls. The optimization problem is defined as follows:
\begin{eqnarray}
\max\limits& \ \mathrm{NPV}=J(\bf v, u) \\
\mathrm{subject \ to} & \ \mathbf{v_{lb}} \leq \mathbf{v} \leq \mathbf{v_{ub}} \\
&\ \mathbf{u_{lb}} \leq \mathbf{u} \leq \mathbf{u_{ub}} ,
\end{eqnarray}
where $\mathbf{v}$ and $\mathbf{u}$ denote the well placement and control vectors.

Two procedures are commonly used for joint well placement and control optimization---a simultaneous procedure or a sequential procedure. The simultaneous procedure optimizes well locations and controls simultaneously, hence the number of optimization variables, $2N+N_tN$ for vertical wells and $6N+N_tN$ for 3-D angled wells, is larger than the individual problems, which makes the optimization more difficult.

In the sequential procedure, well placement is optimized first using some reasonable control scheme. The controls are then optimized for the wells using the best locations found in the first stage. These two stages may be repeated. The sequential procedure decouples the joint problem into two separate subproblems, and the difficulty for each subproblem is decreased. The number of optimization variables for the well placement stage is either $2N$ or $6N$ (for vertical or angled well), and is $N_tN$ for the well control stage.
The joint problem is worth studying because the optimal location of each well depends on how the well is operated and vice-versa \citep{humphries_simultaneous_2013,isebor_generalized_2014,humphries_joint_2015}. 
Furthermore, \citep{forouzanfar_two-stage_2010} point out that the results obtained also depend directly on the specified reservoir life time, which makes the joint optimization problem even more complicated. We use a predefined reservoir life for examples in this paper.

\section{Optimization Methods Considered}
\label{sec:3_solver}

\subsection{Multilevel Coordinate Search}
\label{sec:3_mcs}

MCS, first proposed by Huyer and Neumaier \citep{huyer_global_1999}, was inspired by DIRECT \citep{jones_lipschitzian_1993}. MCS is a branching scheme which searches by recursively splitting hyperrectangles.
Like DIRECT, MCS is a mathematical programming approach which provides a systematic search within the bound constraints (the bounds can be infinite for MCS). MCS builds upon DIRECT by introducing a multilevel mechanism which allows a balanced global and local search. DIRECT has no local search capability.

Levels are assigned as an increasing function of the number of times a box has been split. The global search portion of the algorithm is accomplished by splitting boxes that have not been searched often -- those with a low level. Within a level the boxes with the best function values are selected to complete a local search. The local search builds a quadratic model, determines a promising search direction and performs a line search. This allows for quicker convergence while the global part of the algorithm identifies a region near the global optimum. 

MCS allows for a more irregular splitting than DIRECT, giving preference to regions with low function values. Convergence to a global optima is guaranteed as the number of levels goes to infinity if the objective function is continuous around the global optimizer.

\citep{huyer_global_1999} reports that MCS works well in problems where the global optimum can be constrained by finite bound constraints. 
\citep{posik_comparison_2012} report very good performance in the early search phase with a small budget of objective function evaluations.

MCS provides numerous heuristic enhancements over DIRECT. Consider a $n$-dimensional bound constrained minimization problem
\begin{eqnarray}
\min& \ f(\bf x) \\
\mathrm{subject \ to} & \ \bf u \leq \bf x \leq \bf v. 
\end{eqnarray}
where $f$ denotes the objective function; ${\bf x}=[x_1,\cdots,x_n]$ denotes the optimization variables; ${\bf u}=[u_1,\cdots,u_n]$ and ${\bf v}=[v_1,\cdots,v_n]$ are the lower and upper bound, respectively.
The pseudocode of the basic steps of MCS are described in Alg. \ref{algomcs}. 
A complete description of the algorithm is quite complex and can be found in \citep{huyer_global_1999}. 
For the experiments in this paper, we used the implementation of MCS from \citep{software_mcs}.
The algorithm was originally designed to minimize a function $f$. To maximize a function $f$ we simplify minimize $-f$.


During the initialization portion of the algorithm, MCS accepts an initialization list which is used to produce an initial set of boxes partitioning the search space.
For $i$th coordinate, an initialization list stores $L_i$ points $x_i^1,x_i^2,\cdots,x_i^{L_i}\in [u_i,v_i]$ ($L_i\geq 3$). Users can incorporate their knowledge and experience to choose points with a high possibility of obtaining good solution.
MCS continues to process and split boxes until some boxes reach level $s_{\max}$. 
A box can be split either by rank or by expected gain, depending on the relationship between the box level $s$ and number of split times at each coordinate. 
If a box has been split many times and reached a high level, MCS selects the coordinate which has been split the fewest times and splits along this coordinate (split by rank).
If the level of a box is not high, MCS splits along a coordinate with the maximum expected gain according to a quadratic model obtained by fitting function values (split by expect gain).
The parameter $s_{\max}$ controls the precision of the global search phase before any local search would be attempted. MCS also has the option to turn off the local search phase.

\begin{algorithm}
\SetKwInOut{Input}{Input}
\SetKwInOut{Output}{Output}

\Input{Objective function $f$, bound constraints $\mathbf u, \mathbf v$}
\Input{Initialization list $x_i^j\ (j=1,\cdots,L_i,i=1,\cdots,n)$, maximum level $s_{\max}$, local search state: on/off}
\For{$i=1$ \KwTo $n$}{
  $\mathbf x\leftarrow$ the best of $\{\mathbf x^j\}_{j=1}^{L_i}$ (in terms of objective function value), where $\mathbf x^j$ is $\mathbf x$ with $x_i$ changed to $x_i^j$\;
  Split the current box $B$ along the $i$th coordinate at $x_i^j$ and the points determined by the golden ratio\;
  $B\leftarrow$ the one has best function value of the boxes containing $\mathbf x$\;
}
\While{there are boxes of level $s<s_{\max}$}{
  \For{all non-empty levels $s=2$ \KwTo $s_{\max}-1$}{
    Choose the box $B$ at the level $s$ with the best function value\;
    $i\leftarrow$ number of split times the coordinate used least often when producing box $B$\;
    \If(\tcp*[h]{Split by rank}){$s>2n(i+1)$}{
      Split the box $B$ along the $i$th coordinate\;
    }
    \ElseIf(\tcp*[h]{Split by expected gain}){$s\leq 2n(i+1)$}{
      Determine the most promising splitting coordinate $i$\;
      Compute the minimal expected function value $f_{\mathrm{exp}}$ at new point\;
      \eIf{$f_{\mathrm{exp}}<f_{\mathrm{best}}$}{
        Split $B$ along the $i$th coordinate\;
      }{
        Tag $B$ as not promising and increasing its level by 1\;
      }
    }
  }
  \For{Base points, $\mathbf x$, of all the new boxes at level $s_{\max}$}{
    Start a local search from $\mathbf x$ if improvement is expected\;
  }
}
\Output{$\mathbf x^{\mathrm{best}}$,$f_{\mathrm{best}}$}
\caption{The MCS algorithm \citep{posik_comparison_2012}\label{algomcs}.}
\end{algorithm}

\begin{figure}[htb]
  \centering 
  \subfigure[Initialization]{ 
    \label{fig:subfig:mcs_init} 
    \includegraphics[width=0.45\textwidth]{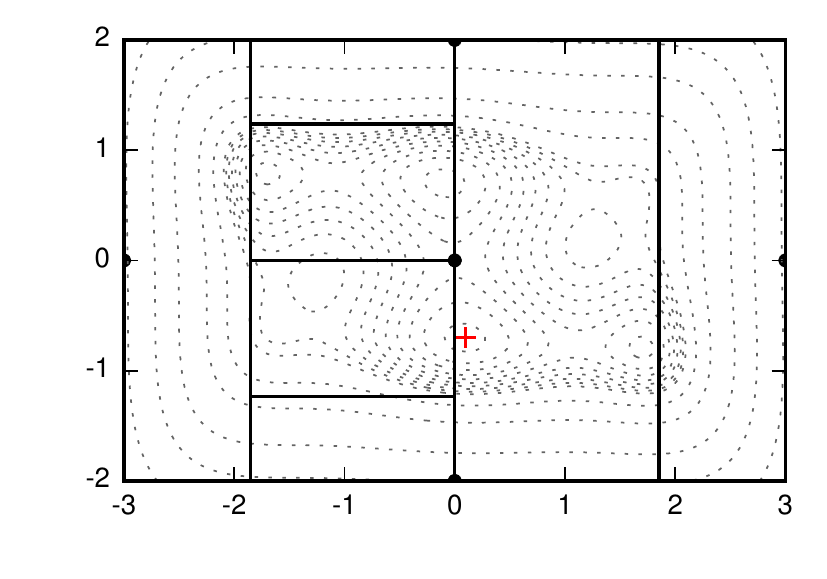}} 
  \subfigure[Splitting]{ 
    \label{fig:subfig:mcs_split} 
    \includegraphics[width=0.45\textwidth]{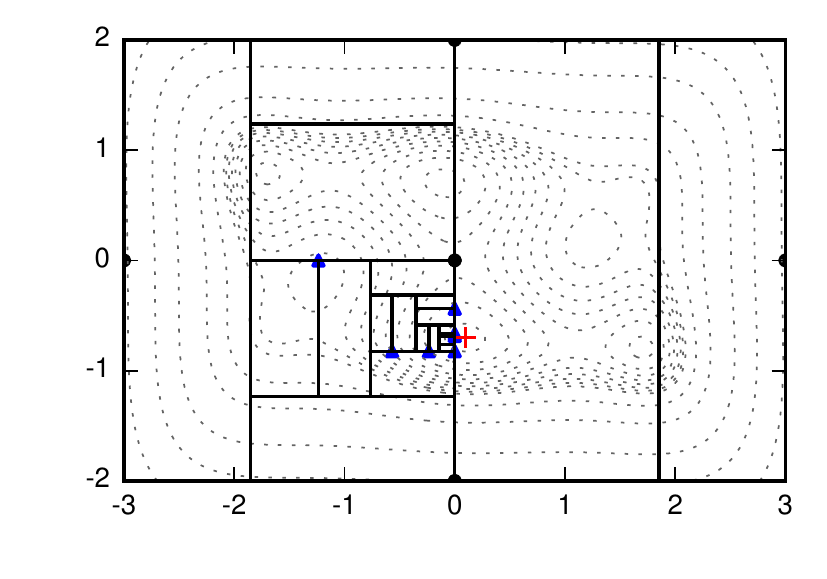}} 
  \subfigure[Local search]{ 
    \label{fig:subfig:mcs_ls} 
    \includegraphics[width=0.45\textwidth]{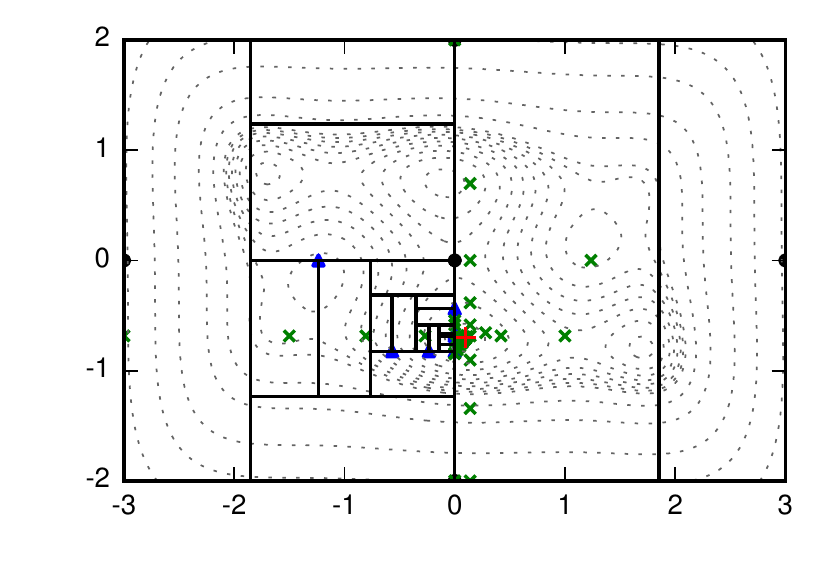}}  
  \caption{The boxes and test points of MCS for the six-hump camel function. The dashed lines are the contour lines of the function. The global optima is known and is indicated by `+'. The dots are points tested during the initialization procedure, while the triangle marker and the `x' marker are the points tested during the splitting procedure and the local search procedure, respectively.} 
  \label{fig:mcs_procedure} 
\end{figure}

We provide a simple example to see how MCS works. We consider the objective function $f=x_1^2(4-2.1x_1^2+x_1^4/3)+x_1 x_2-4x_2^2(1-x_2^2)$, which is a six-hump camel function with 2 unknowns. 
The bounds for the 2 unknowns are $x_1\in[-3,3]$ and $x_2\in[-2,2]$.
The global minimizer for this function is $[0.0898,-0.7127]$ and the global minimum value is $-1.0316$. 
We choose the default parameter settings for MCS: 
$L_i=3$ for $i=1,2$, 
the initialization list $x_i^1=u_i$, $x_i^2=(u_i+v_i)/2$, $x_i^3=v_i$, a maximum number of level $s_{\max}=5n+10=20$, and we turn the local search phase on. Fig.\ \ref{fig:mcs_procedure} shows a loop of MCS for this problem. 

\begin{figure}[htb]
\centering
  \includegraphics[width=0.45\textwidth]{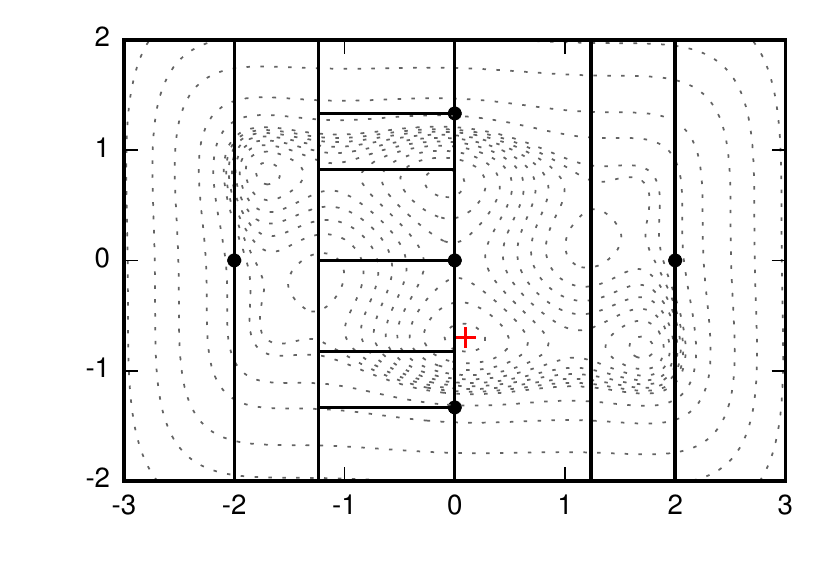}
  \caption{The boxes obtained after the initialization procedure of MCS for the six-hump camel function using the initialization list $x_i^1=(5u_i+v_i)/6$, $x_i^2=(u_i+v_i)/2$, $x_i^3=(u_i+5v_i)/6$.} 
  \label{fig:mcs_algoinit2} 
\end{figure}

Fig.\ \ref{fig:subfig:mcs_init} presents the boxes after the initialization procedure (lines 1--4 of Alg. \ref{algomcs}). By using the default initialization list, MCS first splits the root box along the $x$-coordinate at the midpoint, the two boundary points, and the points between determined by the golden ratio. Then MCS chooses the new box that has the highest estimated variability and splits it along the $y$-coordinate. We note that the initialization list can also be specified by the user. Different initialization lists results in a different split of the boxes. One other commonly used initialization list is $x_i^1=(5u_i+v_i)/6$, $x_i^2=(u_i+v_i)/2$, $x_i^3=(u_i+5v_i)/6$. By using this initialization list, the boxes after the initialization procedure are shown in Fig.\ \ref{fig:mcs_algoinit2}.

The initialization list can also be generated automatically with the aid of line searches. 
Starting from the point $\mathbf{x}_0=min |[\mathbf{u},\mathbf{v}]|$, for $i$th coordinate, we complete line searches along this coordinate to find up to $nloc$ minima within $smaxls$ function evaluations. All local minimizers found by line searches are put into the initialization list. If the number of minima is less than $nloc$, then we use the points closest to $u_i$ and $v_i$ obtained by line searches to supplement the initialization list. The best point is taken as the start point for line searches for the next coordinate.
We choose $nloc=5$ and $smaxls=25$ as recommend by \citep{huyer_global_1999} for all examples in this paper. With this setting, MCS needs up to $25n$ function evaluations to generate the initialization list. 
This stage will slow down the convergence of the optimization in the early stages, however, such an initialization list can ultimately improve the performance of MCS significantly. 
Hence, to get the advantages of the line search, the total number of function evaluations needs set to a number much larger than $25n$.

\begin{figure}[htb]
\centering
  \includegraphics[width=0.45\textwidth]{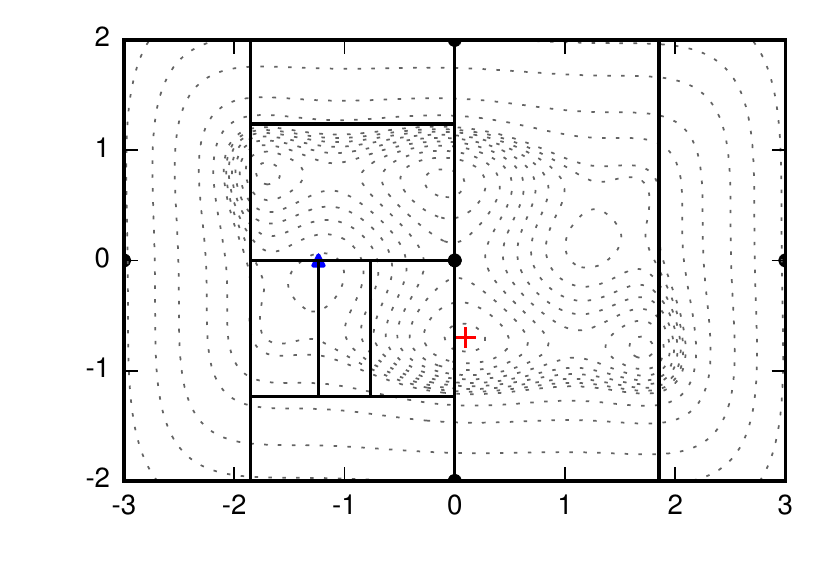}
  \caption{The boxes obtained after the splitting procedure of MCS for the six-hump camel function with $s_{\max}=5n$.} 
  \label{fig:mcs_split2} 
\end{figure}

After the initialization procedure, the search space will be further split until one of the boxes reaches the maximum level $s_{\max}$. Fig.\ \ref{fig:subfig:mcs_split} shows the boxes after the splitting procedure. As mentioned above, $s_{\max}$ decides the depth to which MCS explores a region and hence controls the precision of the global search phase. If $s_{\max}=5n$, then the boxes obtained are shown in Fig.\ \ref{fig:mcs_split2}.

Once a box reaches the maximum level, a local search starts from its base point. Fig.\ \ref{fig:subfig:mcs_ls} shows the points evaluated by the local search. 
The local search stops if the maximum number of steps in local search is reached, or the stopping criteria, $|f-f_{0}|< \gamma |f_{0}|$, where $f$ and $f_{0}$ are the current and the smallest values of objective function, respectively, is triggered. The maximum number of steps is 50, and $\gamma=0.01$, by default.
After that, MCS will cycle back to the split procedure.

MCS was originally designed for either unconstrained optimization problems or problems with bound constraints only. 
In general, there are many natural ways to attempt to extend unconstrained optimization algorithms to handle constraints.  For example, the penalty function method \citep{bouzarkouna_well_2012,ciaurri_derivative-free_2011}, Lagrange multiplier method \citep{chen_robust_2012,brouwer_dynamic_2004}, and so on. 
To the best of our knowledge, there has been no such paper which has addressed this topic for MCS. 

\subsection{Configurations of MCS considered}
\label{sec:3_conf}
In order to analyse the sensitivity of the parameters in the MCS algorithm, we apply MCS, with 7 different settings of the parameters, to our examples. The 7 settings used are:

\begin{itemize}
\item[\textbullet] MCS-1: MCS with its default settings from \citep{huyer_global_1999}. A simple initialization list is used consisting of midpoints and boundary points, i.e. $x_i^1=u_i,\ x_i^2=(u_i+v_i)/2,\ x_i^3=v_i$. The number of levels is chosen as $s_{\max}=5n+10$, where $n$ is the dimension of the problem. The maximal number of visits in the local search is 50, and the acceptable relative accuracy for local search is $\gamma=0.01$.
\item[\textbullet] MCS-2: MCS with the initialization list $x_i^1=(5u_i+v_i)/6,\ x_i^2=(u_i+v_i)/2,\ x_i^3=(u_i+5v_i)/6$. Unlike the initialization list in MCS-1, the points here are uniformly spaced but do not include the boundary points. The other settings are same as in MCS-1.
\item[\textbullet] MCS-3: MCS with an auto-generated initialization list. In MCS-3, we first perform a sequence of line searches along all coordinate directions to generate the initialization list. The other settings are same as in MCS-1.
\item[\textbullet] MCS-4: MCS with the initialization list $x_i^1=u_i,\ x_i^2=x_0,\ x_i^3=v_i$. Unlike the initialization list in MCS-1, we use an user defined initial guess $x_0$ instead of the midpoints. The other settings are same as in MCS-1.
\item[\textbullet] MCS-5: MCS with the initialization list $x_i^1=(5u_i+v_i)/6,\ x_i^2=x_0,\ x_i^3=(u_i+5v_i)/6$. Unlike the initialization list in MCS-2, we use an user defined initial guess $x_0$ instead of the midpoints. The other settings are same as in MCS-2.
\item[\textbullet] MCS-6: MCS with a larger maximum number of levels, $s_{\max}=10n$. This is chosen to attempt to improve the global search phase. The other settings are same as in MCS-4.
\item[\textbullet] MCS-7: MCS without the local search phase. In MCS-7, we set the maximal number of visits in the local search to 0. The other settings are same as in MCS-4.
\end{itemize}



\subsection{Other algorithms considered}
\label{sec:3_otheralgo}

For comparison
three algorithms -- generalized pattern search (GPS), particle swarm optimization (PSO) and covariance matrix adaptation evolution strategy (CMA-ES) -- are used. 

Generalized pattern search (GPS)  \citep{audet_analysis_2002,torczon_convergence_1997,yin_extended_2000} is a deterministic local search algorithm. It does not require gradients and hence, it can be used on problems that are not continuous or differentiable. 
For the parameter settings, we use a $2n$ positive spanning set, where $n$ is the dimension of the search space. The expansion factor is set to 2, and the contraction is set to 0.5 \citep{torczon_convergence_1997,audet_analysis_2002,kolda_optimization_2003}. 

Particle swarm optimization (PSO)  \citep{kennedy_particle_2011,vaz_particle_2007} is a population-based stochastic search method. PSO's search mechanism mimics the social behavior of biological organisms such as a flock of birds.
PSO can search a very large space of candidate solutions, which reduces the chance of getting trapped at an unsatisfactory local optimum. 

The performance of PSO depends on the values assigned to the algorithm parameters. Following the work of \citep{perez_particle_2007}, our implementation of PSO uses the population size of 50, and the weighting parameters $\omega =0.9$, $c_1=0.5$, and $c_2=1.25$. The best parameters are usually problem dependent. Further tuning, for a specific problem, will likely yield superior performance. Further discussion on this issue can be found in \citep{clerc_stagnation_2006} and \citep{onwunalu_optimization_2010}. 

The covariance matrix adaptation strategy (CMA-ES)  \citep{hansen_evaluating_2004,loshchilov_cma-es_2013,auger_restart_2005} is a population-based stochastic optimization algorithm. Unlike a genetic algorithm (GA), PSO, and other classic population-based stochastic search algorithms, candidate solutions of CMA-ES are sampled from a probability distribution which are updated iteratively. 
For CMA-ES, we use the settings from \citep{hansen_evaluating_2004} (See Table \ref{tab:cmaes_para}). In fact, according to their work, CMA-ES doses not require significant parameter tuning for its application.

\begin{table}[htb]
\caption{Strategy parameter values used in CMA-ES.}
\label{tab:cmaes_para}
\centering
\begin{tabular}{cc}
\hline\noalign{\smallskip}
Parameter & Value \\
\noalign{\smallskip}\hline\noalign{\smallskip}
$\lambda$ & $4+\lfloor 3 \ln (n) \rfloor $\\
$\mu$ & $\lfloor \lambda/2 \rfloor$ \\
$c_c $ & $ \frac{4}{n+4} $ \\
$c_\sigma $ & $\frac{\mu_{\rm{eff}}+2}{n+\mu_{\rm{eff}}+3} $ \\
$d_\sigma $ & $1+2\max\left(0,\sqrt{\frac{\mu_{\rm{eff}}-1}{n+1}}-1\right)+c_\sigma $ \\
$\mu_{\rm{cov}} $ & $\mu_{\rm{eff}} $ \\
$c_{\rm{cov}} $ & $\frac{1}{\mu_{\rm{cov}}}\frac{2}{(n+\sqrt{2})^2}+\left(  1-\frac{1}{\mu_{\rm{cov}}}\right)\min\left(1,\frac{2\mu_{\rm{eff}}-1}{(n+2)^2 +\mu_{\rm{eff}}} \right) $ \\
\noalign{\smallskip}\hline
\end{tabular}
\end{table}

PSO and CMA-ES are stochastic algorithms and the result of each trial is different. Thus, in order to assess the overall performance of these algorithms, we run each of the algorithms many times for each test example. 

The objective function evaluations for the well placement and control optimization problems require the evaluation of a numerical  oil reservoir simulator. The cost of the total 
optimization process is completely dominated by the cost of each simulation evaluation. 
The time spent in the nuts and bolts of the optimization algorithm itself can be neglected. 
As a result we choose to use the number of simulation runs as our performance indicator to compare the optimization strategies. This is a widely used measure in the well placement and control optimization literature 
\citep{isebor_generalized_2014,forouzanfar_joint_2014,brouwer_dynamic_2004,humphries_simultaneous_2013}.

It is also worth mentioning that, GPS and CMA-ES need an initial guess to start the optimization processes. PSO can generate an initial population automatically.
In our four examples, we use a physically reasonable initial guess for GPS, PSO, and CMA-ES.
For MCS, the initial point is determined by the initialization list.
Either the initialization list includes the initial guess (MCS-4, MCS-5) or an ordinary initialization list is used (MCS-1, MCS-2, MCS-3).

\subsection{Algorithm combinations for the joint optimization problem}
\label{sec:3_algocombforjoint}

For the joint well placement and control optimization problem, we consider both a simultaneous procedure and a sequential procedure. The simultaneous procedure optimizes over the well locations and controls simultaneously.
To solve the joint problem with a simultaneous approach we consider the 7 different configurations of MCS, and also GPS, PSO, and CMA-ES.

The sequential procedure divides the optimization process into a well placement optimization stage and a well control optimization stage. Each stage is an independent optimization problem and can be optimized using the same or different algorithms. We label the approaches used for the sequential procedure in the form \emph{Algorithm1-Algorithm2}, where Algorithm1 denotes the algorithm used for the well placement optimization stage and Algorithm2 denotes the algorithm used for the well control optimization stage. Many such combinations are possible.

The combinations considered in this paper are divided into 3 groups. The first group includes MCS-MCS, GPS-GPS, PSO-PSO, and CMA-ES-CMA-ES, which use the same algorithm in both the well placement stage and the well control stage. The second group includes MCS-GPS, MCS-PSO, and MCS-CMA-ES, which use MCS for the well placement stage. The third group GPS-MCS, PSO-MCS, and CMA-ES-MCS  uses MCS for the well control stage.

\section{Examples and Results}
\label{sec:4_exam_resu}

\subsection{Example 1, vertical well placement optimization}
\label{sec:41_vwop}

\subsubsection{Reservoir model description}

\begin{figure*}[htb]
\centering 
\subfigure[Tops]{ 
    \label{fig:subfig:punqs3_tops} 
    \includegraphics[width=0.45\textwidth]{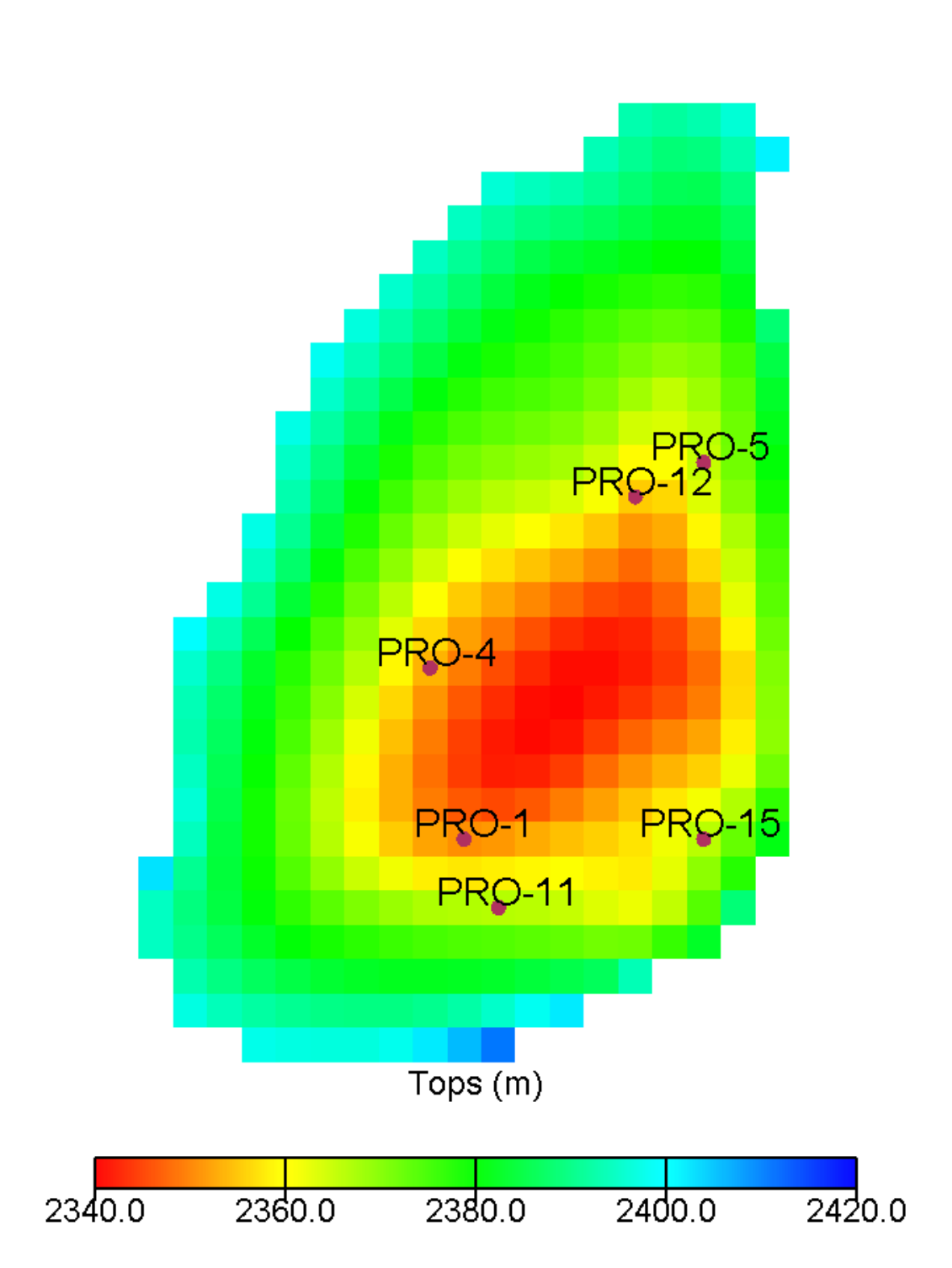}
    }
\subfigure[Permeability]{ 
    \label{fig:subfig:punqs3_perm} 
    \includegraphics[width=0.45\textwidth]{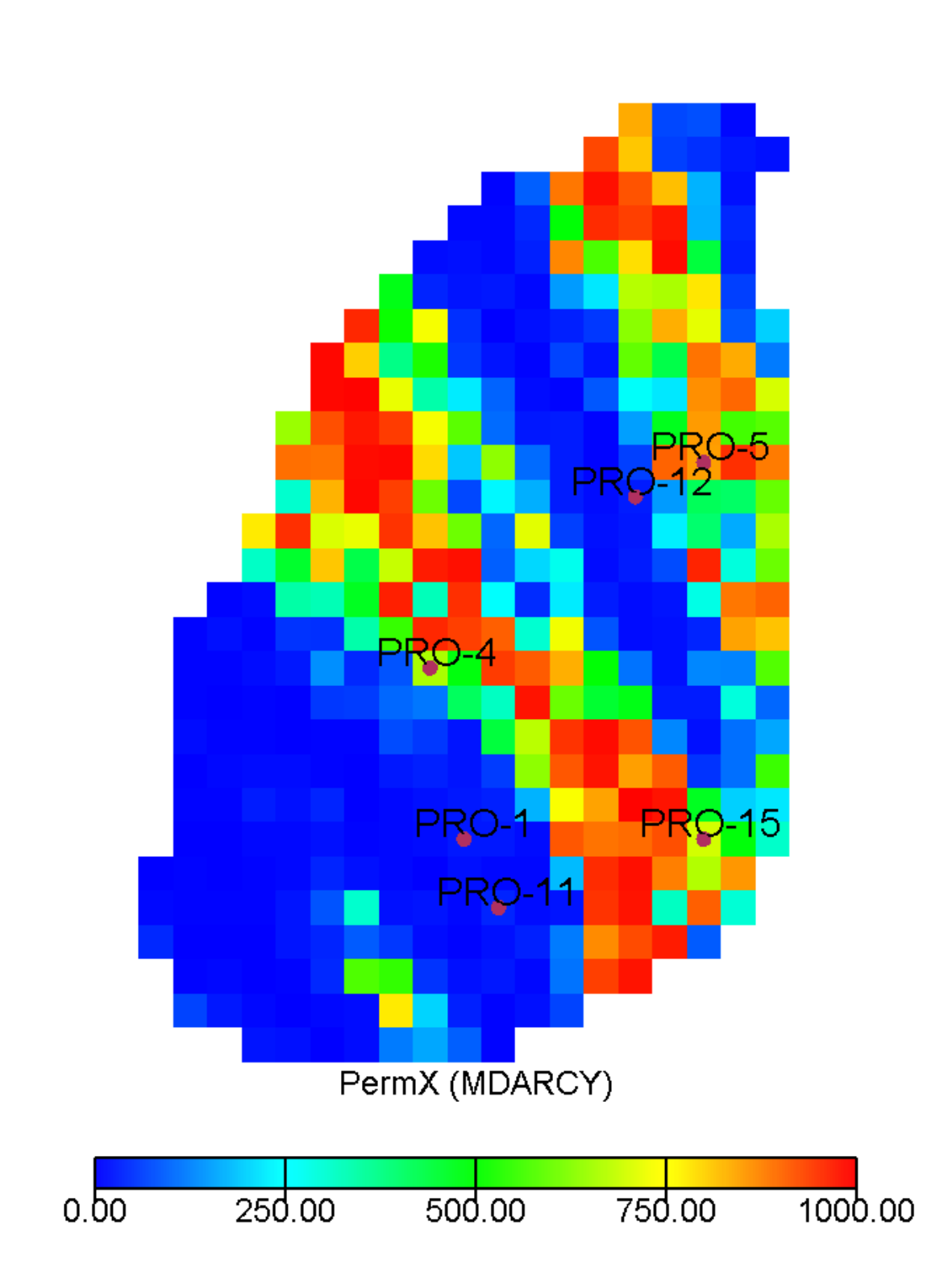}
    }
\caption{Properties and initial well locations for the PUNQ-S3 model used in Example 1.}
\label{fig:punqs3}
\end{figure*}

The first example uses the PUNQ-S3 model, which is a small reservoir model based on an actual North Sea reservoir \citep{gao_quantifying_2006}. The model uses $19 \times 28 \times 5$ grid blocks with $\Delta x= \Delta y =180 $m and 1761 active grid blocks. The simulation model involves a three-phase gas-oil-water flow. The field initially contains 6 production wells and no injection wells due to the strong aquifer. Fig.\ \ref{fig:punqs3} shows the depth of the top face and permeability field, together with the initial well locations of PUNQ-S3 model. The reservoir production time is 20 years, the bottom hole pressure of each well is fixed at 200 bar. 

We seek to optimize the well locations of all 6 wells. The formulation of the well placement optimization problem is given in Section \ref{sec:22placementop}. The objective function is the net present value (NPV). The simulator used to predict the production dynamics (the flow rates of the gas, oil, water produced and water injected) is Eclipse \citep{geoquest_eclipse_2014}, a commercial reservoir simulator from Schlumberger Ltd. The economic parameters used to calculate NPV are given in Table \ref{tab:eco_set_1}. 

Every well has two positional variables which gives a total of 12 optimization variables. Only bound constraints are considered in this example. 
We force $1\leq x \leq 19$ and $1\leq y \leq 28$ for all 6 wells.

The optimization problem was solved by using GPS, PSO, CMA-ES, and all 7 configurations of MCS. The maximum number of simulation runs is set to 600.


\begin{table}[htb]
\caption{Economic parameters used for Example 1.}
\label{tab:eco_set_1}
\centering
\begin{tabular}{cc}
\hline\noalign{\smallskip}
Parameter & Value \\
\noalign{\smallskip}\hline\noalign{\smallskip}
Gas revenue (USD/$\rm m^3$) & 0.5  \\
Oil revenue (USD/$\rm m^3$) & 500.0 \\
Water-production cost (USD/$\rm m^3$) & 80.0 \\
Annual discount rate & 0 \\
\noalign{\smallskip}\hline
\end{tabular}
\end{table}

\subsubsection{Results and discussion}


The results of Example 1 are shown in Table \ref{tab:e1_result}. In this table, the final NPV after 600 simulation runs for each algorithm is given. Moreover, for PSO and CMA-ES, the maximum, minimum, mean, median, and standard deviation of NPV are given.
From the table we can see that GPS obtains the highest NPV value after the 600 simulation runs. Though the maximum NPV for PSO and CMA-ES are slightly higher than MCS in some cases, MCS generally performs better than PSO and CMA-ES when compared to the mean and the median NPV. 

\begin{table}[htb]
\caption{Results for Example 1. Values shown are NPV in $\$ \times 10^9$ USD. }
\label{tab:e1_result}
\centering
\subtable[Deterministic algorithms (MCS, GPS)]{
\label{tab:e1_result_det}
\begin{tabular}{p{2cm}|p{8cm}}
\hline
Algorithm & \hfil NPV \\
\hline
MCS-1  & \hfil 2.10 \\
MCS-2  & \hfil 2.16 \\
MCS-3  & \hfil 2.21 \\
MCS-4  & \hfil 2.17 \\
MCS-5  & \hfil 2.15 \\
MCS-6  & \hfil 2.20 \\
MCS-7  & \hfil 2.11 \\
GPS    & \hfil 2.32 \\
\hline
\end{tabular}
}

\subtable[Stochastic algorithms (PSO, CMA-ES)]{        
\label{tab:e1_result_sto}
\begin{tabular}{p{2cm}|p{1cm}|p{1cm}p{1cm}p{1cm}p{1cm}p{1cm}}
\hline
Algorithm & Trials & Max & Min & Mean & Median & Std. \\
\hline
PSO    &10 &2.18&2.00&2.07&2.04&0.07\\
CMA-ES &10 &2.22&1.91&2.03&2.00&0.12\\
\hline
\end{tabular}
}
\end{table}

\begin{figure}[htb]
\centering
  \includegraphics[width=0.45\textwidth]{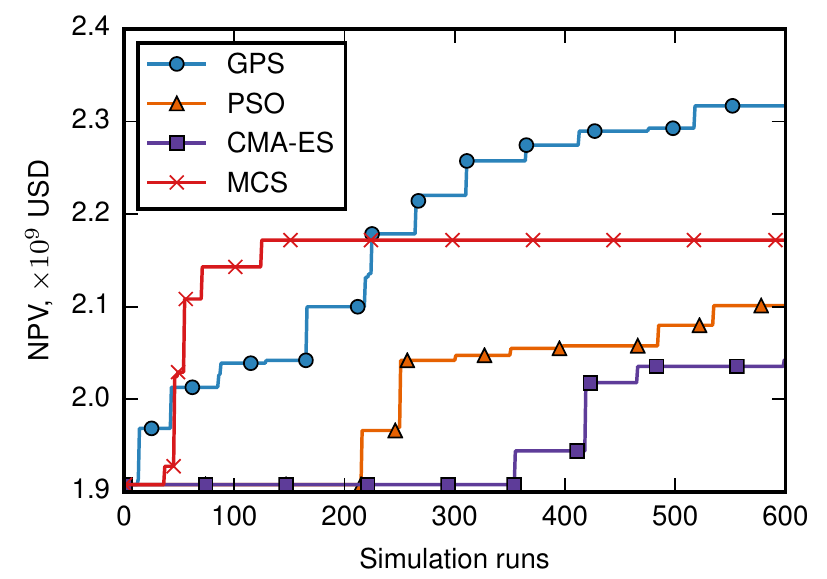}
\caption{Optimization performance for Example 1. For PSO and CMA-ES, the solid lines depict the median NPV. MCS here is a label for MCS-4.}
\label{fig:e1_algo}
\end{figure}

Plots of the NPV for the four algorithms versus the number of simulation runs are shown in Fig.\ \ref{fig:e1_algo}. Note that for PSO and CMA-ES, 10 trials are performed and the solid lines depict the median NPV over all 10 trials.
Since GPS and MCS
are deterministic algorithms, only one trial is performed.

From Fig.\ \ref{fig:e1_algo} we can see that MCS showed excellent convergence speed at the early stage of optimization (simulation runs $<$200). GPS converges slower than MCS at an early stage, but eventually GPS obtains the highest NPV. The two stochastic algorithms, PSO and CMA-ES, showed slow convergence speed at an early stage of the optimization process.

\begin{figure}[htb]
  \centering 
  \subfigure[PSO]{ 
    \label{fig:subfig:e1_pso} 
    \includegraphics[width=0.45\textwidth]{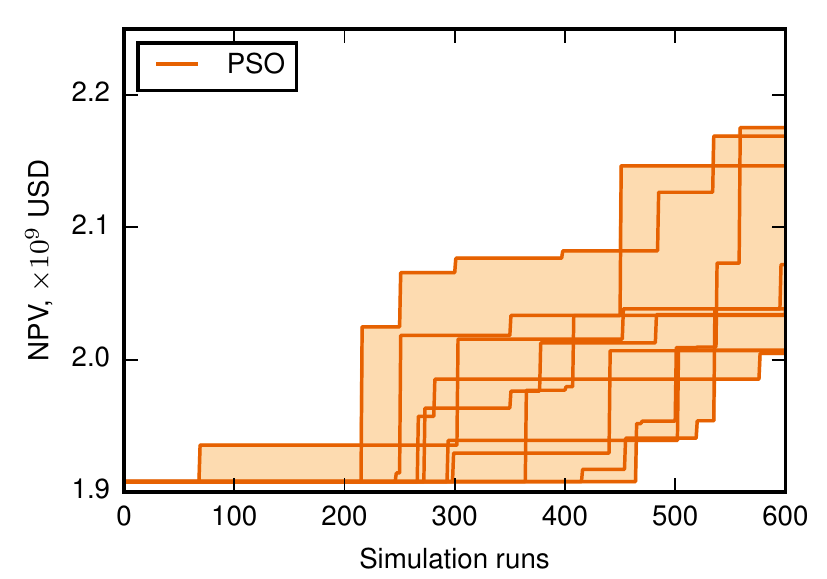}} 
  \subfigure[CMA-ES]{ 
    \label{fig:subfig:e1_cmaes} 
    \includegraphics[width=0.45\textwidth]{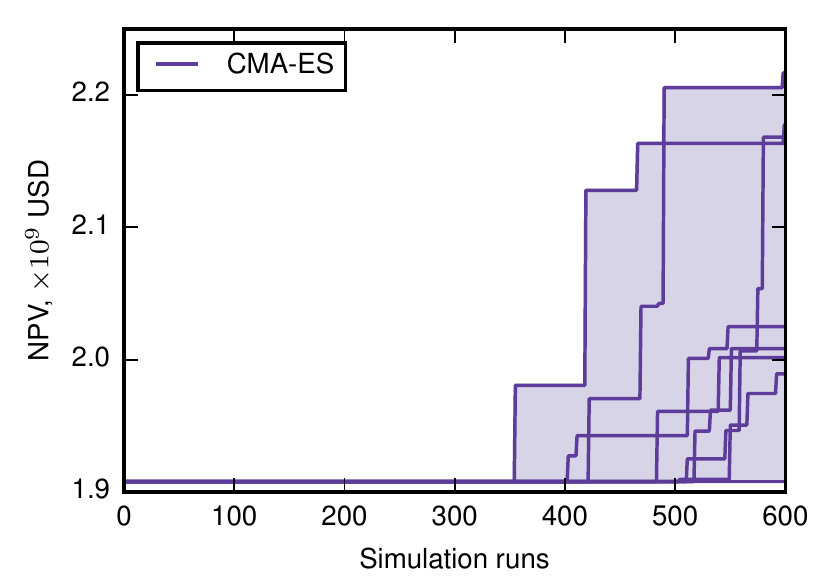}} 
  \caption{The range of NPV found amongst the trials of PSO and CMA-ES for Example 1. Each solid line represents a trial.} 
  \label{fig:e1_psocmaes} 
\end{figure}

Since PSO and CMA-ES are stochastic algorithms, the performance is different for each trial. Fig.\ \ref{fig:e1_psocmaes} shows the range of NPV found amongst the trials of PSO and CMA-ES. In this figure, the areas between the maximum and minimum NPV are shaded for PSO and CMA-ES. It is clear that the NPV obtained by PSO and CMA-ES has a high variation for this example. This suggests that when solving this problem by PSO or CMA-ES, a single trial has a high risk to obtain an unsatisfactory NPV.

\begin{figure}[htb]
  \centering 
  \subfigure[Initialization list]{ 
    \label{fig:subfig:e1_mcs_init} 
    \includegraphics[width=0.45\textwidth]{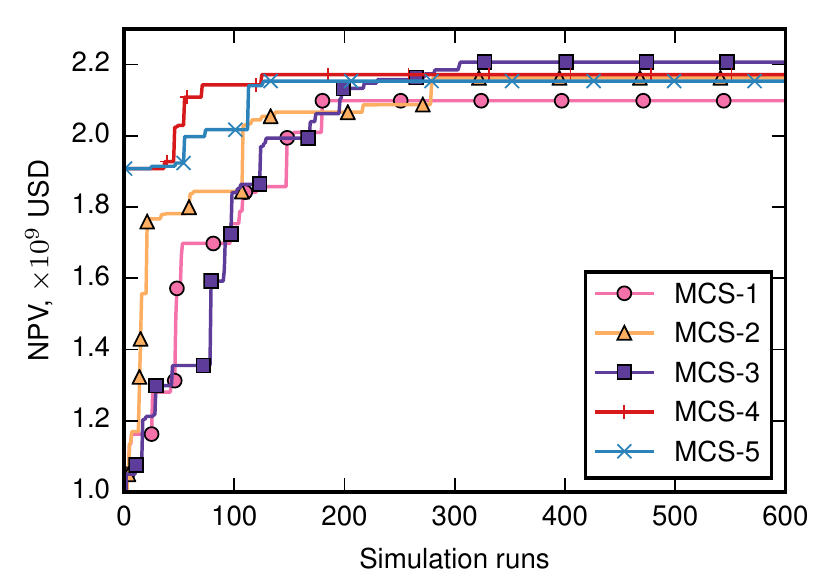}} 
  \subfigure[Levels]{ 
    \label{fig:subfig:e1_mcs_smax} 
    \includegraphics[width=0.45\textwidth]{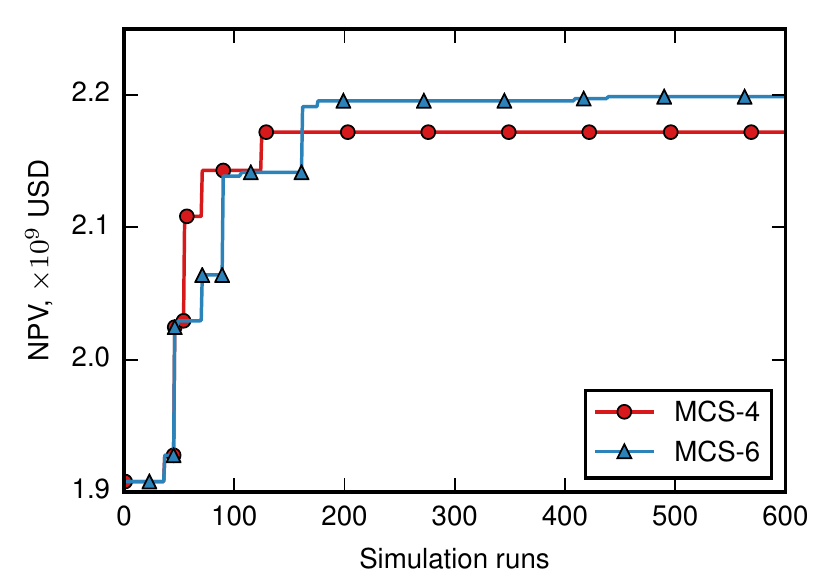}} 
  \subfigure[Local search]{ 
    \label{fig:subfig:e1_mcs_ls} 
    \includegraphics[width=0.45\textwidth]{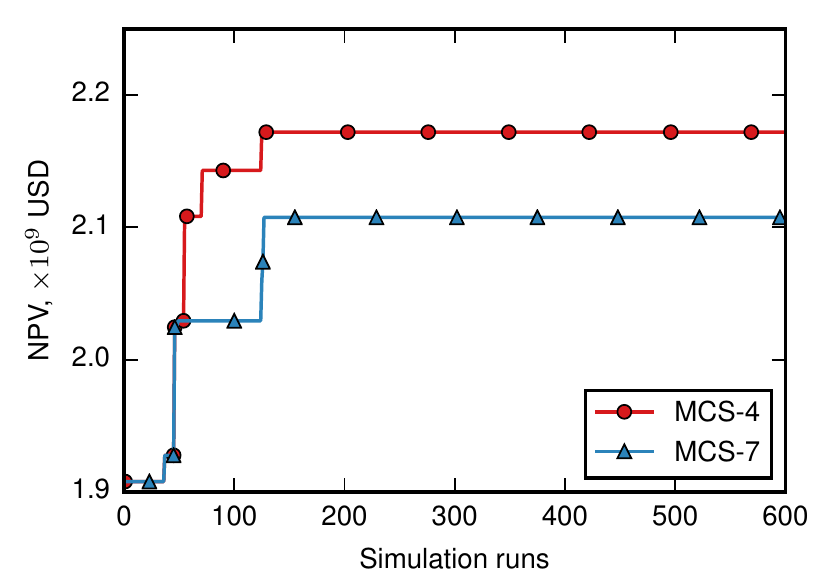}}   
  \caption{The performance of different configurations of MCS for Example 1. Three important parameters in MCS: initialization list, number of levels, and local search.} 
  \label{fig:e1_mcs} 
\end{figure}

The 7 different MCS configurations are also tested with this example. Detailed results are shown in Fig.\ \ref{fig:e1_mcs}. 
These algorithms are divided into 3 groups. The first group (MCS-1, MCS-2, MCS-3, MCS-4, and MCS-5) uses different initialization lists. This allows us to check the impact of the initialization list. The second group (MCS-4, MCS-6) uses a different maximum number of levels, $s_{\max}$. The higher $s_{\max}$, the better the global search ability. The third group (MCS-4, MCS-7) is used to analyse the role of local search in MCS. 

From the first group  
MCS-3 ultimately achieves the highest NPV, followed by MCS-4 and MCS-2, and then MCS-5 and MCS-1. 
The ultimate difference in NPV between the seven configurations is small (about 5\%). 
Starting from a relatively low NPV, MCS-1 obtains a NPV slightly smaller than MCS-4.
The NPVs obtained by MCS-2 and MCS-5 are similar.
This shows that MCS with a good initial guess in the initialization list has an advantage over the others, but the advantage is very slight with a large computational budget. MCS-4 and MCS-5, starts the optimization with a relatively high NPV, and has a significant advantage over MCS-1 and MCS-2 when the computational budget is limited.

\begin{figure}[htb]
\centering
  \includegraphics[width=0.45\textwidth]{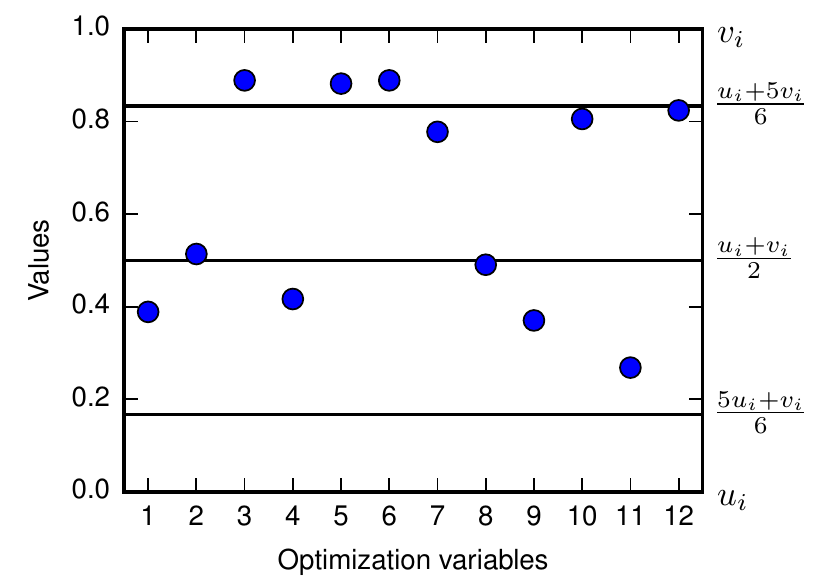}
\caption{Normalized boundary, initialization lists, and the optimal solution for Example 1.}
\label{fig:e1_ila}
\end{figure}

Comparing MCS-1 and MCS-2, MCS-2 converges faster than MCS-1, and obtains higher NPV than MCS-1 ultimately. This indicates that the uniformly spaced initialization list without boundary points is more suitable. To explore this, we normalize the search space to the $[0,1]$-interval and map the initialization lists and the global optima to the normalized search space in Fig.\ \ref{fig:e1_ila}.
It is clear that the optimal solution is aligned better with initialization list II (MCS-2), which explains its better performance for this problem.
The optimal solution is not known a priori in most cases, so although a suitable initialization list can improve the performance of MCS, it is difficult to choose between MCS-1 and MCS-2 a priori. 
For MCS-4 and MCS-5, the difference in the ultimate NPV is small. This indicates that with a good initial guess, the importance of the initialization list decreases.

The second group, MCS-4 and MCS-6, compares the performance of MCS with different specified maximum levels $s_{\max}$. MCS with a larger number of maximum levels, namely $s_{\max}=10n$ ultimately obtains a higher NPV than $s_{\max}=5n+10$. 


The performance of MCS with and without local search are compared in the third group with MCS-4 and MCS-7. The convergence speed of MCS without local search is severely decreased and the maximum NPV found is reduced.

\begin{figure}[htb]
  \centering 
  \subfigure[Initial locations]{ 
    \label{fig:subfig:e1_ini_loc} 
    \includegraphics[width=0.45\textwidth]{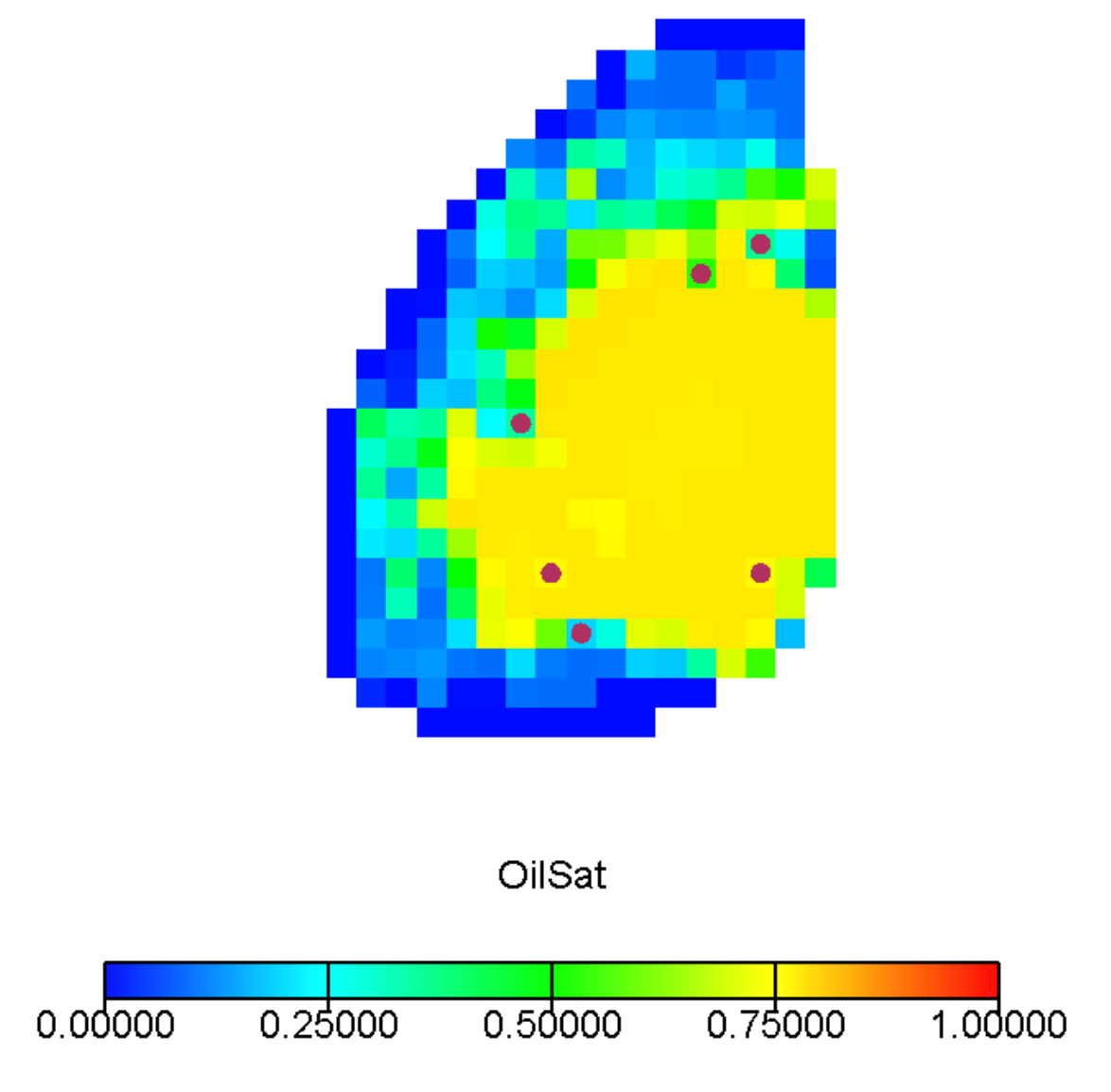}} 
  \subfigure[Optimized locations]{ 
    \label{fig:subfig:e1_op_loc} 
    \includegraphics[width=0.45\textwidth]{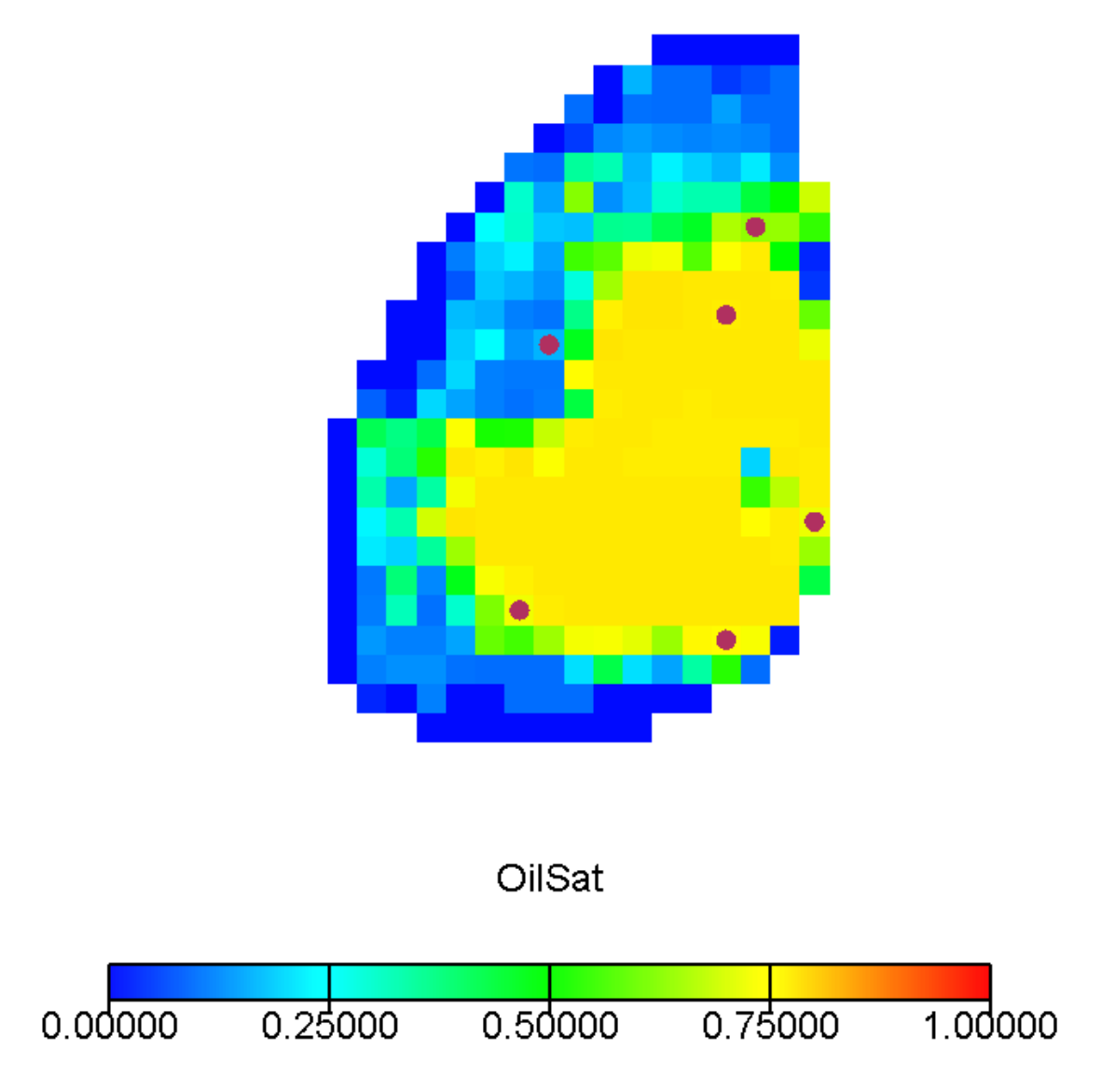}} 
  \caption{The initial well locations and the optimized well locations for Example 1. The base map shows the oil saturation of layer 4 at the end of production. Well locations are indicated by red circles.} 
  \label{fig:e1_loc} 
\end{figure}
\begin{figure}[htb]
\centering
  \includegraphics[width=0.45\textwidth]{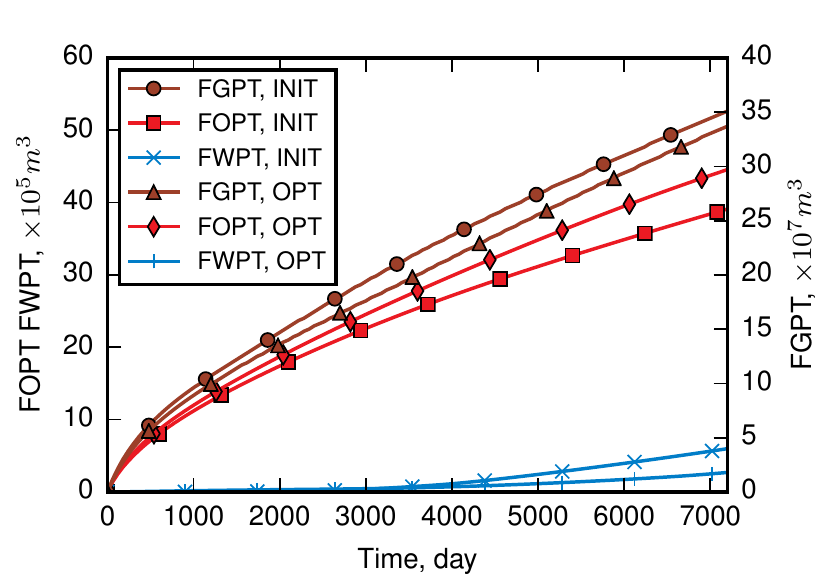}
\caption{Cumulative gas (FGPT), oil (FOPT), and water (FWPT) production with the initial well locations (INIT) and the optimized locations (OPT), versus time, for Example 1.}
\label{fig:e1_fopt}
\end{figure}

The initial well locations and the optimized well locations are shown in Fig.\ \ref{fig:e1_loc}. In this example, the initial locations are chosen from reasonable positions given by industry -- locations used in actual production \citep{gao_quantifying_2006}. As we can see, the wells are drilled around the gas cap. The optimized well locations are still located around the gas cap. 
This is reasonable from a petroleum engineering perspective since the gas cap can keep the pressure up and drive oil to the well bores. Fig.\ \ref{fig:e1_fopt} shows the cumulative gas, oil, and water production using the initial well locations and optimized locations versus time. It is clear that the optimized well locations can produce more oil and less water. The cumulative gas for optimized well locations is lower, this can keep the reservoir pressure higher and drive more oil.

\subsection{Example 2, 3-D angled well placement optimization}
\label{sec:42_awop}

\subsubsection{Reservoir model description}

\begin{figure}[htb]
\centering
  \includegraphics[width=0.45\textwidth]{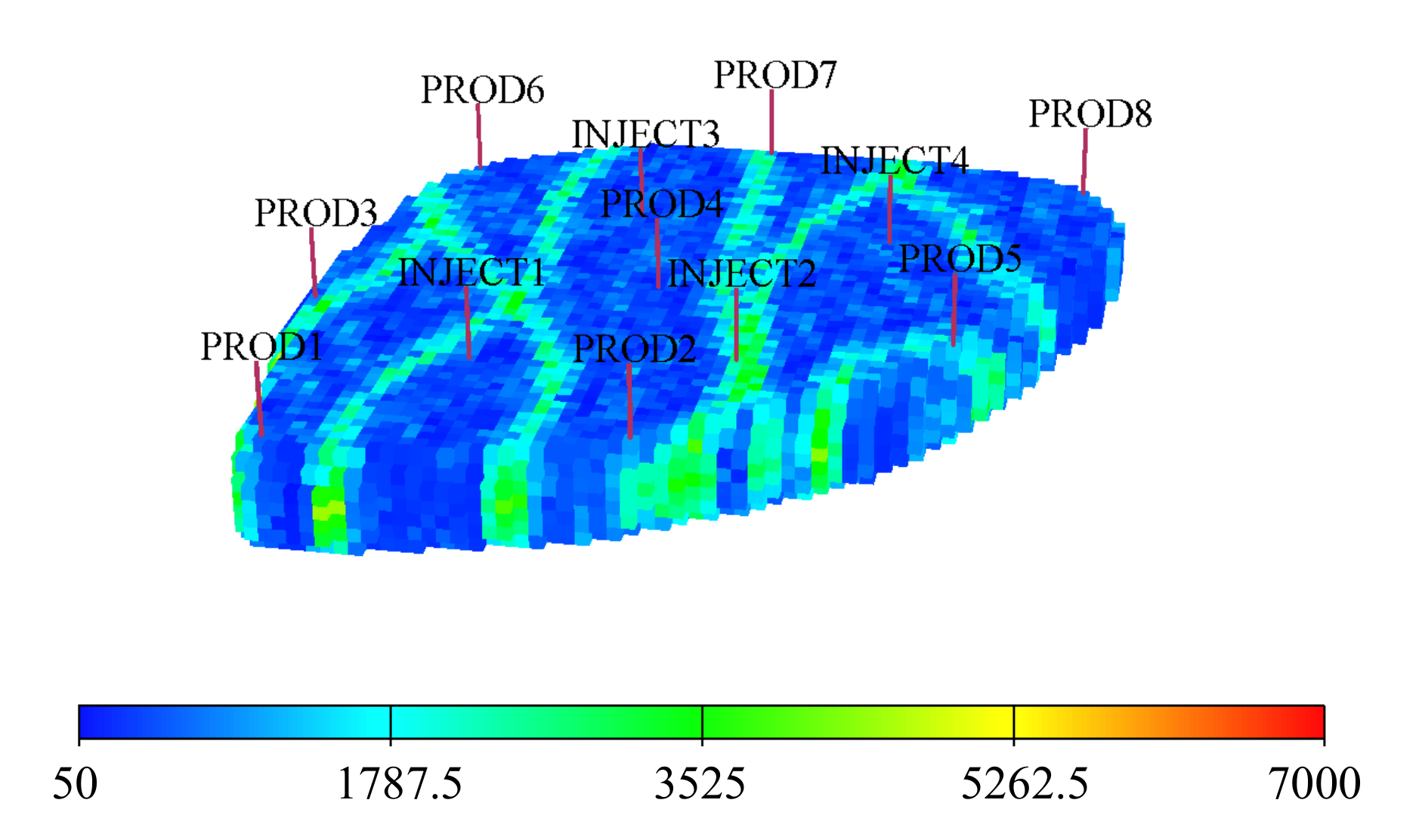}
  \caption{Egg model displaying the permeability and the default placement of wells for Example 2.} 
  \label{fig:eggmodel} 
\end{figure}

This example uses the Egg model which has been used in numerous papers related to well placement and control optimization \citep{zandvliet_adjoint-based_2008,fonseca_robust_2014,siraj_model_2015}. The model uses $60\times 60\times 7=25,200$ grid cells of which 18,553 cells are active. The details of the geological and fluid parameter settings of egg model can be found in \citep{jansen_egg_2014}. Fig.\ \ref{fig:eggmodel} shows the reservoir model, displaying the permeability and the default placement of wells. Note that the model was modified slightly to make it more suitable for production with 
horizontal wells and angled wells. The grid block size is set to $\rm 30 m\times 30 m\times 10m$, and the net to gross thickness ratio is set to 0.2.

We optimize the placement of 12 3-D angled wells (8 producers and 4 injectors) for this example. The total number of variables is 72.
We use the default well placement setting from \citep{jansen_egg_2014} as the initial guess for our optimization problem. That is, the initial $x$ and $y$ are obtained from the default egg model as shown in Fig.\ \ref{fig:eggmodel}. The initial $z$ is set to 1 which means the heel of each well lies in the top layer.  The initial well length, $l$, is set to 60 m. Initially, we choose $\theta=0$ and $\phi=\pi/2$, that is each well is vertical initially. The constraints for this problem include:
\begin{enumerate}
\item $x$ and $y$, the coordinates of each well, are between 1 to 60;
\item $z$, the depth of the well heel, is between 1 to 7;
\item $l$, the length of the well, is between 50 to 300 m;
\item The angle $\theta$ lies between 0 to $2\pi$;
\item The angle $\phi$ of each well lies between 0 and $\pi/2$.
\end{enumerate}

The reservoir simulation time is 20 years. All wells are controlled by bottom hole pressure: 395 bar for production wells and 410 bar for injection wells. The objective function is the NPV and the related economic parameters are given in Table \ref{tab:eco_egg}.

\begin{table}[htb]
\caption{Economic parameters used for Egg model in Example 2.}
\label{tab:eco_egg}
\centering
\begin{tabular}{cc}
\hline\noalign{\smallskip}
Parameter & Value \\
\noalign{\smallskip}\hline\noalign{\smallskip}
Base drill cost (USD/well) & 25M \\
Drilling cost (USD/m) & 50,000 \\
Gas revenue (USD/$\rm m^3$) & 0.5 \\
Oil revenue (USD/$\rm m^3$) & 500.0 \\
Water-production cost (USD/$\rm m^3$) & 80.0 \\
Water-injection cost (USD/$\rm m^3$) & 80.0 \\
Annual discount rate & 0 \\
\noalign{\smallskip}\hline
\end{tabular}
\end{table}


\subsubsection{Results and discussion}


The results of Example 2 are shown in Table \ref{tab:e22_result}. In this table, the ultimate NPV after 10000 simulation runs for each algorithm is given. Moreover, for PSO and CMA-ES, the maximum, minimum, mean, median, and standard deviation of NPV are shown.
From the table we can see that unlike the simple Example 1, PSO eventually outperforms the other algorithms. This is because the search space for this problem is much larger than the simple examples, and PSO has good ability explore the entire search space. MCS outperforms the other algorithms with a limited computational budget and without the variability inherent in PSO.


\begin{table}[htb]
\caption{Results for Example 2. Values shown are NPV in $\$ \times 10^8$ USD. }
\label{tab:e22_result}
\centering
\subtable[Deterministic algorithms (MCS, GPS)]{
\label{tab:e22_result_det}
\begin{tabular}{p{2cm}|p{8cm}}
\hline
Algorithm & \hfil NPV \\
\hline
MCS-1  & \hfil 3.89 \\
MCS-2  & \hfil 8.97 \\
MCS-3  & \hfil 7.72 \\
MCS-4  & \hfil 8.47 \\
MCS-5  & \hfil 8.70 \\
MCS-6  & \hfil 8.53 \\
MCS-7  & \hfil 7.79 \\
GPS    & \hfil 8.19 \\
\hline
\end{tabular}
}

\subtable[Stochastic algorithms (PSO, CMA-ES)]{        
\label{tab:e22_result_sto}
\begin{tabular}{p{2cm}|p{1cm}|p{1cm}p{1cm}p{1cm}p{1cm}p{1cm}}
\hline
Algorithm & Trials & Max & Min & Mean & Median & Std. \\
\hline
PSO    &5 &11.35&8.37&9.07&8.71&1.14\\
CMA-ES &5 &8.86&8.64&8.72&8.69&0.09\\
\hline
\end{tabular}
}
\end{table}

\begin{figure}[htb]
\centering
  \includegraphics[width=0.45\textwidth]{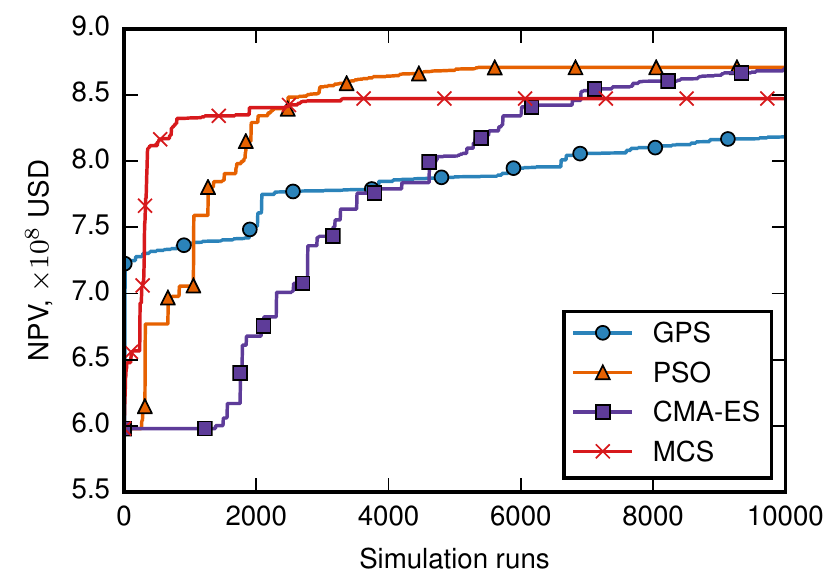}
  \caption{Optimization performance for Example 2. For PSO and CMA-ES, the solid lines depict the median NPV. MCS here is a label for MCS-4.} 
  \label{fig:uncwop_algo} 
\end{figure}

Plots of the NPV for the four algorithms versus the number of simulation runs are shown in Fig.\ \ref{fig:uncwop_algo}. Note that for PSO and CMA-ES, 5 trials are performed and the solid lines depict the median NPV over all 5 trials.
From the figure we can see that MCS showed excellent convergence speed at the early stage of optimization (simulation runs $<$2000). This is useful for optimization given a limited computational budget.

\begin{figure}[htb]
\centering
  \includegraphics[width=0.45\textwidth]{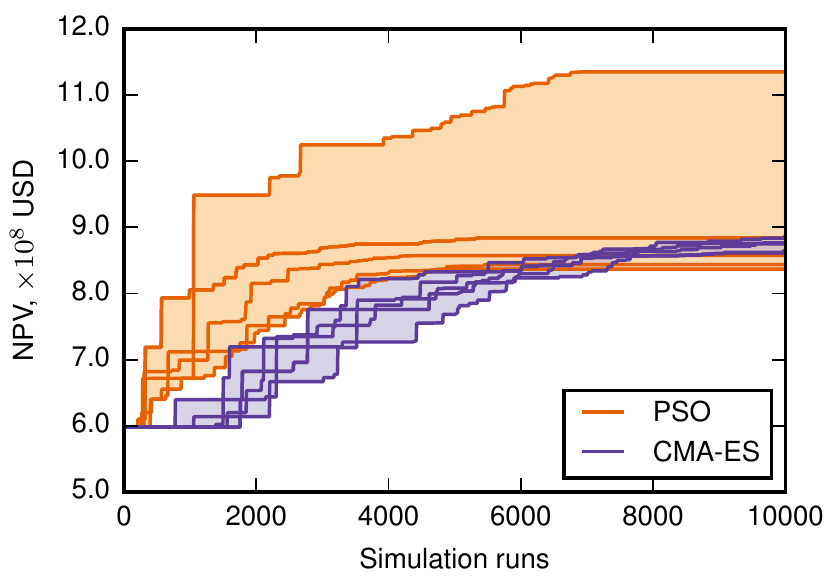}
  \caption{The range of NPV found amongst the trials of PSO and CMA-ES for Example 2. Each solid line represents a trial.} 
  \label{fig:uncwop_psocmaes} 
\end{figure}

Fig.\ \ref{fig:uncwop_psocmaes} shows the range of NPV found amongst the trials of PSO and CMA-ES. For PSO and CMA-ES, most trials converge to a NPV between 8 to 9$\times 10^8$USD. PSO is able to 
find a much higher NPV (about 11$\times 10^8$USD). This shows that the problem is hard and most algorithms have only converged to a local optima. It also shows that PSO potentially has a better ability for this type of hard problem. 

\begin{figure}[htb]
  \centering 
  \subfigure[MCS-1]{ 
    \label{fig:subfig:e22_mcs1_initanalysis} 
    \includegraphics[width=0.45\textwidth]{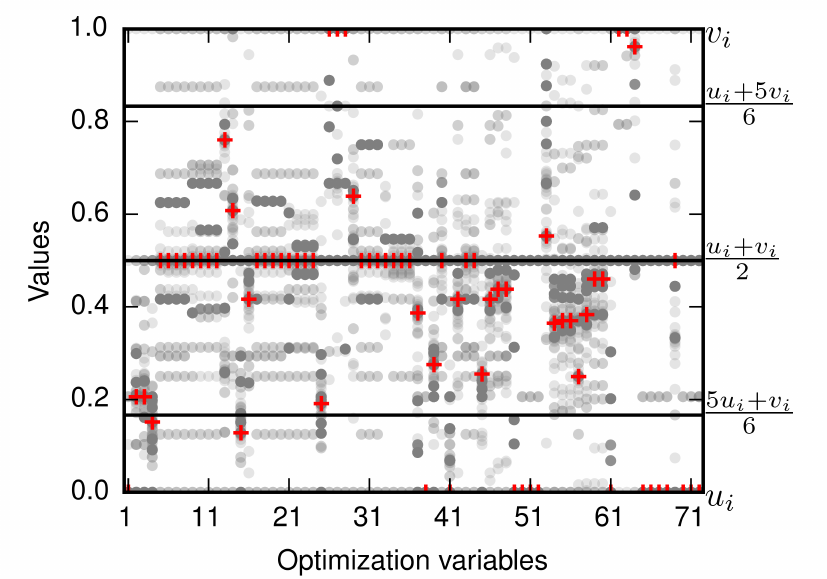}} 
  \subfigure[MCS-2]{ 
    \label{fig:subfig:e22_mcs2_initanalysis} 
    \includegraphics[width=0.45\textwidth]{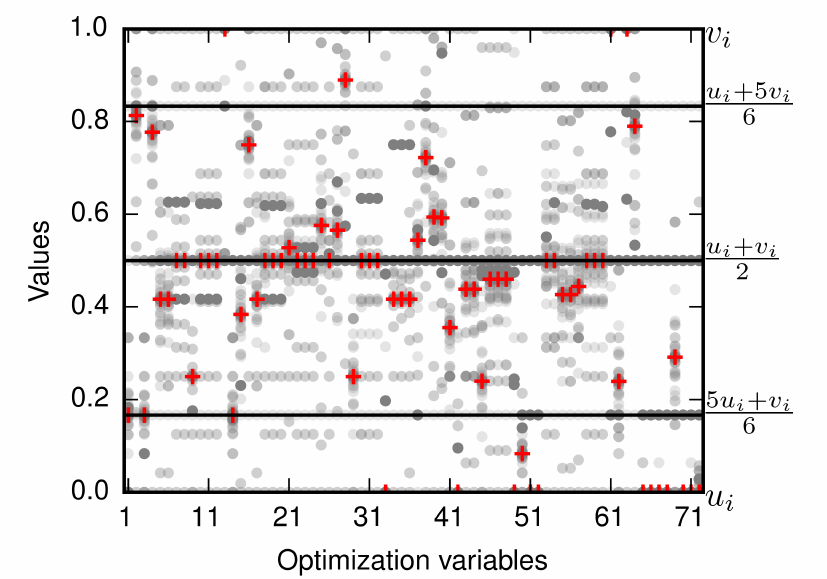}}   
  \caption{Normalized boundary, initialization lists, all optimization candidates, and the optimal solution for Example
2. Circles denote optimization candidates. Cross signs denote the optimal solutions.} 
  \label{fig:e22_mcs_initanalysis} 
\end{figure}

\begin{figure}[htb]
  \centering 
  \subfigure[Initialization list]{ 
    \label{fig:subfig:e22_mcs_init} 
    \includegraphics[width=0.45\textwidth]{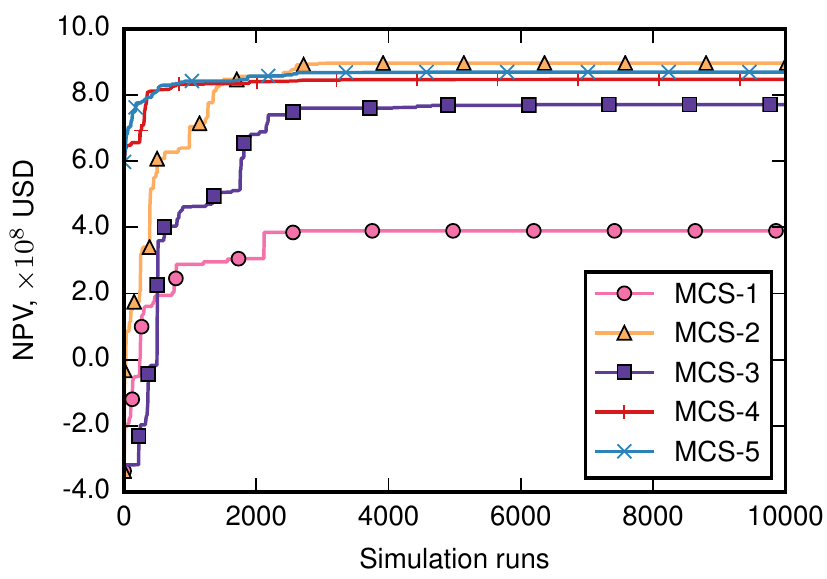}} 
  \subfigure[Levels]{ 
    \label{fig:subfig:e22_mcs_smax} 
    \includegraphics[width=0.45\textwidth]{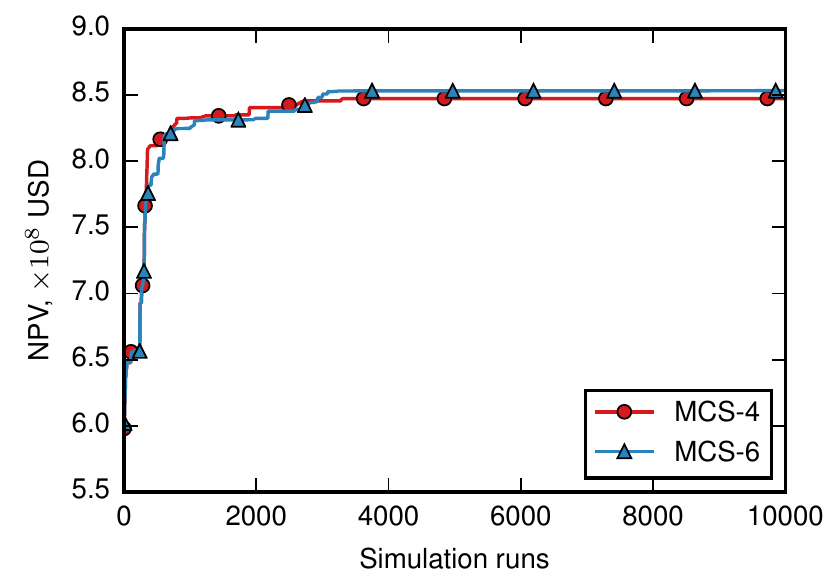}} 
  \subfigure[Local search]{ 
    \label{fig:subfig:e22_mcs_ls} 
    \includegraphics[width=0.45\textwidth]{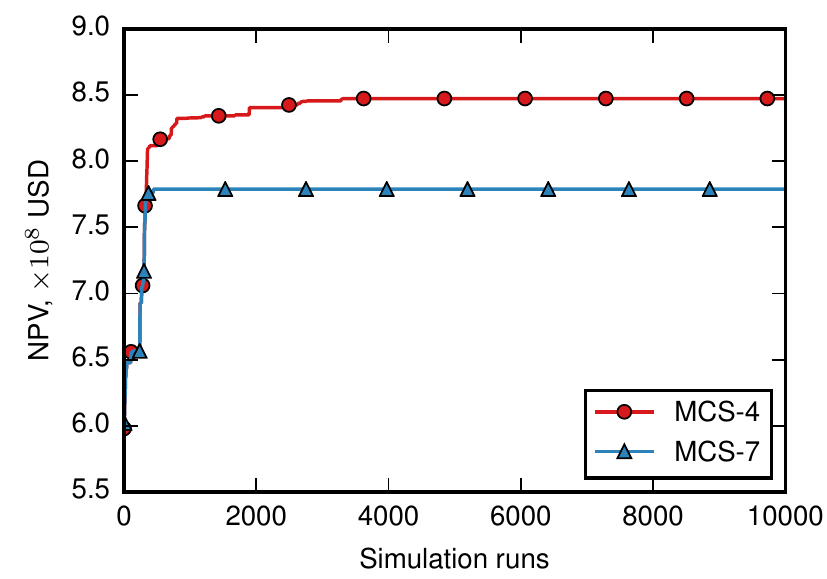}}   
  \caption{The performance of different configurations of MCS for Example 2.} 
  \label{fig:e22_mcs} 
\end{figure}

Fig.\ \ref{fig:e22_mcs} shows the performance of the different configurations of MCS. Using the default initialization list, as in MCS-1 and MCS-2 gives an unreasonable initial guess. MCS-4 and MCS-5 use a reasonable initial guess leading to a much improved performance. Despite this MCS-2 eventually obtains the highest NPV.
In Fig.\ \ref{fig:e22_mcs_initanalysis} we normalize the search space to the [0,1]-interval and map
the all optimization candidates (denoted by circles) and the optimal solutions (denoted by crosses) of MCS-1 and MCS-2 to the normalized search space. From the figure we can see that 46 out of 72 optimization variables in the optimal solution obtained by MCS-1 after 10000 function evaluations are still located in the positions defined by the initialization list. This occurs for 23 out of 72 variables for MCS-2. 
For this problem, an uniformly spaced initialization list, not containing any boundary points, is more suitable than an initialization list with boundary points.

MCS-3 uses line search to generate the initialization list automatically. The optimization starts with a very low NPV, but eventually obtains a high NPV. This configuration is highly recommended for problems where the users can not provide a good initial guess.

MCS-6, using a larger maximum number of levels ($s_{\max}=10n$), performs better than MCS-4 with $s_{\max}=5n+10$. For a large scale optimization problem, a larger maximum number of levels is recommended. MCS-7, which turns off the local search phase, performs a little worse than MCS-4.

\begin{figure}[htb]
  \centering 
  \subfigure[Initial well placement]{ 
    \label{fig:subfig:e22_well_init} 
    \includegraphics[width=0.45\textwidth]{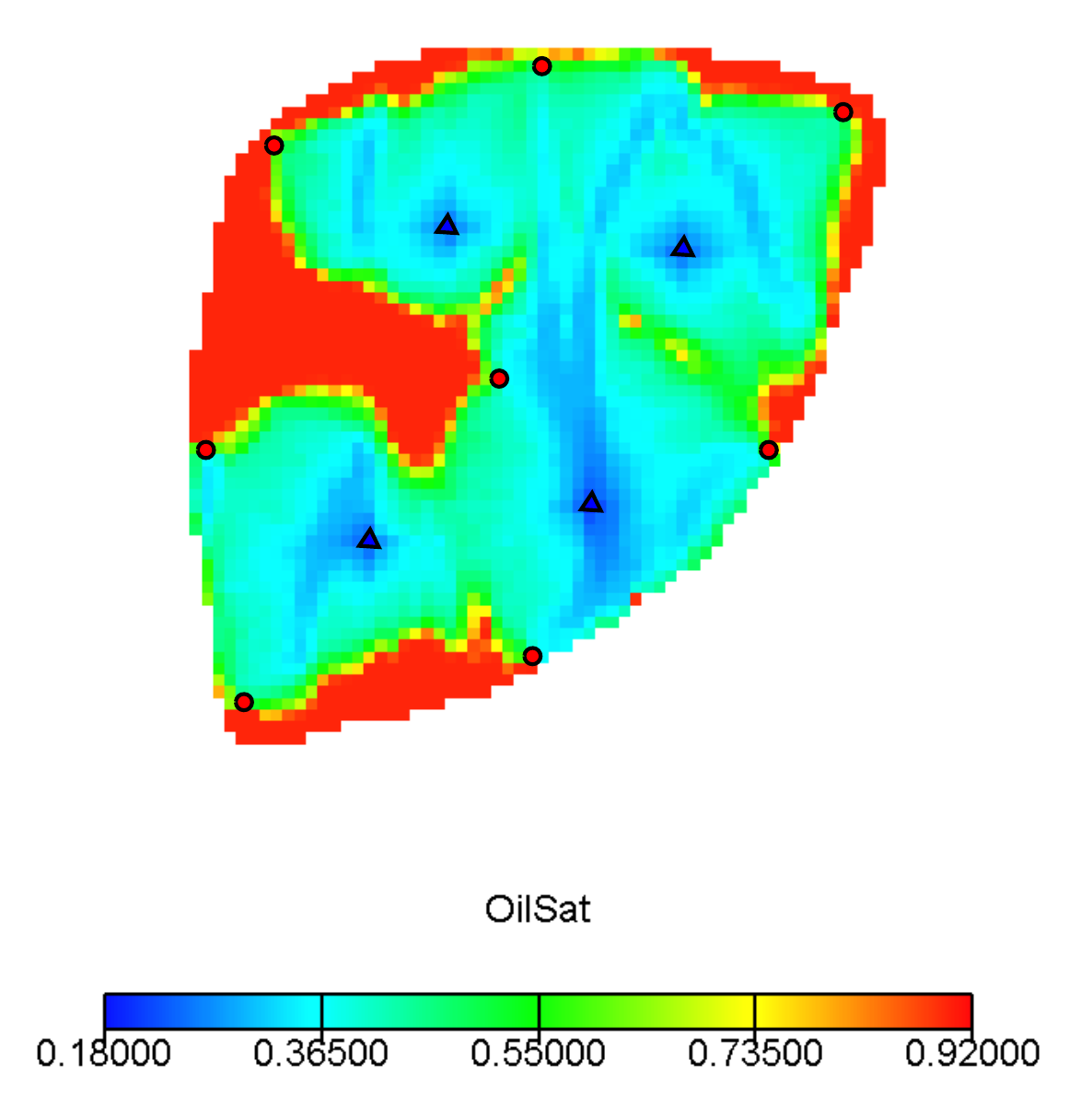}} 
  \subfigure[Optimal well placement]{ 
    \label{fig:subfig:e22_well_op} 
    \includegraphics[width=0.45\textwidth]{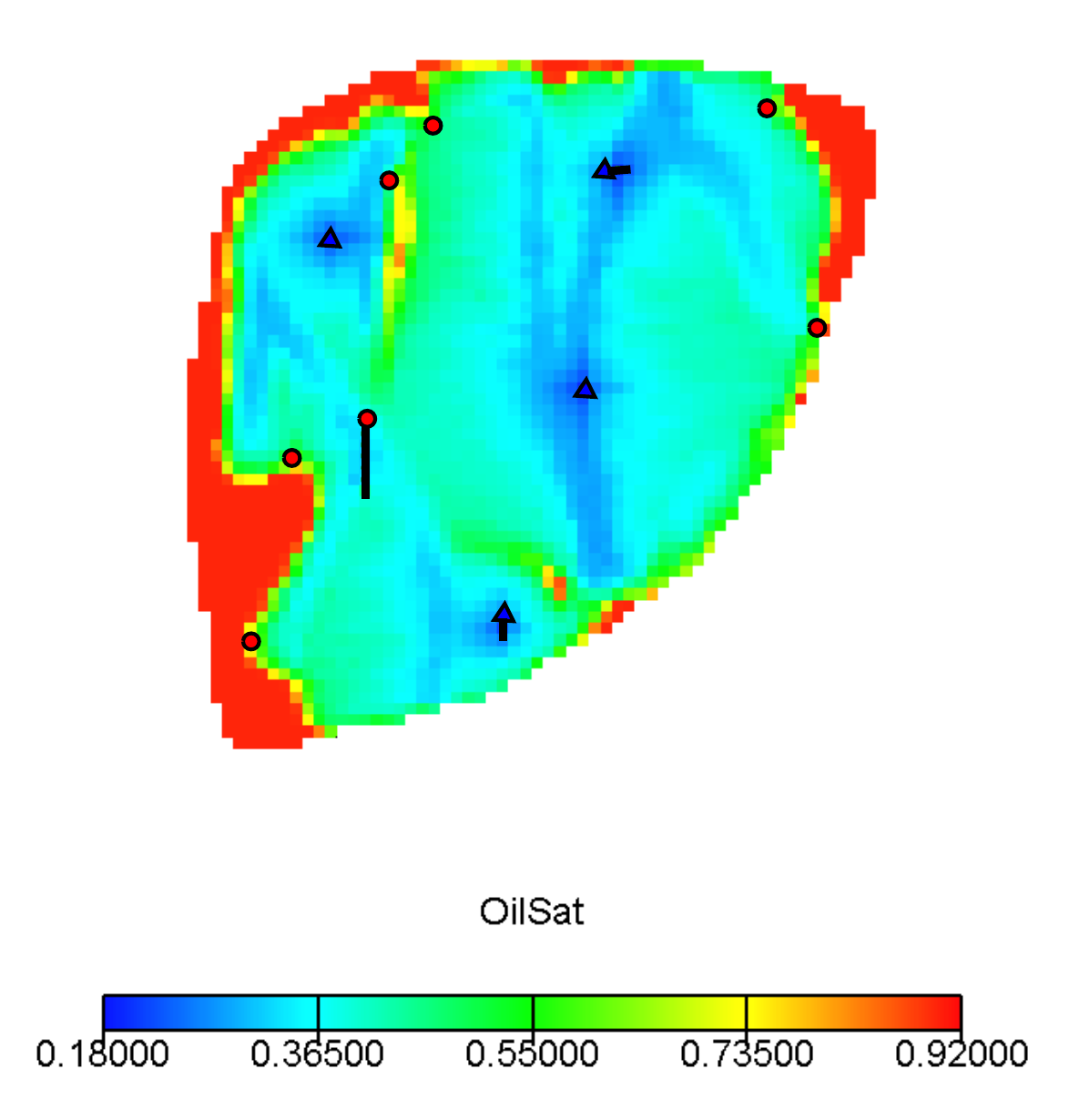}}  
  \caption{ The well placement and the final oil saturation distribution for Example 2.} 
  \label{fig:e22_well} 
\end{figure}

The final oil saturation distribution for the initial well placement and the optimal well placement are shown in Fig.\ \ref{fig:e22_well}. All wells are vertical for the initial well placement. For the optimal well placement obtained by the optimization algorithm, some wells are angled.  From the final oil saturation distribution, we can see that the well placement obtained by
the optimization algorithm gives a larger sweep area, and obtains higher production performance eventually.

\subsection{Example 3, well control optimization}
\label{sec:43_wcop}

\subsubsection{Reservoir model description}

\begin{figure}[htb]
\centering
  \includegraphics[width=0.5\textwidth]{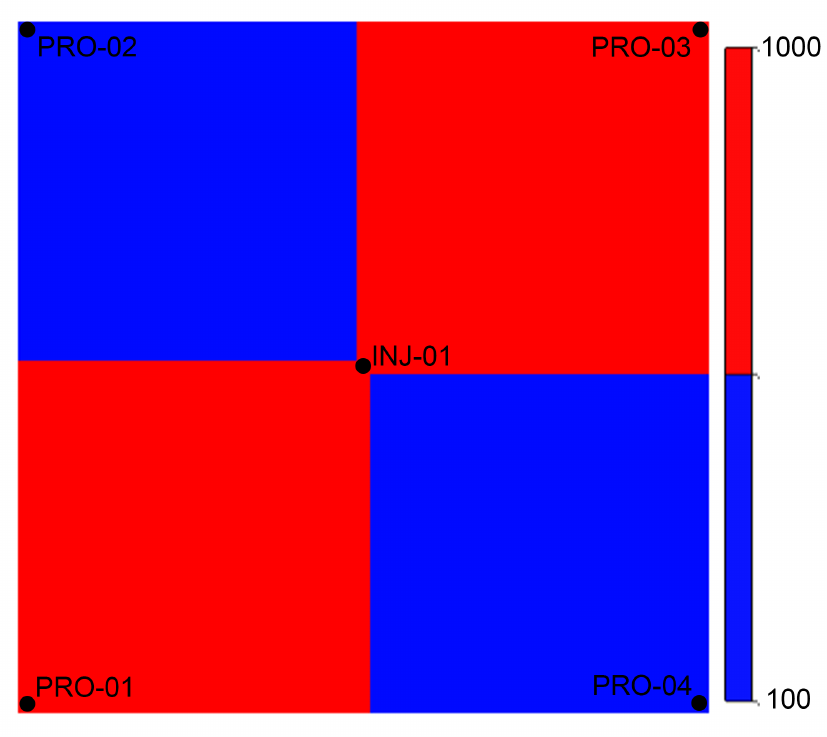}
\caption{Permeanbility field (mD) for the five-spot model used in Example 3.}
\label{fig:e1_perm}
\end{figure}

This example from \citep{oliveira_adaptive_2014} uses a single-layer reservoir model with $51\times 51$ uniform grid blocks with $\Delta x=\Delta y=10$m and $\Delta z=5$m. The model consists of four production wells and one injection well. The wells form a five-spot well pattern. We consider an oil-water two phase flow in this model. The permeability field and well placements are shown in Fig.\ \ref{fig:e1_perm}. There are two high permeability zones and two low permeability zones in the model. The permeabilities are $\rm 1000 mD$ and $\rm 100 mD$, respectively. Detailed reservoir information is given in Table \ref{tab:e2_Rpara}.

\begin{table}[htb]
\caption{Reservoir parameters used in Example 3.}
\label{tab:e2_Rpara}
\centering
\begin{tabular}{cc}
\hline\noalign{\smallskip}
Parameter & Value \\
\noalign{\smallskip}\hline\noalign{\smallskip}
Reservoir grid & $51\times 51\times 1$ \\
Grid size (m) & $\rm 10\times 10\times 5$ \\
Porosity & 0.2 \\
Net-to-gross ratio & 0.2 \\
Initial oil saturation & 0.8 \\
Initial pressure ($\rm bar$) & 200 \\
Oil viscosity ($\rm mPa\cdot s$) & 0.42 \\
Water viscosity ($\rm mPa\cdot s$) & 1.7 \\
\noalign{\smallskip}\hline
\end{tabular}
\end{table}

The reservoir lifetime is set to 720 days. With a fixed injection rate of $\rm 240m^3/d$ for well INJ-01, we seek to optimize the liquid rates of four production wells. Two variations of this optimization problem are considered.
\begin{itemize}
 \item Case 1: each well is produced under a liquid rate throughout its lifetime. This gives 4 optimization variables in total.
 \item Case 2: the liquid rate for each well is updated every 90 days (8 control periods). This gives 32 optimization variables in total. 
 \end{itemize} 

The objective function we use for this example is NPV and the corresponding economic parameters are the same as in Example 1 and are given in Table \ref{tab:eco_set_1}. Only bound constraints are considered and the detailed optimization parameters are given in Table \ref{tab:e2_Opara}.

\begin{table}[htb]
\caption{Optimization parameters used in Example 3.}
\label{tab:e2_Opara}
\centering
\begin{tabular}{ccc}
\hline\noalign{\smallskip}
Parameter & Case 1 & Case 2 \\
\noalign{\smallskip}\hline\noalign{\smallskip}
Variables & 4 & 32 \\
Initial rate (PRO-01, PRO-03) ($\rm m^3/d$) & 20 & 20 \\
Initial rate (PRO-02, PRO-04) ($\rm m^3/d$) & 40 & 40 \\
Minimum rate ($\rm m^3/d$) & 0 & 0 \\
Maximum rate (PRO-01, PRO-03) ($\rm m^3/d$) & 40 & 40 \\
Maximum rate (PRO-02, PRO-04) ($\rm m^3/d$) & 80 & 80 \\
\noalign{\smallskip}\hline
\end{tabular}
\end{table}

\subsubsection{Results and discussion}

The results for the two cases of Example 3 are shown in Table \ref{tab:e20_result} and Table \ref{tab:e23_result}. Case 1 uses a maximum of 400 simulations to optimize 4 variables and Case 2 uses a maximum of 3200 simulations to optimize 32 variables. 

The initial guess for the optimal control strategy for GPS, PSO, and CMA-ES is the point half way between the lower and upper bounds. Coincidentally, MCS with its default settings also uses the middle value as its start point. So for this example, MCS-1 and MCS-4 are identical, and MCS-2 and MCS-5 are identical. Hence, we omit MCS-4 and MCS-5 while analyzing the results.

From Table \ref{tab:e20_result} we can see that for Case 1, all algorithms GPS, PSO, CMA-ES and MCS (except for configuration MCS-7) are able to obtain a high NPV value at the end of the optimization, and the ultimate difference between the algorithms is small. The mean and median NPV found by PSO is slightly smaller than the other algorithms. For Case 2, similar conclusions can be drawn from Table \ref{tab:e23_result}. After 3200 simulation runs, GPS obtains the highest NPV. CMA-ES and MCS (again except for configuration MCS-7) are in the middle, while PSO performs the worst.

\begin{table}[htb]
\caption{Results for Case 1 of Example 3. Values shown are NPV in $\$ \times 10^6$ USD. }
\label{tab:e20_result}
\centering
\subtable[Deterministic algorithms (MCS, GPS)]{
\label{tab:e20_result_det}
\begin{tabular}{p{2cm}|p{8cm}}
\hline
Algorithm & \hfil NPV \\
\hline
MCS-1  & \hfil 5.29 \\
MCS-2  & \hfil 5.30 \\
MCS-3  & \hfil 5.30 \\
MCS-6  & \hfil 5.28 \\
MCS-7  & \hfil 4.85 \\
GPS    & \hfil 5.31 \\
\hline
\end{tabular}
}

\subtable[Stochastic algorithms (PSO, CMA-ES)]{        
\label{tab:e20_result_sto}
\begin{tabular}{p{2cm}|p{1cm}|p{1cm}p{1cm}p{1cm}p{1cm}p{1cm}}
\hline
Algorithm & Trials & Max & Min & Mean & Median & Std. \\
\hline
PSO    &10 &5.27&5.12&5.22&5.23&0.04\\
CMA-ES &10 &5.31&5.30&5.30&5.30&0.00\\
\hline
\end{tabular}
}
\end{table}

\begin{table}[htb]
\caption{Results for Case 2 of Example 3. Values shown are NPV, $\times 10^6$ USD. }
\label{tab:e23_result}
\centering
\subtable[Deterministic algorithms (MCS, GPS)]{
\label{tab:e23_result_det}
\begin{tabular}{p{2cm}|p{8cm}}
\hline
Algorithm & \hfil NPV \\
\hline
MCS-1  & \hfil 11.99 \\
MCS-2  & \hfil 12.22 \\
MCS-3  & \hfil 12.19 \\
MCS-6  & \hfil 11.67 \\
MCS-7  & \hfil 10.37 \\
GPS    & \hfil 12.35 \\
\hline
\end{tabular}
}

\subtable[Stochastic algorithms (PSO, CMA-ES)]{        
\label{tab:e23_result_sto}
\begin{tabular}{p{2cm}|p{1cm}|p{1cm}p{1cm}p{1cm}p{1cm}p{1cm}}
\hline
Algorithm & Trials & Max & Min & Mean & Median & Std. \\
\hline
PSO    &10 &12.24&11.03&11.91&11.97&0.37\\
CMA-ES &10 &12.35&12.27&12.34&12.35&0.02\\
\hline
\end{tabular}
}
\end{table}

\begin{figure}[htb]
 \centering 
  \subfigure[Case 1]{ 
    \label{fig:subfig:e2_algo0} 
    \includegraphics[width=0.45\textwidth]{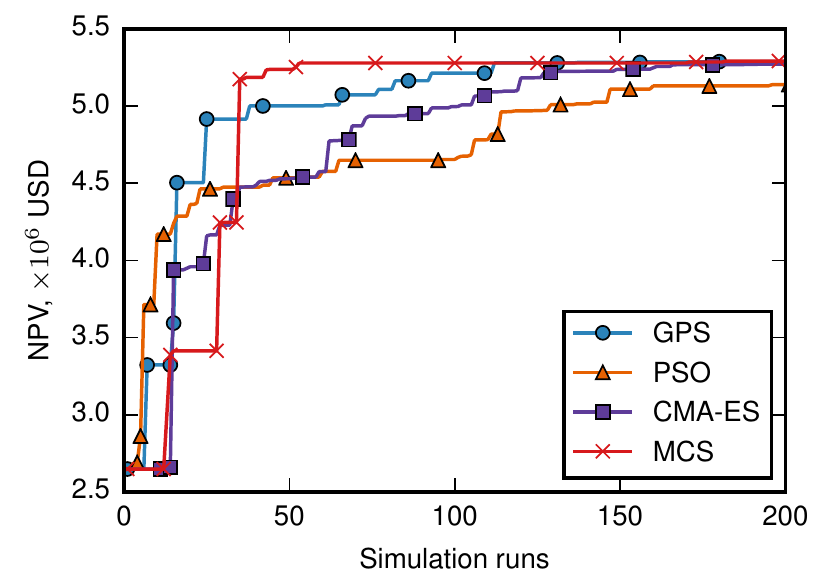}} 
  \subfigure[Case 2]{ 
    \label{fig:subfig:e2_algo3} 
    \includegraphics[width=0.45\textwidth]{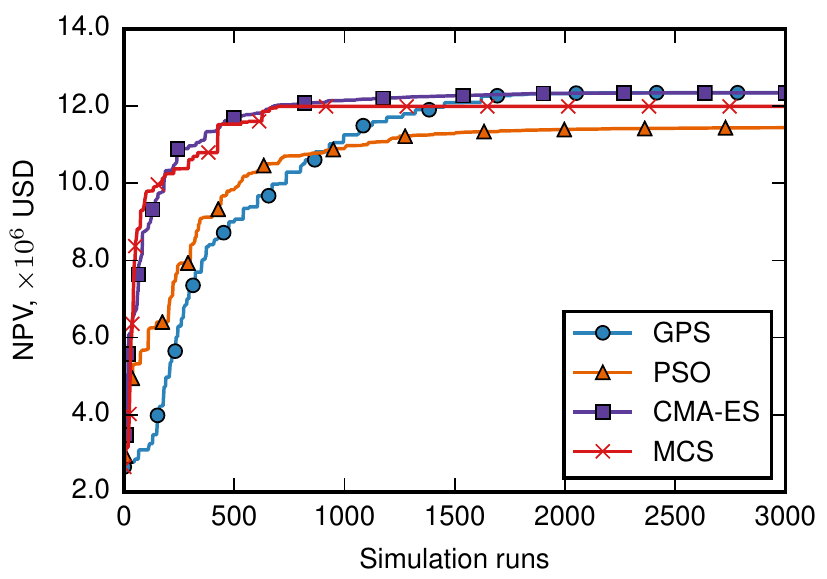}} 
\caption{Optimization performance for Example 3. For PSO and CMA-ES, the solid lines depict the median NPV. MCS here is a label for MCS-4.}
\label{fig:e2_algo}
\end{figure}

Plots of the NPV of the four algorithms versus the number of simulation runs are shown in Fig.\ \ref{fig:e2_algo}. As in Example 1, 10 trials are performed for PSO and CMA-ES, and the solid lines depict the median NPV over all 10 trials of these two algorithms.

From Fig.\ \ref{fig:e2_algo} we can see that, the final NPV obtained by MCS is not the highest over all algorithms tested, however, once again MCS outperforms when the number of simulation runs is limited and is not affected by the variability of the other algorithms. When the number of simulation runs is limited to 15\% of the final number of simulation runs (60 simulation runs for Case 1 and 480 simulation runs for Case 2), the NPV obtained by each algorithm is given in Table \ref{tab:e2_result_limited}. Note that in this table we use the median NPV of 10 trials for PSO and CMA-ES. We use the median instead of the mean because it is less sensitive to outliers in the data. When the total number of simulation runs is limited, MCS showed significant advantages over PSO, GPS, and CMA-ES. Again MCS-7 provides poor results -- showing the importance of the local search feature within MCS. This table shows the potential of MCS with a low computational budget.

\begin{table}[htb]
\caption{Results of Example 3 with a limited number of simulation runs. Values shown are NPV, $\times 10^6$ USD obtained after 15\% of the maximum simulation runs.}
\label{tab:e2_result_limited}
\centering
\subtable[Case 1]{
\begin{tabular}{lc}
\hline
Algorithm & NPV \\
\hline
MCS-1  &5.28\\
MCS-2  &5.28\\
MCS-3  &5.26\\
MCS-6  &5.28\\
MCS-7  &3.88\\
GPS    &5.01\\
PSO    &4.79\\
CMA-ES &4.92\\
\hline
\end{tabular}
}
\qquad
\subtable[Case 2]{        
\begin{tabular}{lc}
\hline
Algorithm & NPV \\
\hline
MCS-1  &11.53\\
MCS-2  &11.83\\
MCS-3  &11.96\\
MCS-6  &10.63\\
MCS-7  &9.80\\
GPS    &8.80\\
PSO    &9.67\\
CMA-ES &11.52\\
\hline
\end{tabular}
}
\end{table}

As we progress from Case 1 to Case 2, the number of optimization variables increases from 4 to 32. The performance of GPS with a low number of simulation runs decreases. In Case 2, the maximum NPV found by GPS is less than the other 3 algorithms when the number of simulation runs is limited to 1000. After 1000 simulation runs, GPS is able to find a higher NPV than PSO. The early stage of the optimization process mainly reflects the global search phase, and the later stage of the optimization process includes the effect of the local search phase for the algorithms tested. In Case 2, it is clear that PSO performs better than GPS at an early stage, but GPS outperforms later. Overall, MCS, which includes both a global search phase and a local search phase, showed a better convergence rate than GPS and PSO. 

\begin{figure}[htb]
  \centering 
  \subfigure[Case 1, PSO]{ 
    \label{fig:subfig:e20_pso} 
    \includegraphics[width=0.45\textwidth]{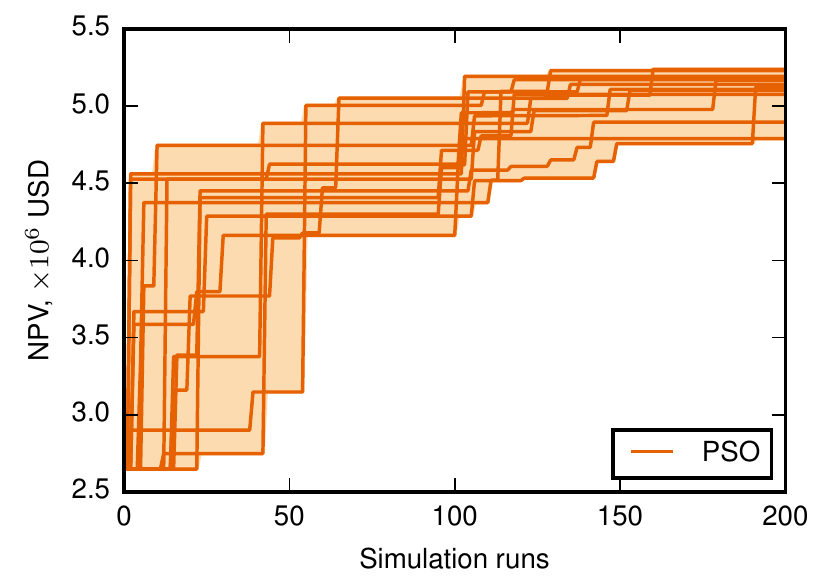}} 
  \subfigure[Case 1, CMA-ES]{ 
    \label{fig:subfig:e20_cmaes} 
    \includegraphics[width=0.45\textwidth]{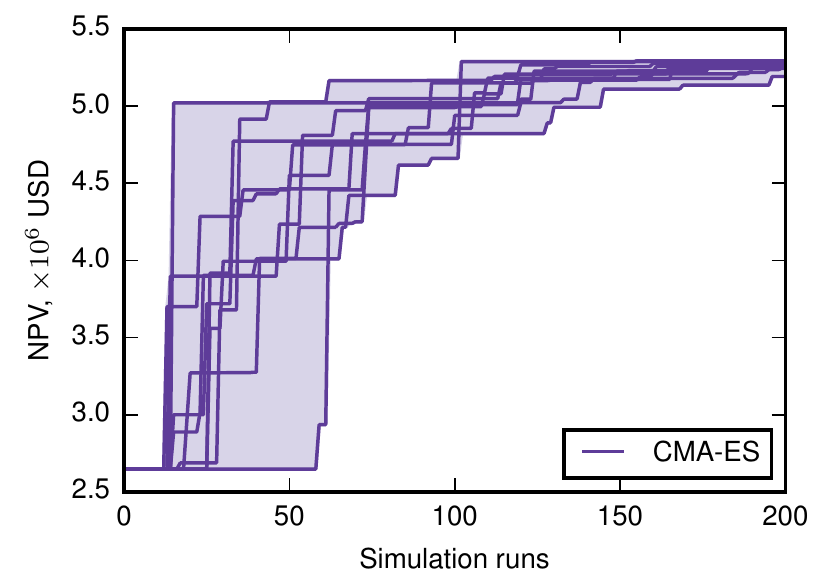}} 

  \subfigure[Case 2, PSO]{ 
    \label{fig:subfig:e23_pso} 
    \includegraphics[width=0.45\textwidth]{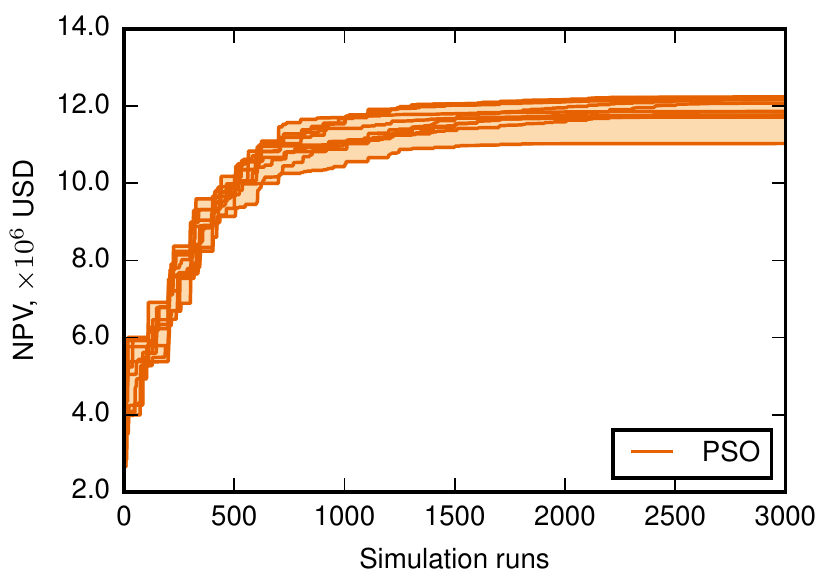}} 
  \subfigure[Case 2, CMA-ES]{ 
    \label{fig:subfig:e23_cmaes} 
    \includegraphics[width=0.45\textwidth]{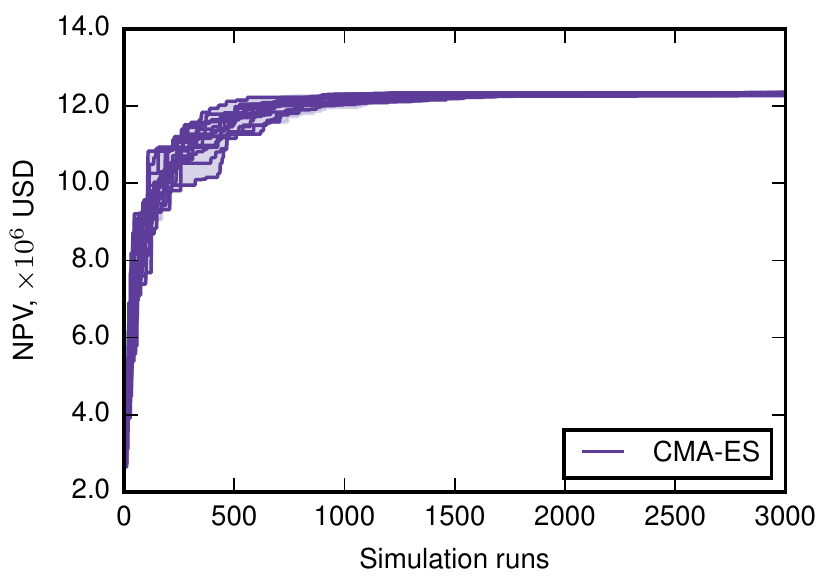}} 
  \caption{The range of NPV found amongst the trials of PSO and CMA-ES for Example 3. Each solid line represents a trial.} 
  \label{fig:e2_psocmaes} 
\end{figure}

Fig.\ \ref{fig:e2_psocmaes} shows the range of NPV for the trials for PSO and CMA-ES for Example 3. In this figure, the areas between the maximum and minimum NPV are shaded for PSO and CMA-ES. From this figure we can see that for Case 1, the range of the best NPV is large initially and then the range decreases for both PSO and CMA-ES. CMA-ES has a small variability near convergence. For Case 2, with a larger number of optimization variables than Case 1, the range of NPV does not decrease for PSO. Each trial falls into a local optima and has a difficult time to escape. The range for CMA-ES decreases to near zero. This indicates that for PSO and CMA-ES, a large computational budget can decrease the performance variability for this example. Compared to CMA-ES, PSO more easily falls into a local optima for problems with a large number of optimization variables.

\begin{figure}[htb]
  \centering 
  \subfigure[Initialization list]{ 
    \label{fig:subfig:e2_0_mcs_init} 
    \includegraphics[width=0.45\textwidth]{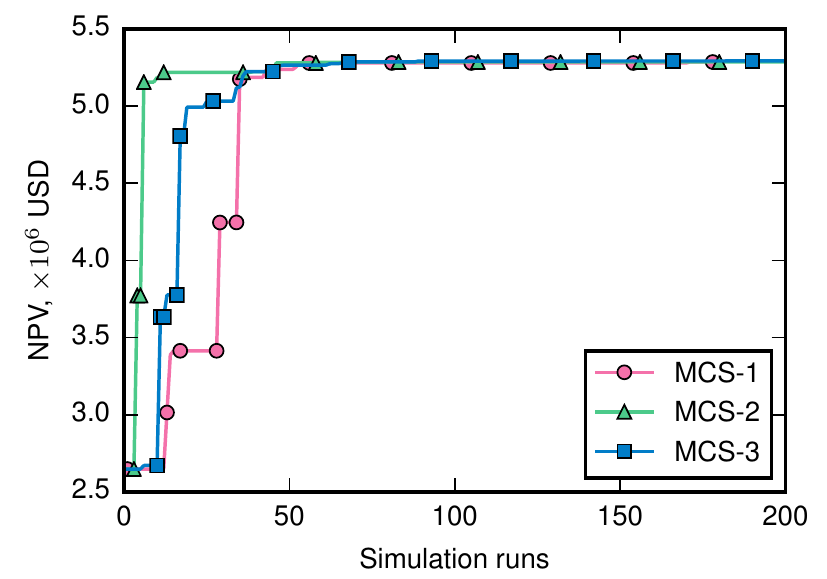}} 
  \subfigure[Levels]{ 
    \label{fig:subfig:e2_0_mcs_smax} 
    \includegraphics[width=0.45\textwidth]{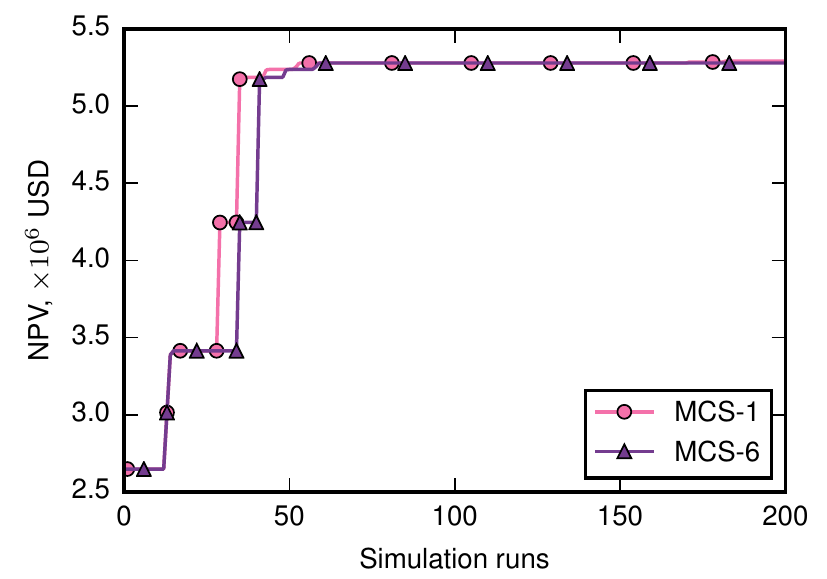}} 
  \subfigure[Local search]{ 
    \label{fig:subfig:e2_0_mcs_ls} 
    \includegraphics[width=0.45\textwidth]{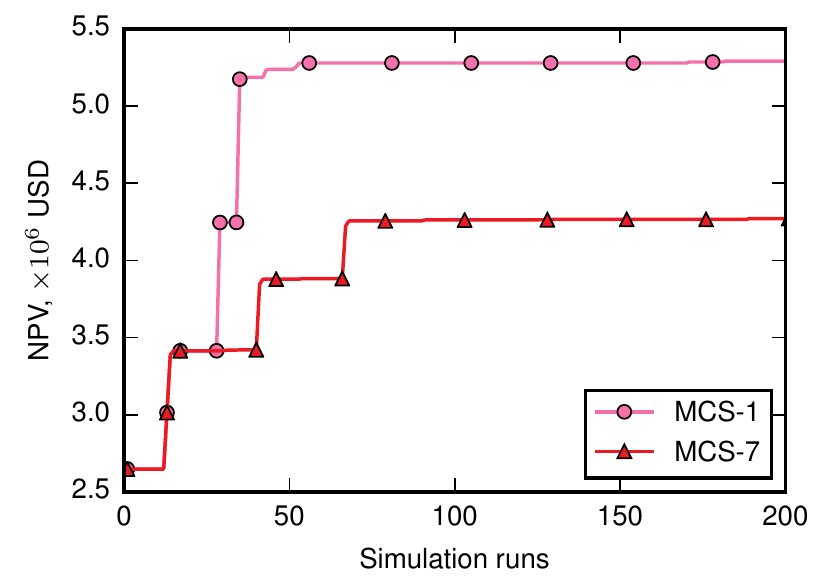}}   
  \caption{The performance of different configurations of MCS for Case 1 of Example 3.} 
  \label{fig:e2_0_mcs} 
\end{figure}

\begin{figure}[htb]
  \centering 
  \subfigure[Initialization list]{ 
    \label{fig:subfig:e2_3_mcs_init} 
    \includegraphics[width=0.45\textwidth]{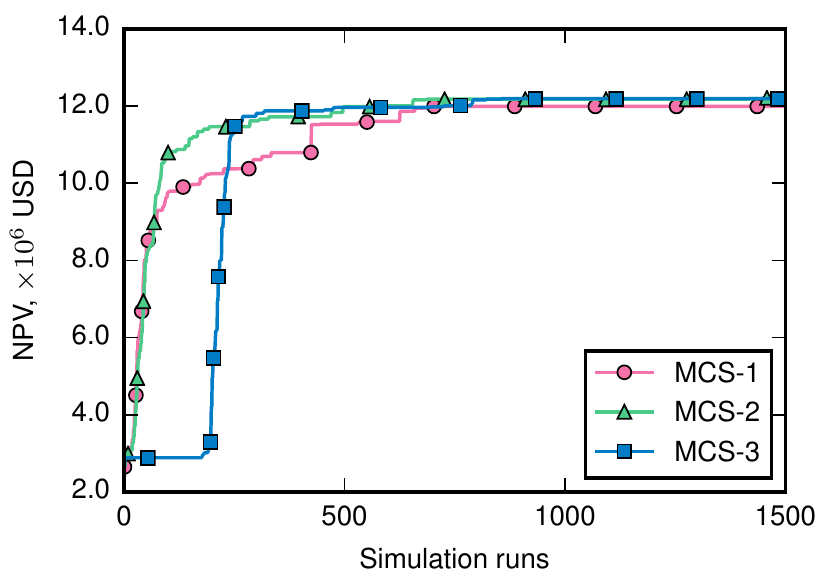}} 
  \subfigure[Levels]{ 
    \label{fig:subfig:e2_3_mcs_smax} 
    \includegraphics[width=0.45\textwidth]{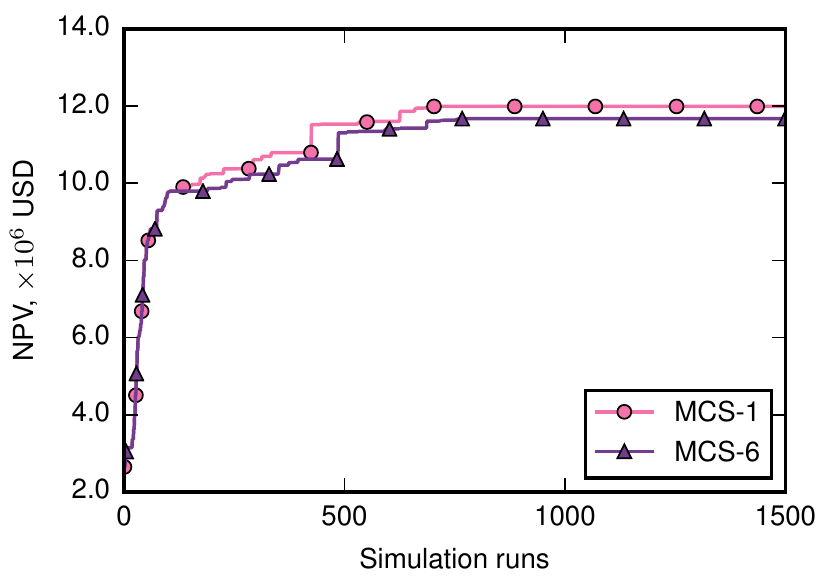}} 
  \subfigure[Local search]{ 
    \label{fig:subfig:e2_3_mcs_ls} 
    \includegraphics[width=0.45\textwidth]{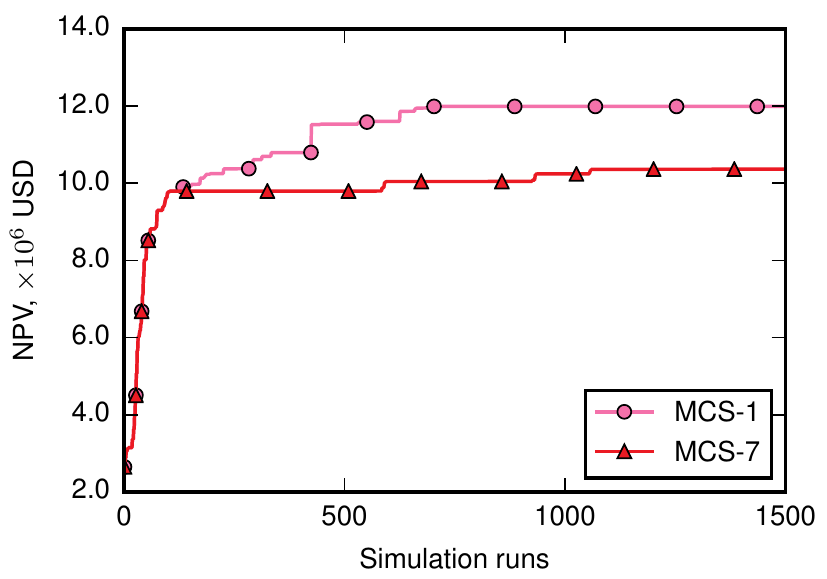}}   
  \caption{The performance of different configurations of MCS for Case 2 of Example 3.} 
  \label{fig:e2_3_mcs} 
\end{figure}

As in Example 1, we tested different MCS configurations and divided them into 3 groups to do further analysis. The results are shown in Fig.\ \ref{fig:e2_0_mcs} and Fig.\ \ref{fig:e2_3_mcs}.


Fig.\ \ref{fig:subfig:e2_0_mcs_init} and Fig.\ \ref{fig:subfig:e2_3_mcs_init} compare the performance of MCS with different initialization lists for the two cases of Example 3. For Case 1, the convergence rate of MCS-2 is the fastest, followed by MCS-3 and MCS-1. MCS-2 and MCS-3 ultimately obtain the highest NPV.

For Case 2, MCS-1 and MCS-2 give a similar convergence rate at an early stage in the optimization process, then MCS-1 falls behind MCS-2. 
MCS-3 shows a very slow rate of convergence at an early stage of the optimization process, but it obtains the highest NPV finally. MCS-3 generates the initialization list by using a line search. This takes a few additional simulation runs before the splitting and local search steps. This explains the slow convergence initially.

The effect of the maximum number of levels is shown in Fig.\ \ref{fig:subfig:e2_0_mcs_smax} and Fig.\ \ref{fig:subfig:e2_3_mcs_smax}. For Case 1, using $s_{\max}=5n+10$ (MCS-1) performs similarly to using $s_{\max}=10n$ (MCS-6). For Case 2, which has 32 variables, using $s_{\max}=5n+10$ (MCS-1) converges slightly faster than using $s_{\max}=10n$ (MCS-6), and finally obtains a higher NPV. This indicates that a small number of levels is enough for these cases. 

Fig.\ \ref{fig:subfig:e2_0_mcs_ls} and Fig.\ \ref{fig:subfig:e2_3_mcs_ls} show that local search plays an important role in MCS, without it the convergence speed decreases significantly.

\begin{figure}[htbp]
\centering 
\subfigure[PRO-01, Case 1]{ 
    \label{fig:subfig:rate1_1} 
    \includegraphics[width=0.3\textwidth]{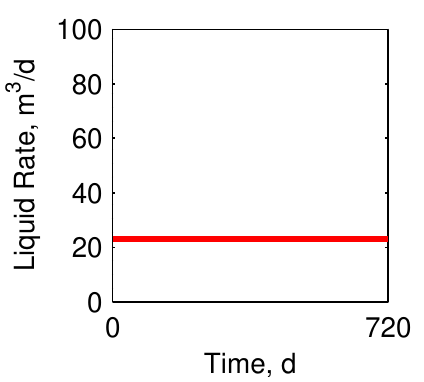}
    }
\subfigure[PRO-01, Case 2]{ 
    \label{fig:subfig:rate1_3} 
    \includegraphics[width=0.3\textwidth]{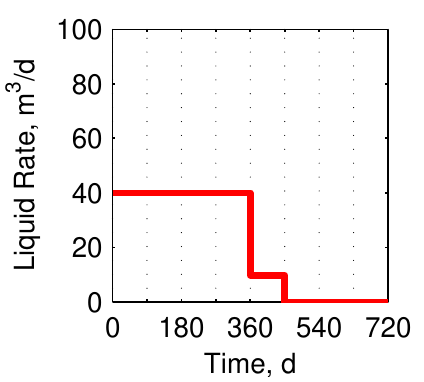}
    }\\
\subfigure[PRO-02, Case 1]{ 
    \label{fig:subfig:rate2_1} 
    \includegraphics[width=0.3\textwidth]{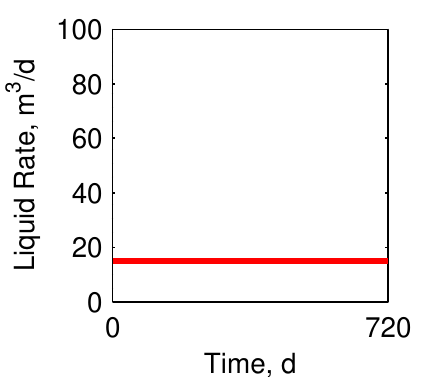}
    }
\subfigure[PRO-02, Case 2]{ 
    \label{fig:subfig:rate2_3} 
    \includegraphics[width=0.3\textwidth]{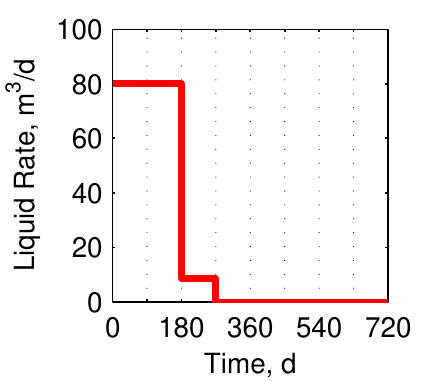}
    }
  \caption{The optimal well control strategies found for Case 1 and Case 2 of Example 3.} 
  \label{fig:opt_control} 
\end{figure}

Fig.\ \ref{fig:opt_control} presents the optimum controls for wells PRO-01 and PRO-02 under different control frequencies. We omit the results for well PRO-03 and PRO-04 because the reservoir is symmetric. The optimum controls become more like a bang-bang solution for all wells with an increase in the number of control steps. It is worth noting that the optimum controls for Case 1 are significantly different that those for Case 2. This reflects the different production strategies for wells using a static rate compared to using dynamic well controls in water flooding reservoirs. 

\subsection{Example 4, joint well placement and control optimization}
\label{sec:44_jop}

\subsubsection{Reservoir model description}

\begin{figure*}[htb]
\centering 
\subfigure[Permeability($\log_{10}$mD)]{ 
    \label{fig:subfig:both_perm} 
    \includegraphics[width=0.45\textwidth]{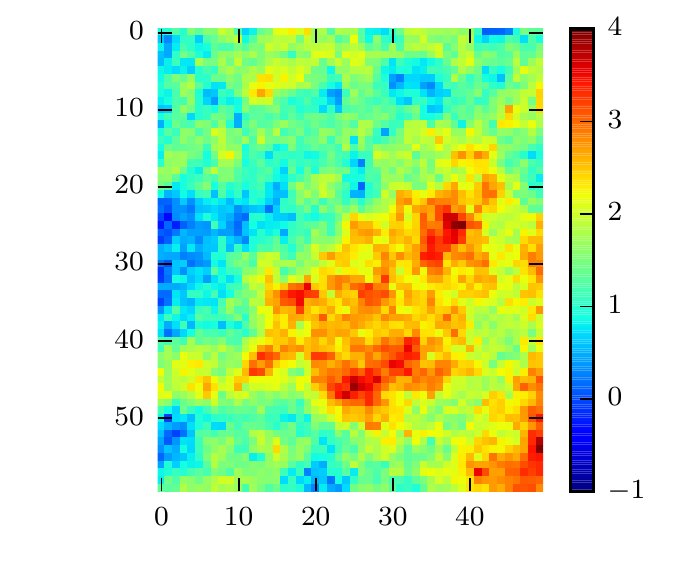}
    }
\subfigure[Porosity]{ 
    \label{fig:subfig:both_poro} 
    \includegraphics[width=0.45\textwidth]{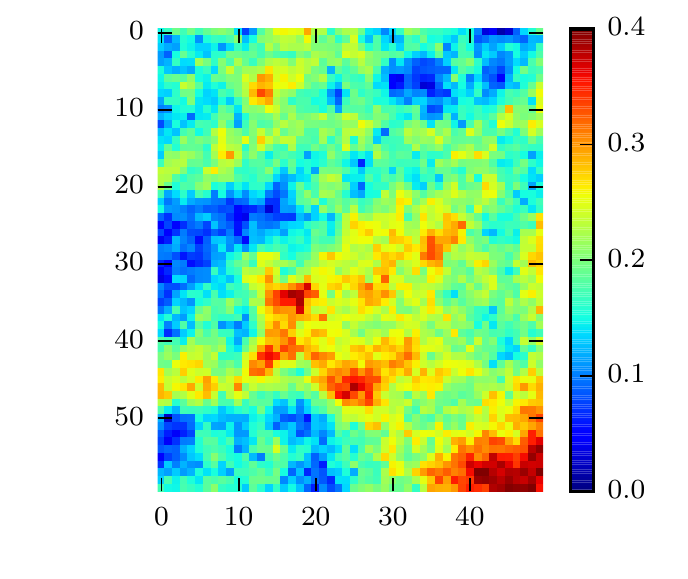}
    }
\caption{Permeability and porosity fields of the SPE10 model used in Example 4.}
\label{fig:both_spe10}
\end{figure*}

This example use a 2D reservoir model with the permeability and porosity fields taken from the third layer of the SPE10 benchmark model \citep{christie_tenth_2001}. It consists of $60\times 50$ grid cells and the size of each grid cell is $\rm 32m\times 32m \times 10m$. We consider an oil-water two phase flow in this model and the initial oil saturation is 0.8. Fig.\ \ref{fig:both_spe10} shows the permeability and porosity fields of the model.

The optimization problem is to place four wells in the reservoir, including two production wells (P1, P2) and two injection wells (I1, I2). All wells are controlled via BHP that is updated every two years. The production period for this example is 10 years. Thus, there are two location variables and five control variables per well and 28 variables in total. Only bound constraints are considered in this example. The economic and optimization parameters are summarized in Table \ref{tab:e3_Epara} and Table \ref{tab:e3_Opara} respectively. 
Both the simultaneous procedure and the sequential procedure are used for this example.

\begin{table}[htb]
\caption{Economic parameters used in Example 4.}
\label{tab:e3_Epara}
\centering
\begin{tabular}{cc}
\hline\noalign{\smallskip}
Parameter & Value \\
\noalign{\smallskip}\hline\noalign{\smallskip}
Oil revenue ($\rm{USD/ m^3}$) & 503.2 \\
Water-production cost ($\rm{USD/ m^3}$) & 75.5 \\
Water-injection cost ($\rm{USD/ m^3}$) & 50.3 \\
Annual discount rate & 0 \\
\noalign{\smallskip}\hline
\end{tabular}
\end{table}

\begin{table}[htb]
\caption{Optimization parameters used in Example 4.}
\label{tab:e3_Opara}
\centering
\begin{tabular}{ccccc}
\hline\noalign{\smallskip}
Parameter & P1 & P2 & I1 & I2 \\
\noalign{\smallskip}\hline\noalign{\smallskip}
Initial location & (60,25) & (1,25) & (30,1) & (30,50) \\
Initial BHP ($\rm bar$) & 175 & 175 & 362.5 & 362.5 \\
Minimum BHP ($\rm bar$) & 100 & 100 & 275 & 275 \\
Maximum BHP ($\rm bar$) & 250 & 250 & 450 & 450 \\
\noalign{\smallskip}\hline
\end{tabular}
\end{table}

\subsubsection{Simultaneous procedure}

\begin{figure}[htb]
\centering
  \includegraphics[width=0.45\textwidth]{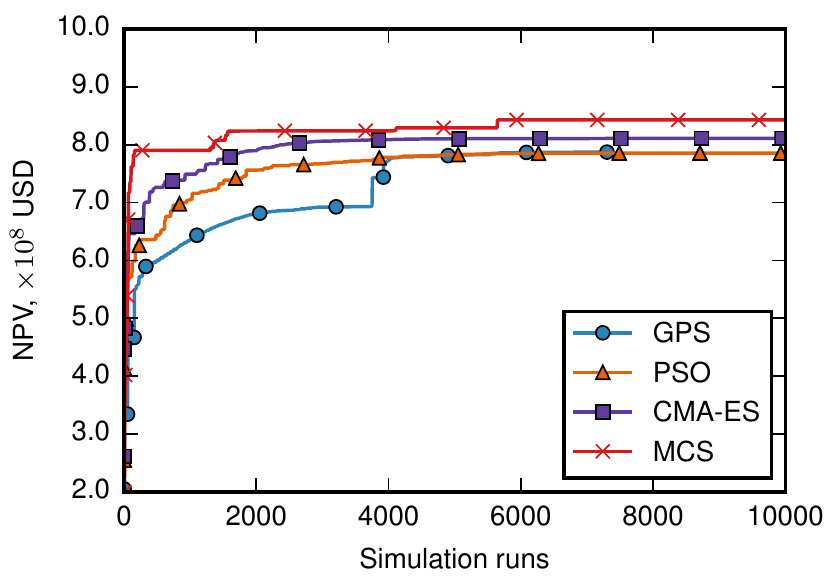}
\caption{Optimization performance of simultaneous procedure for Example 4. For PSO and CMA-ES, the solid lines depict the median NPV. MCS here is a label for MCS-1.}
\label{fig:e3_algo}
\end{figure}

The simultaneous procedure optimizes over the well locations and controls simultaneously.
For this problem, we optimize the locations and control parameters of the 4 wells. Each well has 2 location variables and 5 control variables, Thus there are 28 variables in total.
Given this problem's complexity, we set the maximum number of simulation runs for this example  to be 10000.
The maximum, minimum, mean, median, and standard deviation of NPV for each algorithm is given in Table \ref{tab:e3_result}. From the table we can see that MCS-1 obtains the highest NPV value after 10000 simulation runs.
The average NPV for PSO and CMA-ES are in the middle, while GPS performs the worst.
Plots of the NPV of the four algorithms versus the number of simulation runs are shown in Fig.\ \ref{fig:e3_algo}. 

\begin{table}[htb]
\caption{Results of simultaneous procedure for Example 4. Values shown are NPV in $\$ \times 10^8$ USD. }
\label{tab:e3_result}
\centering
\subtable[Deterministic algorithms (MCS, GPS)]{
\label{tab:e3_result_det}
\begin{tabular}{p{2cm}|p{8cm}}
\hline
Algorithm & \hfil NPV \\
\hline
MCS-1  & \hfil 8.43 \\
MCS-2  & \hfil 7.50 \\
MCS-3  & \hfil 8.48 \\
MCS-4  & \hfil 7.62 \\
MCS-5  & \hfil 7.69 \\
MCS-6  & \hfil 7.94 \\
MCS-7  & \hfil 6.81 \\
GPS    & \hfil 7.87 \\
\hline
\end{tabular}
}

\subtable[Stochastic algorithms (PSO, CMA-ES)]{        
\label{tab:e3_result_sto}
\begin{tabular}{p{2cm}|p{1cm}|p{1cm}p{1cm}p{1cm}p{1cm}p{1cm}}
\hline
Algorithm & Trials & Max & Min & Mean & Median & Std. \\
\hline
PSO    &10 &8.50&7.34&8.08&8.22&0.48\\
CMA-ES &10 &8.45&7.58&8.13&8.16&0.28\\
\hline
\end{tabular}
}
\end{table}

\begin{figure}[htb]
\centering 
  \subfigure[PSO]{ 
    \label{fig:subfig:e3_pso} 
    \includegraphics[width=0.45\textwidth]{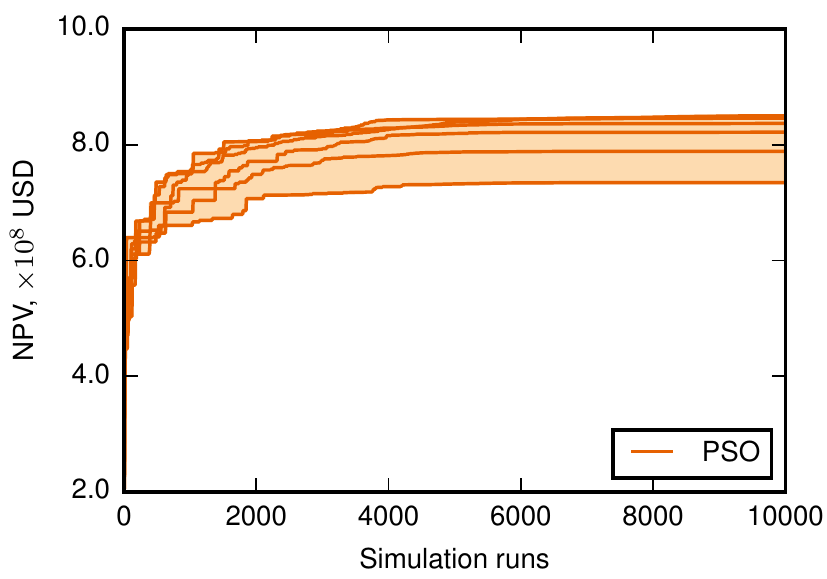}} 
  \subfigure[CMA-ES]{ 
    \label{fig:subfig:e3_cmaes} 
    \includegraphics[width=0.45\textwidth]{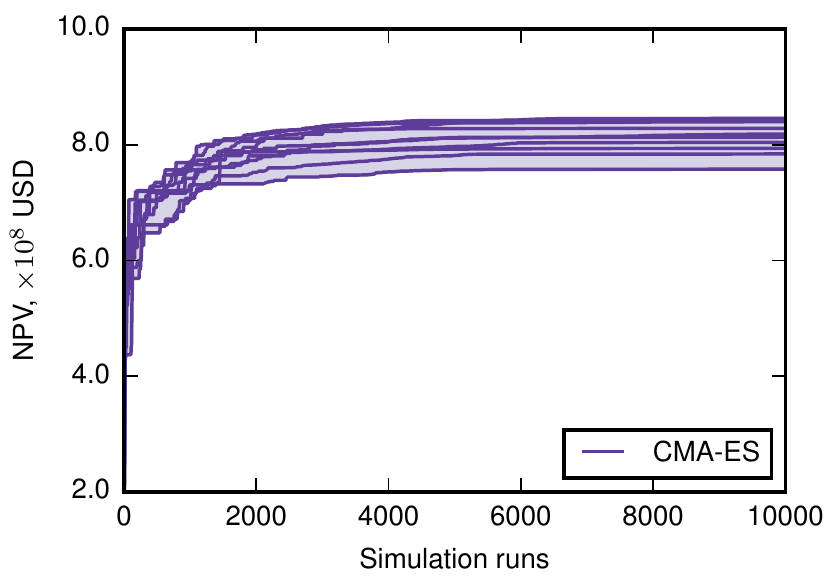}} 
  \caption{The range of NPV found amongst the trials of PSO and CMA-ES for Example 4. Each solid line represents a trial.} 
  \label{fig:e3_psocmaes} 
\end{figure}

Fig.\ \ref{fig:e3_algo} shows that MCS converges fastest, followed by CMA-ES, PSO, and GPS in that order. Unlike Example 1 and 2, the convergence speed of GPS is slowest among all algorithms. The NPV of GPS has a jump at about 4000 simulation runs. It appears that at this point GPS jumps from a local optima.

Fig.\ \ref{fig:e3_psocmaes} shows the range of NPV for the trials of PSO and CMA-ES. In this figure, the areas between the maximum and minimum NPV are shaded for PSO and CMA-ES. It is clear that the NPV obtained by PSO and CMA-ES has a high variation for this example. 

\begin{figure}[htb]
\centering 
  \subfigure[Initialization list]{ 
    \label{fig:subfig:e3_mcs_init} 
    \includegraphics[width=0.45\textwidth]{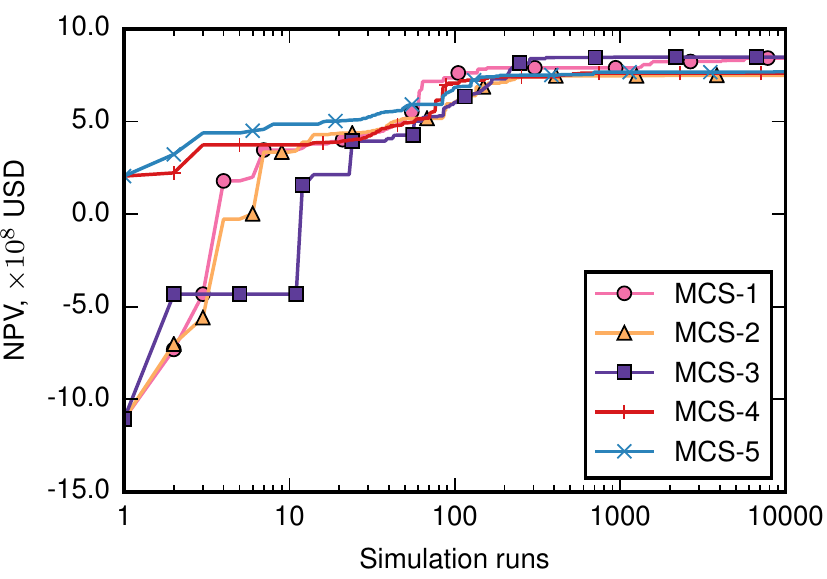}} 
  \subfigure[Levels]{ 
    \label{fig:subfig:e3_mcs_smax} 
    \includegraphics[width=0.45\textwidth]{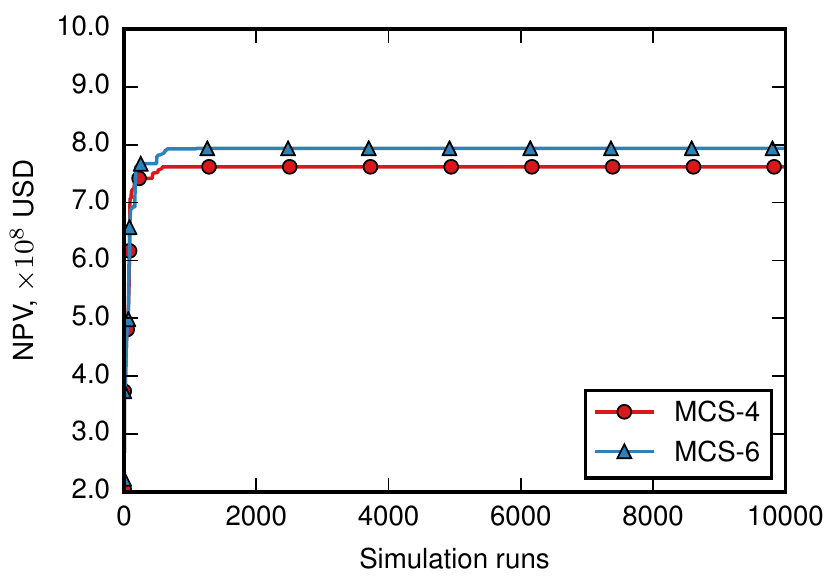}} 
  \subfigure[Local search]{ 
    \label{fig:subfig:e3_mcs_ls} 
    \includegraphics[width=0.45\textwidth]{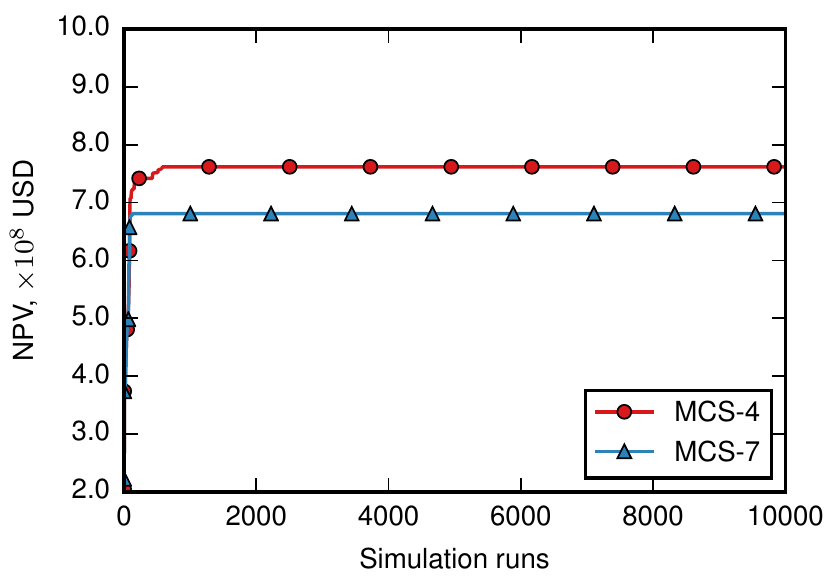}}   
  \caption{The performance of different configurations of MCS for Example 4.} 
  \label{fig:e3_mcs} 
\end{figure}

As with Examples 1, 2, and 3, we use 7 MCS configurations, and the results are compared within the 3 groups in Fig.\ \ref{fig:e3_mcs}. Fig.\ \ref{fig:subfig:e3_mcs_init} shows the performance of MCS with different initialization lists. We use a semilog plot to make this figure clearer. 
The initialization list with a good initial guess (MCS-4 and MCS-5), starts its search from a relatively high NPV, but obtains a NPV lower than MCS-1, which uses the default initialization list with a boundary point. For this example, MCS-1, MCS-2, and MCS-3 recover quite quickly from the bad initial guess, and are able to converge more quickly.

\begin{figure}[htb]
\centering
  \includegraphics[width=0.45\textwidth]{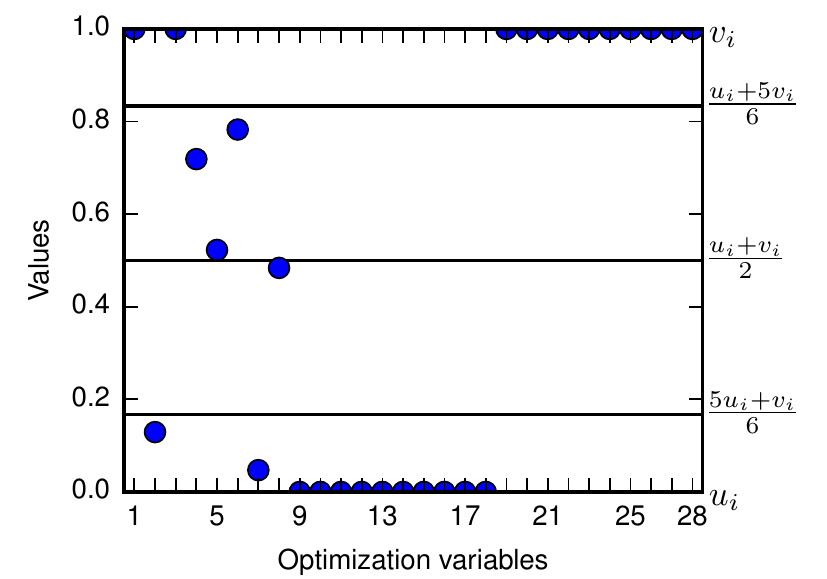}
\caption{Normalized boundary, initialization lists, and the optimal solution for Example 4.}
\label{fig:e3_ila}
\end{figure}

We can see that for this example, the initialization list without boundary points (MCS-2) performs unsatisfactorily both in terms of the convergence rate and the final NPV. This is because the optimal solution for this example lies near the boundary, as shown in Fig.\ \ref{fig:e3_ila}.
Using a line search to generate the initialization list (MCS-3) ultimately obtains the highest NPV for this example. 

Fig.\ \ref{fig:subfig:e3_mcs_smax} shows the performance of MCS with different numbers of maximum levels. Using $s_{\max}=10n$ (MCS-6) outperforms choosing $s_{\max}=5n+10$ (MCS-4) for this example.

The performance of MCS with and without local search is shown in Fig.\ \ref{fig:subfig:e3_mcs_ls}. MCS without local search (MCS-7) is clearly inferior.

\subsubsection{Sequential procedure}
The sequential procedure decouples the joint problem into two separate subproblems. 
For the well placement optimization subproblem, we optimize the locations of the 4 wells under an assumed control scheme. This gives 8 optimization variables. For the well control optimization subproblem, we optimize the well controls of the 4 wells with assumed well locations. Each well has 5 control steps, this gives 20 optimization variables in total. 
The maximum number of simulation runs for each well placement optimization stage is 60, while for each well control optimization stage the maximum number of simulation runs is set to 140. And the maximum number of simulation runs for the problem in total is 5000. 
This allows us to iterate 25 times between the well placement and well control optimization phases.
Based on the results of the previous section we use MCS-1 in the sequential procedure.

\begin{figure}[htb]
\centering 
  \subfigure[]{ 
    \label{fig:subfig:e3_seq1} 
    \includegraphics[width=0.45\textwidth]{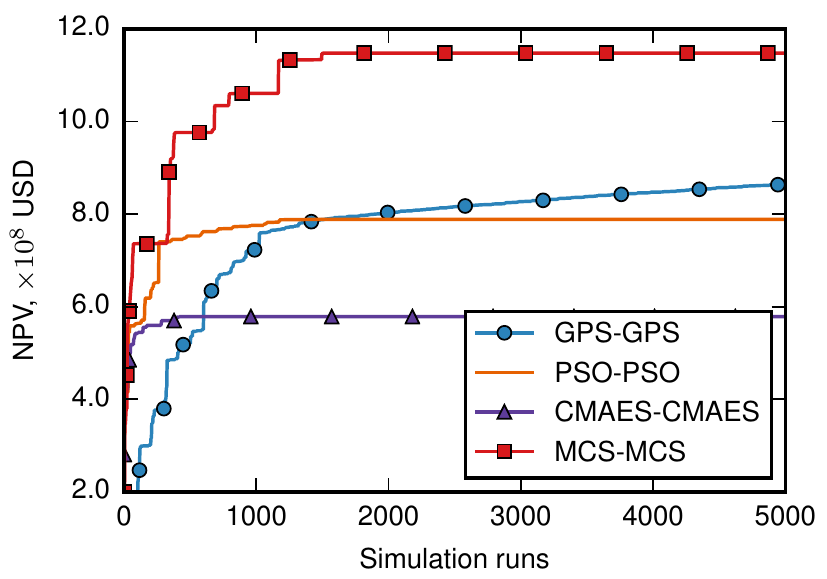}} 
  \subfigure[]{ 
    \label{fig:subfig:e3_seq2} 
    \includegraphics[width=0.45\textwidth]{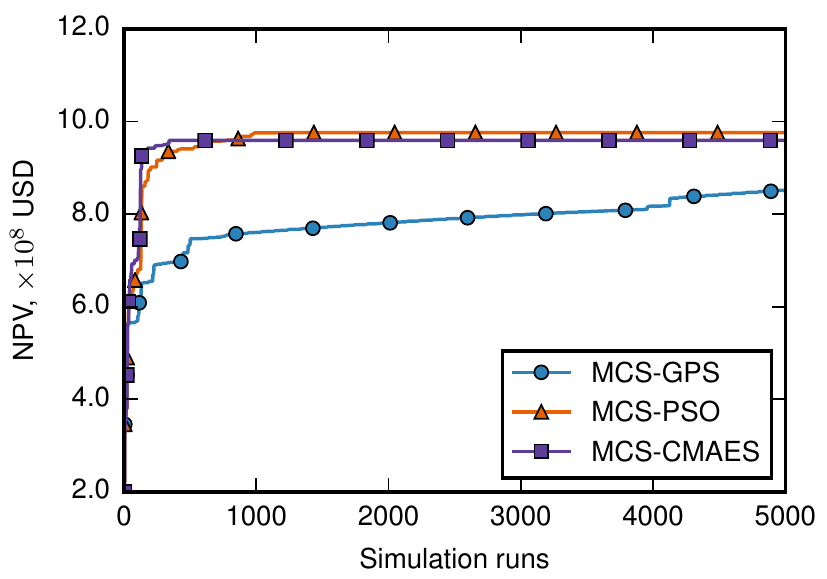}} 
  \subfigure[]{ 
    \label{fig:subfig:e3_seq3} 
    \includegraphics[width=0.45\textwidth]{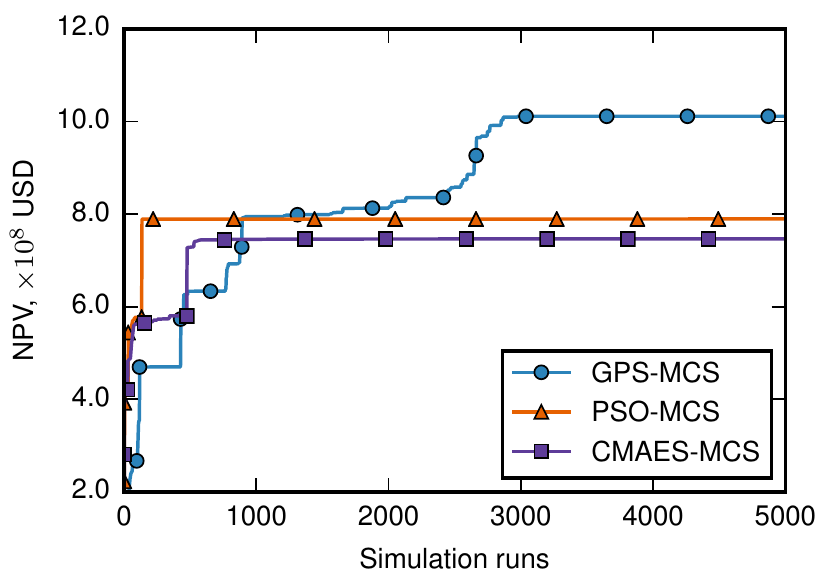}}   
  \caption{Optimization performance of the sequential procedure for Example 4 using different algorithm combinations.} 
  \label{fig:e3_seq} 
\end{figure}

The maximum, minimum, mean, median, and standard deviation of NPV for each algorithm combination is given in Table \ref{tab:e3_result_seq}. Plots of the NPV versus the number of simulation runs for each approach are shown in Fig.\ \ref{fig:e3_seq}. From Table \ref{tab:e3_result_seq} and Fig.\ \ref{fig:e3_seq} we see that MCS-MCS converges faster than the other combinations and obtained the highest NPV value at the end of the optimization. GPS-MCS converges slowly at the early stage, but it ultimately obtains the second highest NPV. The combinations which contain stochastic algorithms, especially CMA-ES, perform unsatisfactorily for this example.

\begin{table}[htb]
\caption{Results of the sequential procedure for Example 4. Values shown are NPV in $\$ \times 10^8$ USD. }
\label{tab:e3_result_seq}
\centering
\begin{tabular}{p{3.5cm}|p{1cm}|p{1cm}p{1cm}p{1cm}p{1cm}p{1cm}}
\hline
Algorithm & Trials & Max & Min & Mean & Median & Std. \\
\hline
MCS-MCS         &1  &11.48&11.48&11.48&11.48&0\\
GPS-GPS         &1  &8.65&8.65&8.65&8.65&0\\
PSO-PSO         &5  &8.11&6.64&7.54&7.89&0.79\\
CMA-ES-CMA-ES   &5  &6.58&5.75&6.04&5.79&0.47\\
\hline
MCS-GPS         &1  &8.53&8.53&8.53&8.53&0\\
MCS-PSO         &5  &10.78&9.55&10.03&9.76&0.66\\
MCS-CMA-ES      &5  &9.99&9.48&9.69&9.59&0.27\\
\hline
GPS-MCS         &1  &10.11&10.11&10.11&10.11&0\\
PSO-MCS         &5  &9.01&6.63&7.84&7.90&1.19\\
CMA-ES-MCS      &5  &8.98&5.80&7.42&7.47&1.59\\
\hline
\end{tabular}
\end{table}

\begin{figure}[htb]
\centering
  \includegraphics[width=0.9\textwidth]{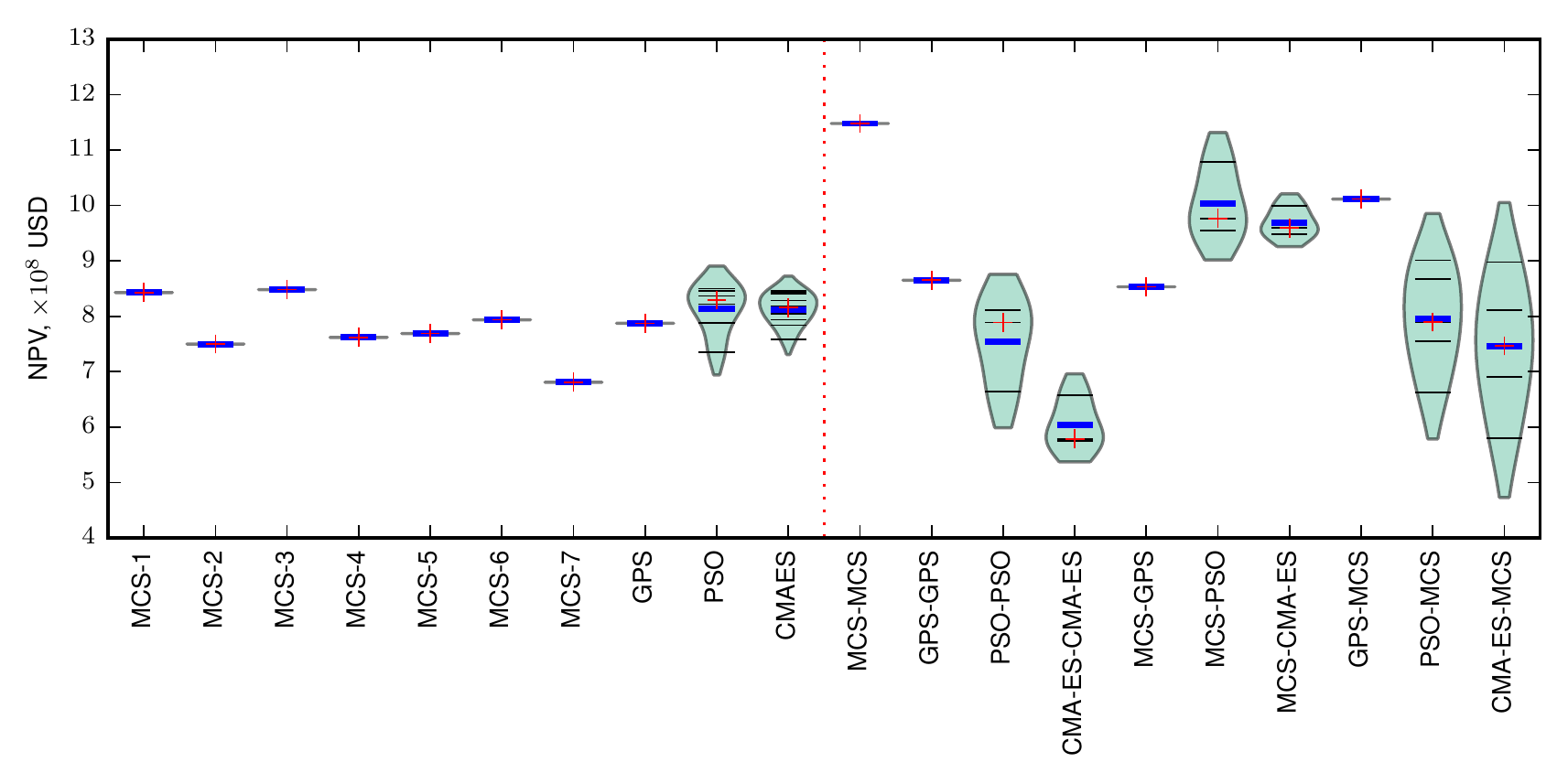}
\caption{Beanplots of the final NPV values for the simultaneous and the sequential procedure for Example 4. The beanplots to left of the dotted vertical line gives the results obtained by algorithms using the simultaneous procedure, and the beanplots to the right gives the results obtained by the sequential procedure. The individual horizontal lines show the NPV obtained by each trial. The horizontal bold line and the marker `+' denote the mean and median of all trials, respectively.}
\label{fig:seqvssim}
\end{figure}

We also compare the optimal NPV obtained using the simultaneous and the sequential procedures using beanplots in Fig.\ \ref{fig:seqvssim}. A beanplot \citep{kampstra_beanplot:_2008} promotes visual comparison of univariate data between groups. In a beanplot, the individual observations are shown as small lines in a one-dimensional scatter plot. In addition, the estimated density of the distributions is visible and the mean (bold line) and median (marker `+') are shown.

From Fig.\ \ref{fig:seqvssim} we can see that, with the simultaneous procedure, the final NPV values obtained by all algorithms are less than $10^9$USD. With the sequential procedure, MCS-MCS, GPS-MCS, and MCS-PSO can obtain a NPV value higher than $10^9$USD.

\begin{figure}[htb]
\centering
  \includegraphics[width=0.5\textwidth]{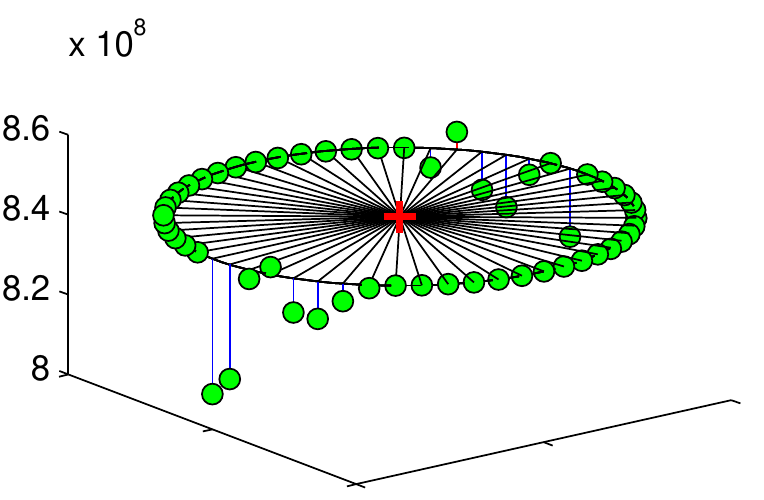}
\caption{Illustration of perturbation around the optimal solution. The marker `+' denotes the optimal point, the bubbles denote the perturbation points. Each line denotes a dimension. The height of each point shows the value of NPV at this point.}
\label{fig:e3_local}
\end{figure}

Indeed the simultaneous algorithm becomes trapped in a sub-optimal solution. In Fig.\ \ref{fig:e3_local} we show the NPV obtained around the candidate solution in each of the 28 dimensions by sampling at $x_i \pm \Delta x_i$ for $i=1,\cdots,28$. For the well location variables we choose $\Delta x_i=1$ and for the well control variables we set $\Delta x_i$ to 1\% of the range of the $i$th variable. We see that in most of directions the NPV remains constant. The NPV is lower in a few directions. And only gives an improved NPV in a few directions. Hence the algorithms have difficulty finding an ascent direction. 

\begin{figure}[htb]
\centering 
  \subfigure[Simultaneous procedure]{ 
    \label{fig:subfig:e3_soil_sim} 
    \includegraphics[width=0.45\textwidth]{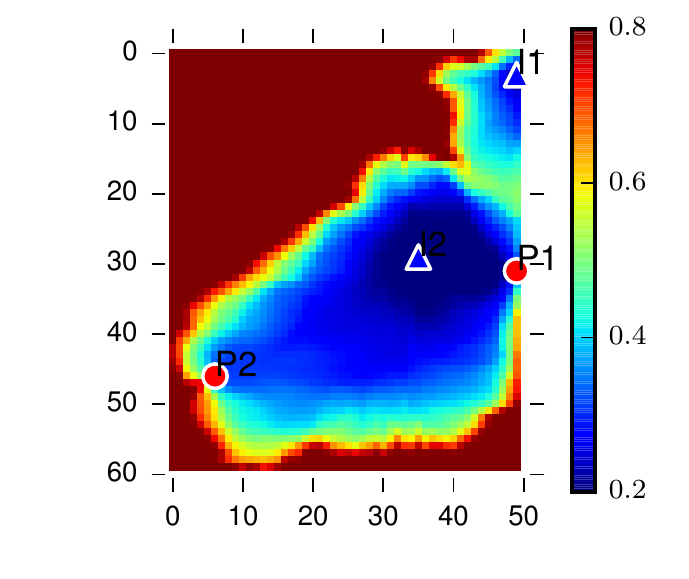}} 
  \subfigure[Sequential procedure]{ 
    \label{fig:subfig:e3_soil_seq} 
    \includegraphics[width=0.45\textwidth]{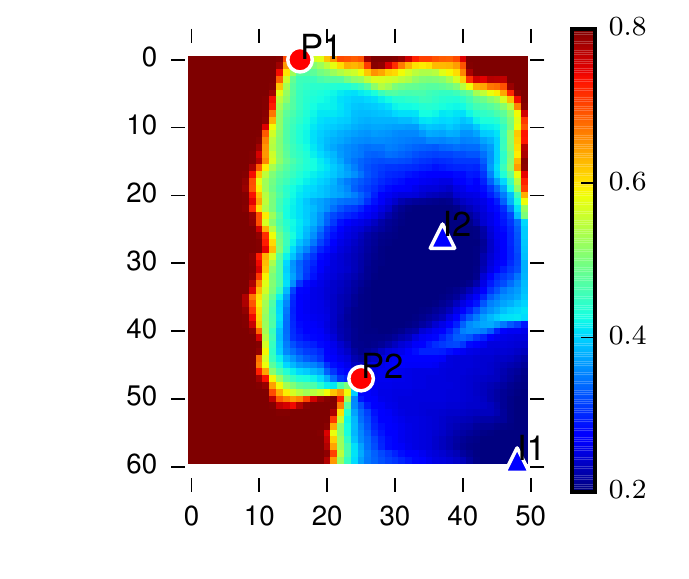}}  
  \caption{The optimal well placement and the final oil saturation distribution for the simultaneous procedure and the sequential procedure for Example 4.} 
  \label{fig:e3_soil} 
\end{figure}

The optimal well placement and the final oil saturation distribution obtained with the simultaneous and sequential procedures are shown in Fig.\ \ref{fig:e3_soil}. The corresponding optimal controls of each well are given in Table \ref{tab:e3_control}. The optimal well locations obtained by the simultaneous procedure are significantly different from the locations obtained by the sequential procedure. From the final oil saturation distribution, we can see that, the locations obtained by the sequential procedure give a larger sweep area. The optimal controls obtained by the simultaneous procedure are similar to those found by the sequential algorithm.

In theory, for a joint well placement and control optimization problem, the simultaneous procedure is able to find the global optima, but this is not guaranteed for the sequential procedure since the optimal location of each well depends on how the well is operated and vice-versa. The simultaneous procedure, with a larger number of optimization variables, makes the joint problem more difficult. It requires a higher computational budget and has a higher risk of falling into a local optima and achieving a suboptimal solution, especially for a larger scale problem. The sequential procedure, decouples a hard joint problem into two easier subproblems, and hopes to approach the global optima iteratively. In general, the sequential procedure is worth considering in practice.

\begin{table}[htb]
\caption{The optimal control of each well for the simultaneous procedure and the sequential procedure for Example 4. Values shown are BHP in bar.}
\label{tab:e3_control}
\centering
\subtable[Simultaneous procedure]{
\label{tab:e3_c_sim}
\begin{tabular}{cc}
\hline
Well & BHP (bar) \\
\hline
P1 & $[100,100,100,100,100]$ \\
P2 & $[100,100,100,100,100]$ \\
I1 & $[450,450,450,450,450]$ \\
I2 & $[450,450,450,450,450]$ \\
\hline
\end{tabular}
}

\subtable[Sequential procedure]{
\label{tab:e3_c_seq}
\begin{tabular}{cc}
\hline
Well & BHP (bar) \\
\hline
P1 & $[103,103,103,103,103]$ \\
P2 & $[103,103,103,103,103]$ \\
I1 & $[450,447,447,447,450]$ \\
I2 & $[447,447,447,447,447]$ \\
\hline
\end{tabular}
}
\end{table}

\subsection{Summary}

Our test results show that MCS is strongly competitive with existing algorithms for well placement, well control, and joint optimization problems. In all 4 examples, MCS offers good convergence speed, especially when the number of simulation runs is limited. 
Moreover, MCS does not suffer from the inherent variability of the stochastic algorithms.
Based on the results of the examples, for placement and control optimization we suggest a MCS configuration which uses a line search to generate the initialization list. The number of levels $s_{\max}=5n+10$ is enough for most problems but a higher $s_{\max}$ should be used for difficult problems. Local search is an important part of MCS, and is highly recommended.

\section{Concluding Remarks}
\label{sec:6_concl}

In this paper, we applied the multilevel coordinate search algorithm for four typical oil field development optimization problems. The problems include well placement optimization, well control optimization, and joint optimization of well placement and control. The performance of MCS has been compared with generalized pattern search, particle swarm optimization, and covariance matrix adaptation evolution strategy through several case studies including both synthetic and real reservoirs. The results presented here demonstrate that the MCS algorithm is strongly competitive, and outperforms the other algorithms in most cases, especially for the joint optimization problem. MCS has significant advantages in solving optimization problems with a limited number of simulation runs. In addition, MCS does not suffer from the inherent variability of the stochastic approaches.

For joint well placement and well control optimization problem, both the simultaneous procedure and the sequential procedure were considered. In our example, the sequential procedure finds the best solution. Although the simultaneous procedure can theoretically obtain the global optima, the sequential procedure is worth considering in practice. The sequential procedure decouples a difficult joint problem to two easier separate subproblem, decreases the number of optimization variables, make the problem easier to solve and decreases the risk of the algorithm falling into a local optima. Among all algorithm combinations considered in this paper, MCS-MCS showed best performance both in terms of convergence speed and final NPV value in the sequential procedure. 

MCS has shown its potential in our work, but more research is needed. Future work includes applying the MCS algorithm to realistic large-scale oil field cases.
This will involve an extension of MCS to handle linearly and nonlinearly constrained problems, 
possibly by a penalty approach.

\section*{Acknowledgments}
The authors acknowledge funding from the Natural Sciences and Engineering Research Council of Canada (NSERC) Discovery Grant Program, the National Science and Technology Major Project of the Ministry of Science and Technology of China (2011ZX05011-002), the Foundation for Outstanding Young Scientist in Shandong Province (Grant no. BS2014NJ011), and the program of China Scholarships Council (No. 201406450017).





\begin{thebibliography}{70}
\expandafter\ifx\csname natexlab\endcsname\relax\def\natexlab#1{#1}\fi
\expandafter\ifx\csname url\endcsname\relax
  \def\url#1{\texttt{#1}}\fi
\expandafter\ifx\csname urlprefix\endcsname\relax\def\urlprefix{URL }\fi

\bibitem[{AlQahtani et~al.(2014)AlQahtani, Alzahabi, Spinner, and
  Soliman}]{alqahtani_computational_2014}
AlQahtani, G.~D., Alzahabi, A., Spinner, T., Soliman, M.~Y., 2014. A
  computational comparison between optimization techniques for wells placement
  problem: Mathematical formulations, genetic algorithms and very fast
  simulated annealing. Journal of Materials Science and Chemical Engineering
  02~(10), 59--73.

\bibitem[{Asadollahi et~al.(2014)Asadollahi, N{\ae}vdal, Dadashpour, and
  Kleppe}]{asadollahi_production_2014}
Asadollahi, M., N{\ae}vdal, G., Dadashpour, M., Kleppe, J., Feb. 2014.
  Production optimization using derivative free methods applied to {Brugge}
  field case. Journal of Petroleum Science and Engineering 114, 22--37.

\bibitem[{Asheim(1988)}]{asheim_maximization_1988}
Asheim, H., 1988. Maximization of water sweep efficiency by controlling
  production and injection rates. In: European {Petroleum} {Conference}.
  Society of Petroleum Engineers.

\bibitem[{Audet and Dennis(2002)}]{audet_analysis_2002}
Audet, C., Dennis, J., Jan. 2002. Analysis of {Generalized} {Pattern}
  {Searches}. SIAM Journal on Optimization 13~(3), 889--903.

\bibitem[{Auger and Hansen(2005)}]{auger_restart_2005}
Auger, A., Hansen, N., Sep. 2005. A restart {CMA} evolution strategy with
  increasing population size. In: The 2005 {IEEE} {Congress} on {Evolutionary}
  {Computation}, 2005. Vol.~2. pp. 1769--1776.

\bibitem[{Bangerth et~al.(2006)Bangerth, Klie, Wheeler, Stoffa, and
  Sen}]{bangerth_on_2006}
Bangerth, W., Klie, H., Wheeler, M.~F., Stoffa, P.~L., Sen, M.~K., Aug 2006. On
  optimization algorithms for the reservoir oil well placement problem.
  Computational Geosciences 10~(3), 303--319.

\bibitem[{Bellout et~al.(2012)Bellout, Ciaurri, Durlofsky, Foss, and
  Kleppe}]{bellout_joint_2012}
Bellout, M.~C., Ciaurri, D.~E., Durlofsky, L.~J., Foss, B., Kleppe, J., Jul.
  2012. Joint optimization of oil well placement and controls. Computational
  Geosciences 16~(4), 1061--1079.

\bibitem[{Bouzarkouna et~al.(2012)Bouzarkouna, Ding, and
  Auger}]{bouzarkouna_well_2012}
Bouzarkouna, Z., Ding, D.~Y., Auger, A., Jan. 2012. Well placement optimization
  with the covariance matrix adaptation evolution strategy and meta-models.
  Computational Geosciences 16~(1), 75--92.

\bibitem[{Brouwer and Jansen(2004)}]{brouwer_dynamic_2004}
Brouwer, D.~R., Jansen, J.~D., 2004. Dynamic {Optimization} of {Waterflooding}
  {With} {Smart} {Wells} {Using} {Optimal} {Control} {Theory}. SPE Journal
  9~(04), 391--402.

\bibitem[{Chavent(1974)}]{chavent_identification_1974}
Chavent, G., 1974. Identification of functional parameters in partial
  differential equations. In: Joint {Automatic} {Control} {Conference}. pp.
  155--156.

\bibitem[{Chen et~al.(2012)Chen, Li, and Reynolds}]{chen_robust_2012}
Chen, C., Li, G., Reynolds, A., 2012. Robust constrained optimization of
  short-and long-term net present value for closed-loop reservoir management.
  SPE Journal 17~(03), 849--864.

\bibitem[{Chen et~al.(1974)Chen, Gavalas, Seinfeld, and
  Wasserman}]{chen_new_1974}
Chen, W.~H., Gavalas, G.~R., Seinfeld, J.~H., Wasserman, M.~L., 1974. A {New}
  {Algorithm} for {Automatic} {History} {Matching}. Society of Petroleum
  Engineers Journal 14~(06), 593--608.

\bibitem[{Christie and Blunt(2001)}]{christie_tenth_2001}
Christie, M.~A., Blunt, M.~J., 2001. Tenth {SPE} comparative solution project:
  {A} comparison of upscaling techniques. In: {SPE} {Reservoir} {Simulation}
  {Symposium}. Society of Petroleum Engineers.

\bibitem[{Ciaurri et~al.(2011)Ciaurri, Mukerji, and
  Durlofsky}]{ciaurri_derivative-free_2011}
Ciaurri, D.~E., Mukerji, T., Durlofsky, L.~J., 2011. Derivative-{Free}
  {Optimization} for {Oil} {Field} {Operations}. In: Yang, X.-S., Koziel, S.
  (Eds.), Computational {Optimization} and {Applications} in {Engineering} and
  {Industry}. No. 359 in Studies in {Computational} {Intelligence}. Springer
  Berlin Heidelberg, pp. 19--55.

\bibitem[{Clerc(2006)}]{clerc_stagnation_2006}
Clerc, M., 2006. Stagnation analysis in particle swarm optimization or what
  happens when nothing happens. Tech. Rep. CSM-460, Department of Computer
  Science, University of Essex.

\bibitem[{Emerick et~al.(2009)Emerick, Silva, Messer, Almeida, Szwarcman,
  Pacheco, and Vellasco}]{emerick_well_2009}
Emerick, A.~A., Silva, E., Messer, B., Almeida, L.~F., Szwarcman, D., Pacheco,
  M. A.~C., Vellasco, M. M. B.~R., 2009. Well placement optimization using a
  genetic algorithm with nonlinear constraints. In: {SPE} reservoir simulation
  symposium. Society of Petroleum Engineers.

\bibitem[{Fonseca et~al.(2014)Fonseca, Stordal, Leeuwenburgh, Van~den Hof, and
  Jansen}]{fonseca_robust_2014}
Fonseca, R.~M., Stordal, A.~S., Leeuwenburgh, O., Van~den Hof, P. M.~J.,
  Jansen, J.~D., 2014. Robust ensemble-based multi-objective optimization. In:
  {ECMOR} {XIV}-14th {European} conference on the mathematics of oil recovery.

\bibitem[{Forouzanfar et~al.(2010)Forouzanfar, Li, and
  Reynolds}]{forouzanfar_two-stage_2010}
Forouzanfar, F., Li, G., Reynolds, A.~C., 2010. A two-stage well placement
  optimization method based on adjoint gradient. In: {SPE} {Annual} {Technical}
  {Conference} and {Exhibition}. Society of Petroleum Engineers.

\bibitem[{Forouzanfar et~al.(2015)Forouzanfar, Poquioma, and
  Reynolds}]{forouzanfar_covariance_2015}
Forouzanfar, F., Poquioma, W.~E., Reynolds, A.~C., 2015. A covariance matrix
  adaptation algorithm for simultaneous estimation of optimal placement and
  control of production and water injection wells. In: {SPE} Reservoir
  Simulation Symposium. Society of Petroleum Engineers ({SPE}).

\bibitem[{Forouzanfar and Reynolds(2014)}]{forouzanfar_joint_2014}
Forouzanfar, F., Reynolds, A.~C., Jul. 2014. Joint optimization of number of
  wells, well locations and controls using a gradient-based algorithm. Chemical
  Engineering Research and Design 92~(7), 1315--1328.

\bibitem[{Forouzanfar et~al.(2012)Forouzanfar, Reynolds, and
  Li}]{forouzanfar_optimization_2012}
Forouzanfar, F., Reynolds, A.~C., Li, G., May 2012. Optimization of the well
  locations and completions for vertical and horizontal wells using a
  derivative-free optimization algorithm. Journal of Petroleum Science and
  Engineering 86–87, 272--288.

\bibitem[{Gao et~al.(2006)Gao, Zafari, and Reynolds}]{gao_quantifying_2006}
Gao, G., Zafari, M., Reynolds, A.~C., Dec. 2006. Quantifying uncertainty for
  the {PUNQ}-s3 problem in a bayesian setting with {RML} and {EnKF}. {SPE}
  Journal 11~(04), 506--515.

\bibitem[{GeoQuest(2014)}]{geoquest_eclipse_2014}
GeoQuest, S., 2014. {ECLIPSE} reference manual. Schlumberger, Houston, Texas.

\bibitem[{Handels et~al.(2007)Handels, Zandvliet, Brouwer, and
  Jansen}]{handels_adjoint-based_2007}
Handels, M., Zandvliet, M., Brouwer, R., Jansen, J.~D., 2007. Adjoint-based
  well-placement optimization under production constraints. In: {SPE} reservoir
  simulation symposium. Society of Petroleum Engineers.

\bibitem[{Hansen and Kern(2004)}]{hansen_evaluating_2004}
Hansen, N., Kern, S., Jan. 2004. Evaluating the {CMA} evolution strategy on
  multimodal test functions. In: Yao, X., Burke, E.~K. (Eds.), Parallel Problem
  Solving from Nature - {PPSN} {VIII}. No. 3242 in Lecture Notes in Computer
  Science. Springer Berlin Heidelberg, pp. 282--291.

\bibitem[{Humphries and Haynes(2015)}]{humphries_joint_2015}
Humphries, T., Haynes, R., Feb. 2015. Joint optimization of well placement and
  control for nonconventional well types. Journal of Petroleum Science and
  Engineering 126, 242--253.

\bibitem[{Humphries et~al.(2013)Humphries, Haynes, and
  James}]{humphries_simultaneous_2013}
Humphries, T.~D., Haynes, R.~D., James, L.~A., Sep. 2013. Simultaneous and
  sequential approaches to joint optimization of well placement and control.
  Computational Geosciences 18~(3-4), 433--448.

\bibitem[{Huyer and Neumaier(1999)}]{huyer_global_1999}
Huyer, W., Neumaier, A., 1999. Global optimization by multilevel coordinate
  search. Journal of Global Optimization 14~(4), 331--355.

\bibitem[{Isebor(2009)}]{isebor_constrained_2009}
Isebor, O.~J., 2009. Constrained production optimization with an emphasis on
  derivative-free methods. Ph.D. thesis, Stanford University.

\bibitem[{Isebor et~al.(2014{\natexlab{a}})Isebor, Ciaurri, and
  Durlofsky}]{isebor_generalized_2014}
Isebor, O.~J., Ciaurri, D.~E., Durlofsky, L.~J., Oct. 2014{\natexlab{a}}.
  Generalized field-development optimization with derivative-free procedures.
  {SPE} Journal 19~(05), 891--908.

\bibitem[{Isebor et~al.(2014{\natexlab{b}})Isebor, Durlofsky, and
  Ciaurri}]{isebor_derivative-free_2014}
Isebor, O.~J., Durlofsky, L.~J., Ciaurri, D.~E., Aug. 2014{\natexlab{b}}. A
  derivative-free methodology with local and global search for the constrained
  joint optimization of well locations and controls. Computational Geosciences
  18~(3-4), 463--482.

\bibitem[{Jansen et~al.(2014)Jansen, Fonseca, Kahrobaei, Siraj, Van~Essen, and
  Van~den Hof}]{jansen_egg_2014}
Jansen, J.~D., Fonseca, R.~M., Kahrobaei, S., Siraj, M.~M., Van~Essen, G.~M.,
  Van~den Hof, P. M.~J., Nov. 2014. The egg model – a geological ensemble for
  reservoir simulation. Geoscience Data Journal 1~(2), 192--195.

\bibitem[{Jones et~al.(1993)Jones, Perttunen, and
  Stuckman}]{jones_lipschitzian_1993}
Jones, D.~R., Perttunen, C.~D., Stuckman, B.~E., Oct. 1993. Lipschitzian
  optimization without the {Lipschitz} constant. Journal of Optimization Theory
  and Applications 79~(1), 157--181.

\bibitem[{Kampstra(2008)}]{kampstra_beanplot:_2008}
Kampstra, P., 2008. Beanplot: A boxplot alternative for visual comparison of
  distributions. Journal of Statistical Software 28~(1).

\bibitem[{Kennedy(2011)}]{kennedy_particle_2011}
Kennedy, J., 2011. Particle {Swarm} {Optimization}. In: Sammut, C., Webb, G.~I.
  (Eds.), Encyclopedia of {Machine} {Learning}. Springer Science + Business
  Media, pp. 760--766.

\bibitem[{Knudsen and Foss(2013)}]{knudsen_shut-based_2013}
Knudsen, B.~R., Foss, B., Nov. 2013. Shut-in based production optimization of
  shale-gas systems. Computers \& Chemical Engineering 58, 54--67.

\bibitem[{Kolda et~al.(2003)Kolda, Lewis, and
  Torczon}]{kolda_optimization_2003}
Kolda, T., Lewis, R., Torczon, V., Jan. 2003. Optimization by {Direct}
  {Search}: {New} {Perspectives} on {Some} {Classical} and {Modern} {Methods}.
  SIAM Review 45~(3), 385--482.

\bibitem[{Lambot et~al.(2002)Lambot, Javaux, Hupet, and
  Vanclooster}]{lambot_global_2002}
Lambot, S., Javaux, M., Hupet, F., Vanclooster, M., Nov. 2002. A global
  multilevel coordinate search procedure for estimating the unsaturated soil
  hydraulic properties. Water Resources Research 38~(11), 6--15.

\bibitem[{Li and Jafarpour(2012)}]{li_variable-control_2012}
Li, L., Jafarpour, B., 2012. A variable-control well placement optimization for
  improved reservoir development. Computational Geosciences 16~(4), 871--889.

\bibitem[{Li et~al.(2012)Li, Jafarpour, and
  Mohammad-Khaninezhad}]{li_simultaneous_2012}
Li, L., Jafarpour, B., Mohammad-Khaninezhad, M.~R., Nov. 2012. A simultaneous
  perturbation stochastic approximation algorithm for coupled well placement
  and control optimization under geologic uncertainty. Computational
  Geosciences 17~(1), 167--188.

\bibitem[{Li et~al.(2003)Li, Reynolds, and Oliver}]{li_history_2003}
Li, R., Reynolds, A.~C., Oliver, D.~S., 2003. History {Matching} of
  {Three}-{Phase} {Flow} {Production} {Data}. SPE Journal 8~(04), 328--340.

\bibitem[{Loshchilov(2013)}]{loshchilov_cma-es_2013}
Loshchilov, I., Jun. 2013. {CMA}-{ES} with restarts for solving {CEC} 2013
  benchmark problems. In: 2013 {IEEE} Congress on Evolutionary Computation.
  Institute of Electrical {\&} Electronics Engineers ({IEEE}).

\bibitem[{Merlini~Giuliani and
  Camponogara(2015)}]{merlini_giuliani_derivative-free_2015}
Merlini~Giuliani, C., Camponogara, E., Apr. 2015. Derivative-free methods
  applied to daily production optimization of gas-lifted oil fields. Computers
  \& Chemical Engineering 75, 60--64.

\bibitem[{Neumaier(2008)}]{software_mcs}
Neumaier, A., 2008. {MCS}: Global optimization by multilevel coordinate search.
  \url{https://www.mat.univie.ac.at/~neum/software/mcs/}.

\bibitem[{Oliveira and Reynolds(2014)}]{oliveira_adaptive_2014}
Oliveira, D.~F., Reynolds, A., Oct. 2014. An adaptive hierarchical multiscale
  algorithm for estimation of optimal well controls. {SPE} Journal 19~(05),
  909--930.

\bibitem[{Onwunalu(2010)}]{onwunalu_optimization_2010}
Onwunalu, J.~E., 2010. Optimization of field development using particle swarm
  optimization and new well pattern descriptions. Ph.D. thesis, Stanford
  University.

\bibitem[{Onwunalu and Durlofsky(2011)}]{onwunalu_new_2011}
Onwunalu, J.~E., Durlofsky, L., Sep. 2011. A new well-pattern-optimization
  procedure for large-scale field development. {SPE} Journal 16~(03), 594--607.

\bibitem[{Onwunalu and Durlofsky(2009)}]{onwunalu_development_2009}
Onwunalu, J.~E., Durlofsky, L.~J., 2009. Development and application of a new
  well pattern optimization algorithm for optimizing large scale field
  development. In: {SPE} Annual Technical Conference and Exhibition. Society of
  Petroleum Engineers ({SPE}).

\bibitem[{Perez and Behdinan(2007)}]{perez_particle_2007}
Perez, R.~E., Behdinan, K., Oct. 2007. Particle swarm approach for structural
  design optimization. Computers \& Structures 85~(19–20), 1579--1588.

\bibitem[{Po{\v{s}}{\'{\i}}k et~al.(2012)Po{\v{s}}{\'{\i}}k, Huyer, and
  P{\'{a}}l}]{posik_comparison_2012}
Po{\v{s}}{\'{\i}}k, P., Huyer, W., P{\'{a}}l, L., Dec. 2012. A comparison of
  global search algorithms for continuous black box optimization. Evolutionary
  Computation 20~(4), 509--541.

\bibitem[{Rios and Sahinidis(2013)}]{rios_derivative-free_2013}
Rios, L.~M., Sahinidis, N.~V., Jul. 2013. Derivative-free optimization: a
  review of algorithms and comparison of software implementations. Journal of
  Global Optimization 56~(3), 1247--1293.

\bibitem[{Sarma et~al.(2005)Sarma, Aziz, and Durlofsky}]{sarma_impl_2005}
Sarma, P., Aziz, K., Durlofsky, L.~J., 2005. Implementation of adjoint solution
  for optimal control of smart wells. In: {SPE} {Reservoir} {Simulation}
  {Symposium}. Society of Petroleum Engineers.

\bibitem[{Sarma and Chen(2008)}]{sarma_efficient_2008}
Sarma, P., Chen, W.~H., 2008. Efficient well placement optimization with
  gradient-based algorithms and adjoint models. In: Intelligent {Energy}
  {Conference} and {Exhibition}. Society of Petroleum Engineers.

\bibitem[{Sarma et~al.(2006)Sarma, Durlofsky, Aziz, and
  Chen}]{sarma_efficient_2006}
Sarma, P., Durlofsky, L.~J., Aziz, K., Chen, W.~H., Mar. 2006. Efficient
  real-time reservoir management using adjoint-based optimal control and model
  updating. Computational Geosciences 10~(1), 3--36.

\bibitem[{Shakhsi-Niaei et~al.(2014)Shakhsi-Niaei, Iranmanesh, and
  Torabi}]{shakhsi-niaei_optimal_2014}
Shakhsi-Niaei, M., Iranmanesh, S.~H., Torabi, S.~A., Jun. 2014. Optimal
  planning of oil and gas development projects considering long-term production
  and transmission. Computers \& Chemical Engineering 65, 67--80.

\bibitem[{Siraj et~al.(2015)Siraj, Van~den Hof, and Jansen}]{siraj_model_2015}
Siraj, M.~M., Van~den Hof, P.~M., Jansen, J.~D., 2015. Model and {Economic}
  {Uncertainties} in {Balancing} {Short}-{Term} and {Long}-{Term} {Objectives}
  in {Water}-{Flooding} {Optimization}. In: {SPE} {Reservoir} {Simulation}
  {Symposium}. Society of Petroleum Engineers.

\bibitem[{Tavallali et~al.(2013)Tavallali, Karimi, Teo, Baxendale, and
  Ayatollahi}]{tavallali_optimal_2013}
Tavallali, M.~S., Karimi, I.~A., Teo, K.~M., Baxendale, D., Ayatollahi, S.,
  Aug. 2013. Optimal producer well placement and production planning in an oil
  reservoir. Computers \& Chemical Engineering 55, 109--125.

\bibitem[{Torczon(1997)}]{torczon_convergence_1997}
Torczon, V., Feb. 1997. On the {Convergence} of {Pattern} {Search}
  {Algorithms}. SIAM Journal on Optimization 7~(1), 1--25.

\bibitem[{Vaz and Vicente(2007)}]{vaz_particle_2007}
Vaz, A. I.~F., Vicente, L.~N., Oct. 2007. A particle swarm pattern search
  method for bound constrained global optimization. Journal of Global
  Optimization 39~(2), 197--219.

\bibitem[{Vlemmix et~al.(2009)Vlemmix, Joosten, Brouwer, and
  Jansen}]{vlemmix_adjoint-based_2009}
Vlemmix, S., Joosten, G., Brouwer, R., Jansen, J.-D., 2009. Adjoint-based
  {Well} {Trajectory} {Optimization} in a {Thin} {Oil} {Rim} ({SPE}-121891).
  In: 71st {EAGE} {Conference} \& {Exhibition}.

\bibitem[{Volkov and Voskov(2014)}]{volkov_effect_2014}
Volkov, O., Voskov, D., 2014. Effect of time stepping strategy on adjoint-based
  production optimization. In: {ECMOR} {XIV} - 14th European conference on the
  mathematics of oil recovery. {EAGE} Publications.

\bibitem[{Wang et~al.(2007)Wang, Li, and Reynolds}]{wang_optimal_2007}
Wang, C., Li, G., Reynolds, A.~C., 2007. Optimal well placement for production
  optimization. In: Eastern {Regional} {Meeting}. Society of Petroleum
  Engineers.

\bibitem[{Wang et~al.(2009)Wang, Li, and Reynolds}]{wang_production_2009}
Wang, C., Li, G., Reynolds, A.~C., Sep. 2009. Production optimization in
  closed-loop reservoir management. {SPE} Journal 14~(03), 506--523.

\bibitem[{Wu et~al.(1999)Wu, Reynolds, and Oliver}]{wu_conditioning_1999}
Wu, Z., Reynolds, A.~C., Oliver, D.~S., 1999. Conditioning {Geostatistical}
  {Models} to {Two}-{Phase} {Production} {Data}. SPE Journal 4~(02), 142--155.

\bibitem[{Yeten et~al.(2002)Yeten, Durlofsky, and
  Aziz}]{yeten_optimization_2002}
Yeten, B., Durlofsky, L.~J., Aziz, K., 2002. Optimization of nonconventional
  well type, location and trajectory. In: {SPE} annual technical conference and
  exhibition. Society of Petroleum Engineers.

\bibitem[{Yin and Cagan(2000)}]{yin_extended_2000}
Yin, S., Cagan, J., 2000. An extended pattern search algorithm for
  three-dimensional component layout. Journal of Mechanical Design 122~(1),
  102--108.

\bibitem[{Zakirov et~al.(1996)Zakirov, Aanonsen, Zakirov, and
  Palatnik}]{zakirov_optimizing_1996}
Zakirov, I., Aanonsen, S.~I., Zakirov, E.~S., Palatnik, B.~M., 1996. Optimizing
  reservoir performance by automatic allocation of well rates. In: 5th
  {European} {Conference} on the {Mathematics} of {Oil} {Recovery}.

\bibitem[{Zandvliet et~al.(2008)Zandvliet, Handels, van Essen, Brouwer, and
  Jansen}]{zandvliet_adjoint-based_2008}
Zandvliet, M., Handels, M., van Essen, G., Brouwer, R., Jansen, J.-D., 2008.
  Adjoint-based well-placement optimization under production constraints. SPE
  Journal 13~(04), 392--399.

\bibitem[{Zhao et~al.(2013)Zhao, Chen, Do, Oliveira, Li, and
  Reynolds}]{zhao_maximization_2013}
Zhao, H., Chen, C., Do, S., Oliveira, D., Li, G., Reynolds, A., 2013.
  Maximization of a {Dynamic} {Quadratic} {Interpolation} {Model} for
  {Production} {Optimization}. SPE Journal 18~(06), 1--012.

\bibitem[{Zhou et~al.(2013)Zhou, Hou, Zhang, Du, Kang, and
  Jiang}]{zhou_optimal_2013}
Zhou, K., Hou, J., Zhang, X., Du, Q., Kang, X., Jiang, S., Aug. 2013. Optimal
  control of polymer flooding based on simultaneous perturbation stochastic
  approximation method guided by finite difference gradient. Computers \&
  Chemical Engineering 55, 40--49.

\end{thebibliography}


\end{document}